\documentclass[twoside,a4paper]{article} % double-sided
\usepackage{amssymb}
\usepackage{a4wide} % gives a4 paper format (compared to US paper)
\usepackage{latexsym} % symbol and maths package
\usepackage{amsmath} % symbol and maths package
\usepackage{theorem} % for theorem environment in maths
\usepackage{mathrsfs} %nice calligraphic characters   %special blackboard bold
\usepackage{xspace}  % useful for macro definitions
\usepackage{fancyhdr} % for header and footer on pages
\usepackage{float} % for graphics/including pictures
\usepackage{listings}% very useful aid for including computer code in the document
\usepackage{enumerate}% useful add on for numbered lists
\usepackage{morefloats}
\usepackage{tikz} 
\usepackage{accents}
\usepackage[compact]{titlesec}
\usepackage{graphicx} 
\usepackage{caption}
\usepackage{rotating}
\usepackage{subcaption}
\definecolor{SussexFlint}{rgb}{.00,.19,.21}
\definecolor{SussexGrey}{rgb}{.51,.58,.49}
\definecolor{SussexOrange}{rgb}{.94,.29,.00}
\definecolor{SussexYellow}{rgb}{1.00,.73,.00}
\definecolor{SussexRed}{rgb}{.94,.01,.49}
\definecolor{SussexPurple}{rgb}{.48,.06,.44}
\definecolor{SussexGreen}{rgb}{.00,.58,.46}
\definecolor{SussexBlue}{rgb}{.00,.58,.65}
\colorlet{a}{SussexOrange}
\colorlet{b}{SussexRed}
\colorlet{c}{SussexYellow}
\colorlet{d}{SussexPurple}
\colorlet{e}{SussexGreen}
\colorlet{f}{SussexBlue}
\colorlet{g}{SussexGrey}
\colorlet{h}{white}
\colorlet{i}{black}
\colorlet{j}{SussexFlint}
\numberwithin{equation}{section}
\newtheorem{theorem}{Theorem}[section]
\newtheorem{lemma}[theorem]{Lemma}
\newtheorem{corollary}[theorem]{Corollary}
\newtheorem{definition}[theorem]{Definition}%
\newtheorem{remark}[theorem]{Remark}
\newcommand{\qed}{\hfill$\square$}

\newcommand{\avg}[1]{\ensuremath{\langle\!\langle#1\rangle\!\rangle} }
\newcommand{\bavg}[1]{\ensuremath{\left\langle\!\!\!\left\langle#1\right\rangle\!\!\!\right\rangle} }
\newcommand{\jump}[1]{\ensuremath{[\![#1]\!]} }
\newcommand{\bjump}[1]{\ensuremath{\left[\!\!\left[#1\right]\!\!\right]} }

\newcommand{\Th}{\ensuremath{\mathscr{T}_h}}
\newcommand{\Ta}{\ensuremath{\mathbf{T}}}
\newcommand{\mc}{\ensuremath{\mathcal{H}}}

\newcommand{\nablaTa}{\ensuremath{\nabla_\Ta}}
\newcommand{\Eb}{\ensuremath{\mathscr{E}^b_h}}
\newcommand{\Ei}{\ensuremath{\mathscr{E}^i_h}}
\newcommand{\X}{\ensuremath{\Omega}}

\newcommand{\Eib}{\ensuremath{\mathscr{E}^{i,b}_h}}

\newcommand{\dg}{\ensuremath{V_{h,p}}}

\newcommand{\diam}{\ensuremath{\operatorname{diam}}}
\newtheorem{assumption}[theorem]{Assumption}

\newcommand{\dgcomp}{\ensuremath{\dg^{\operatorname{comp}}}}

\newcommand{\ms}[1]{\ensuremath{\min\{1,#1\}} }
\title{Finite element theory on curved domains with applications to DGFEMs}
\author{Ellya L. Kawecki\footnote{ELK acknowledges support of the Engineering and Physical Sciences Research Council
[EP/L015811/1].}}
\begin{document}
\maketitle
\begin{abstract}
In this paper we provide key estimates used in the stability and error analysis of discontinuous Galerkin finite element methods (DGFEMs) on domains with curved boundaries. In particular, we review trace estimates, inverse estimates, discrete Poincar\'e--Friedrichs' inequalities, and optimal interpolation estimates in noninteger Hilbert-Sobolev norms, that are well known in the case of polytopal domains. We also prove curvature bounds for curved simplices, which does not seem to be present in the existing literature, even in the polytopal setting, since polytopal domains have piecewise zero curvature. We demonstrate the value of these estimates, by analysing the IPDG method for the Poisson problem, introduced by Douglas and Dupont [\emph{Computing Methods in Applied Sciences, Lecture Notes in Physics, vol 58. Springer, Berlin, Heidelberg}, pages 207--216. Springer, 1976], and by analysing a variant of the $hp$-DGFEM for the biharmonic problem introduced by Mozolevski and S\"{u}li  [\emph{Computer Methods in Applied Mechanics and Engineering}, 196(13-16):1851--1863, 2007]. In both cases we prove stability estimates and optimal a priori error estimates. Numerical results are provided, validating the proven error estimates.
\end{abstract}
\begin{section}{Introduction}
When modelling second- and fourth-order (as well as higher order) elliptic partial differential equations (PDEs), one may be required to consider a domain that cannot be expressed as a finite union of polytopes, for example, the unit ball, $B_1(0):=\{x\in\mathbb R^d:|x|<1\}\subset\mathbb R^d$. This necessity could be driven by the domain considered in the underlying application, where the domain is for example Lipschitz continuous, and piecewise $C^{1,\alpha}$, $\alpha\in(0,1)$, but \emph{not} piecewise smooth, or for example the domain is $C^1$, and thus not polytopal. Such domains arise naturally in the theory of PDEs, for example, a natural assumption for the Monge--Amp\`{e}re equation~\cite{MR0180763,MR1454261,MR718679,trudinger2008monge,caffarelli1984dirichlet,kawecki:lakkis:pryer} is that the domain is uniformly convex~\cite{MR1454261,MR718679,caffarelli1984dirichlet}, and oblique boundary-value problems~\cite{MR3059278,MR1836250,MR2260015,MR1624426} in nondivergence form, with bounded and measurable coefficients, require a $C^2$ boundary assumption~\cite{MR2260015}, both of which rule out the possibility of a polyhedral domain. When it comes to finite element methods (FEMs), it is useful if the domain is polytopal, then since one can discretise the domain, $\X$, \emph{exactly} by polytopes, i.e., there exists a family of shape-regular meshes $(\Th)_{h>0}$ on $\overline{\X}$ for which $\X=\cup_{K\in\Th}\overline{K}$ (the sets $K$ are often $d$-simplices or parallelipeds). 

If the boundary of $\X$ is curved, an \emph{exact} mesh consisting of a finite set of polyhedrons cannot be obtained; one must instead use curved elements. In~\cite{MR1014883}, the author introduces the concept of exact curved domain approximation by curved $d$-simplices, following~\cite{MR2940387,MR842644}, providing an optimal (with respect to the parameter $h$) finite element interpolant (interpolating with and without boundary conditions), with estimates in $W^{m,p}$-norms, $m\in\mathbb N_0$, $p\in[1,\infty]$.

We will see, however, that in order to design and analyse discontinuous Galerkin finite element methods (DGFEMs) for second- and fourth-order elliptic PDEs on domains with curved boundaries, one requires further estimates, in particular: inverse estimates; discrete Poincar\'e--Friedrichs' inequalities; simplicial curvature bounds; and optimal interpolation estimates in noninteger Sobolev norms.
% The final member of the list allows one to potentially utilise the extra piecewise regularity of the weak solution to the PDE that comes the piecewise $C^{k,\alpha}$ smoothness of the domain (assuming that the source term is also sufficiently regular).
Furthermore, since curved domain approximations require the composition of piecewise polynomials with functions that are not piecewise polynomials (the details of this will be made clearer in Section~\ref{sec:3}), applications of the chain rule show that in general, the piecewise derivative no longer maps from the finite element space into itself (as is often seen in discontinuous Galerkin (DG) finite element spaces), complicating the derivation and structure of inverse estimates. For penalty FEMs for fourth order problems, we will see that this leads to the necessity of discrete Poincar\'e--Friedrichs' inequalities. 

One is often motivated to use DGFEMs, other nonconforming FEMs, and mixed FEMs over conforming FEMs, due to the structural and computational challenges that conforming FEMs impose. For conforming FEMs, it is required that the approximation space is a subset of the space of weak solutions to the PDE, examples of this being the spaces $H^1_0(\X)$ and $H^2_0(\X)$ for second- and fourth-order elliptic problems, such as the Poisson problem and biharmonic clamped plate problem, which we shall consider as our model second- and fourth-order problems. In the $H^1_0(\X)$ case, this can be achieved by considering piecewise polynomials that are \emph{globally} continuous, however, for $H^2_0(\X)$, one must also enforce continuity of the gradient across neighbouring elements. An example of this being the Argyris finite element~\cite{MR1930132}, which can be rather expensive to implement, requiring polynomials of degree five on two dimensional simplicial polynomials. In contrast, nonconforming methods weakly enforce this regularity by penalising jumps of the discrete functions, and their derivatives across the edges of neighbouring elements, and as a result, the methods that we consider only require a polynomial degree greater than or equal to the number of derivatives in the weak formulation of the PDE.

The remainder of this paper is organised as follows, in Section~\ref{sec:2} we shall discuss the existence and uniqueness of weak solutions to the Poisson and biharmonic equations, and discuss conforming finite element methods (FEMs), and DGFEMs (the latter of which falls into the category of nonconforming FEMs) on polytopal domains, with the goal of highlighting important features, such as the stability and consistency of such schemes. In Section~\ref{sec:3} we review the key tools from finite element analysis that are well known in the polytopal case, in the context of curved simplicial finite elements. In Section~\ref{sec:4} we will provide the numerical methods for the Poisson and biharmonic problems, and prove that they are stable, and in Section~\ref{sec:5} we prove that the numerical solutions satisfy optimal a priori error estimates in $H^k$-type norms. Finally, in Section~\ref{sec:6} we provide numerical experiments that validate the error estimates of Section~\ref{sec:5}.
\section{Weak formulations, conforming and nonconforming methods}\label{sec:2}
For $k\in\mathbb N$, we denote the standard Hilbert-Sobolev space~\cite{MR2597943} $$H^k_0(K):=\{v\in L^2(K):D^\alpha v\in L^2(K)\,\,\forall\alpha:|\alpha|\le k,D^\beta v|_{\partial K} = 0\,\,\forall\beta:|\beta|\le k-1\},$$
where the restriction to $\partial K$ is considered in the sense of traces. 

Let $\X\subset\mathbb R^d$ be Lipschitz continuous, and consider the following second- and fourth-order elliptic boundary-value problems, for $k=1,2$, find $u_k:\X\to\mathbb R$ such that:
\begin{equation}\label{2and4}
\left\{
\begin{aligned}
(-\Delta)^k u_k & = f\quad\mbox{in }\X,\\
\frac{\partial^ju_k}{\partial n_{\partial\X}^j} & = 0\quad\mbox{on }\partial\X,\,\,0\le j\le k-1,
\end{aligned}
\right.
\end{equation}
where $f\in L^2(\X)$. When $k=1$,~(\ref{2and4}) is the well known Poisson problem, and for $k=2$,~(\ref{2and4}) is the biharmonic clamped plate problem. In particular, one can show that in each case, there exists a unique weak solution $u_k\in H^k_0(\X)$. That is, $u_k$ satisfies
\begin{equation}\label{weakforms}
a_k(u_k,v) = \int_\X fv\quad\forall v\in H^k_0(\X),
\end{equation}
where the bilinear forms $a_k:H^k_0(\X)\times H^k_0(\X)\to\mathbb R$, are given by
\begin{align}
a_1(u,v) &:= \int_\X\nabla u\cdot\nabla v,\quad\forall u,v\in H^1_0(\X),\label{a1def}\\
a_2(u,v) &:= \int_\X\Delta u\,\Delta v,\quad\forall u,v\in H^2_0(\X).\label{a2def}
\end{align}
Note that the existence of such functions follows from applying the Lax--Milgram Theorem~\cite{MR2597943}; in particular one must show that the bilinear forms are coercive in the $H^k$-norm. In the case that $k=1$, this follows from the Poincar\'e inequality~\cite{MR2597943}, and for $k=2$, the following identity (see~(\ref{calderon-1}))
$$\int_\X\Delta u\,\Delta v = \int_\X D^2u\!:\!D^2v,\quad u,v\in H^2_0(\X),$$
implies that $|u|_{H^2(\X)}=\|\Delta u\|_{L^2(\X)}$ if $u\in H^2_0(\X)$, which, coupled with the Poincar\'e inequality, also proves the coercivity of $a_2$.

The derivation of the weak formulations~(\ref{weakforms}) follows from the following integration by parts identities, valid for functions $u,v\in C^\infty(\overline{K})$, where $K\subset\mathbb R^d$ has a Lipschitz boundary, and extendable to $u,v$ in suitable Sobolev spaces by density:
\begin{equation}\label{ibp1}
\int_K(-\Delta u)v = \int_K\nabla u\cdot\nabla v-\int_{\partial K}\frac{\partial u}{\partial n_{\partial K}}v,
\end{equation}
and
\begin{equation}\label{ibp2}
\begin{aligned}
\int_K(\Delta^2u)v &= -\int_K\nabla(\Delta u)\cdot\nabla v+\int_{\partial K}\frac{\partial(\Delta u)}{\partial n_{\partial K}}v,\\
& = \int_K\Delta u\,\Delta v+\int_{\partial K}\frac{\partial(\Delta u)}{\partial n_{\partial K}}v-\Delta u\frac{\partial v}{\partial n_{\partial K}},
\end{aligned}
\end{equation}
where $n_{\partial K}$ is the unit outward normal to ${\partial K}$.
Taking $K=\X$, the choice of $u,v\in H^k_0(\X)$ justifies the lack of the appearance of boundary integrals in~(\ref{a1def})--(\ref{a2def}) (however, for this we utilise the density of $C^\infty_c(\X)$ in $H^k_0(\X)$). 

For a conforming finite element method, one assumes that the finite dimensional space $V_{k,h}\subset H^k_0(\X)$, and so one may obtain a conforming finite element method by directly substituting the finite element functions into the bilinear forms. That is, one seeks $u_{k,h}\in V_{k,h}$ such that
\begin{equation}\label{fem}
a_k(u_{k,h},v_h) = \int_\X fv_h\quad\forall v_h\in V_{k,h}.
\end{equation}
Indeed, since $V_{k,h}\subset H^k_0(\X)$, the properties of the bilinear forms are still valid on $V_{k,h}\times V_{k,h}$, and so the existence and uniqueness of a numerical solution follows in a similar manner to the existence and uniqueness of a weak solution. In particular, the bilinear form $a_k$ is coercive on $V_{k,h}\times V_{k,h}$ in the $H^k(\X)$ norm, and so we obtain the stability estimate
\begin{equation}\label{stab}
a_k(v_h,v_h)\ge C_k\|v_h\|_{H^k(\X)}^2,\quad\forall v_h\in V_{k,h},
\end{equation}
where $C_k$ is a positive constant independent of the approximation parameter $h$. Since the problem~(\ref{fem}) is equivalent to solving a linear system of equations, the stability estimate implies uniqueness, which in turn implies invertibility of the matrix describing the corresponding linear system, which also yields existence. 

Furthermore, we see that the true solutions, $u_k\in H^k_0(\X)$ satisfy
\begin{equation}\label{cons}
a_k(u_k,v_h) = \int_\X fv_h\quad\forall v_h\in V_{k,h},
\end{equation}
and so
\begin{equation}\label{gorthog}
a_k(u_k-u_{k,h},v_h) = 0\quad\forall v_h\in V_{k,h}.
\end{equation}
A finite element method that satisfies~(\ref{cons}) is called \emph{consistent}, and~(\ref{gorthog}) is referred to as Galerkin orthogonality, which, when combined with the stability estimate~(\ref{stab}), yields Cea's Lemma:
$$\|u_k-u_{k,h}\|_{H^k(\X)}\le C_k\inf_{v_h\in V_{k,h}}\|u_k-v_h\|_{H^k(\X)}.$$
One obtains optimal error estimates, by noting that the infimum over $V_{k,h}$ is bounded above by any choice of $z_h\in V_{k,h}$. In particular, assuming that $u_k\in H^s(\X)\cap H^k_0(\X)$, $s\ge k$, one may choose $z_h$ to coincide with a suitable interpolant, yielding
\begin{equation}\label{interpolationbound}
\|u_k-u_{k,h}\|_{H^k(\X)}\le C_k\inf_{v_h\in V_{k,h}}\|u_k-v_h\|_{H^k(\X)}\le C_*C_k\|u_k\|_{H^s(\X)}h^{\min\{p+1,s\}-k}.
\end{equation}

Unlike conforming finite element methods, where the approximating space $V_{k,h}$ is a subset of $H^k_0(\X)$, nonconforming finite element methods involve approximating spaces for which this is not true; in the case of DGFEMs one only has $V_{k,h}\subset L^2(\X)$, and for the $C^0$-interior penalty method proposed in~\cite{MR2142191}, one has $V_{2,h}\subset H^1_0(\X)$, which is nonconforming in the sense that $V_{2,h}$ is not contained in $H^2_0(\X)$. 

For DGFEMs, one also has analogues of stability, consistency, and optimal error estimates. However, since the finite element functions do not have sufficient \emph{global} regularity, one cannot directly substitute $u_h,v_h\in V_{k,h}:=\{v\in L^2(\X):v|_K\in\mathbb P^p(K)\,\forall K\in\Th\}$ (implicitly, we assume $p\ge k$, and that $(\Th)_{h>0}$ is a family of regular simplicial meshes on $\overline{\X}$) into the bilinear forms $a_k$, $k=1,2$. 

Such functions do, however, satisfy a property of piecewise regularity; since, $V_{k,h}\subset H^{2k}(\X;\Th):=\{v\in L^2(\X):v|_K\in H^{2k}(K)\,\forall K\in\Th\}$ (in particular, piecewise polynomials are piecewise smooth) and so, assuming $u_k\in H^k_0(\X)\cap H^{2k}(\X)$ are the weak solutions to the PDE, we can sum the integration by parts identities~(\ref{ibp1}) and~(\ref{ibp2}) over all $K\in\Th$, obtaining (see Definition~\ref{faceandvertexsets:def}, as well as~(\ref{average}) and~(\ref{jump}) for the relevant notational conventions in present in the identities that follow):
\begin{equation}\label{conslap}
 \sum_{K\in\Th}\int_K\nabla u_1\cdot\nabla v_h-\int_{\partial K}\frac{\partial u_1}{\partial n_{\partial K}}v_h = \sum_{K\in\Th}\int_K(-\Delta u_1)v_h = \sum_{K\in\Th}\int_Kfv_h=:\ell(v_h)\quad\forall v_h\in V_{1,h},
\end{equation}
and
\begin{equation}\label{consbiharm}
\sum_{K\in\Th}\int_K\Delta u_2\,\Delta v_h+\int_{\partial K}\frac{\partial(\Delta u_2)}{\partial n_{\partial K}}v_h-\Delta u_2\frac{\partial v_h}{\partial n_{\partial K}} = \sum_{K\in\Th}\int_K((-\Delta )^2u_2)v_h = \sum_{K\in\Th}\int_Kfv_h\quad\forall v_h\in V_{2,h}.
\end{equation}
Since $u_k\in H^{2k}(\X)\cap H^k_0(\X)$, it follows that
\begin{align}
\jump{D^\alpha u_k} & = 0\quad\forall F\in\Eb,\quad|\alpha|\le k-1,\quad k = 1,2,\label{jumps0u1}\\
\jump{D^\alpha u_k} & = 0\quad\forall F\in\Ei,\quad|\alpha|\le 2k-1,\quad k = 1,2.\label{jumps0u2}
\end{align}
Thus, we obtain
\begin{equation}\label{jumpidconslap}
\begin{aligned}
-\sum_{K\in\Th}\int_{\partial K}\frac{\partial u_1}{\partial n_{\partial K}}v_h & = -\sum_{F\in\Ei}\int_F\bjump{\frac{\partial u_1}{\partial n_F}}\avg{v_h}-\sum_{F\in\Eib}\int_F\bavg{\frac{\partial u_1}{\partial n_F}}\jump{v_h}\\
& = -\sum_{F\in\Eib}\int_F\bavg{\frac{\partial u_1}{\partial n_F}}\jump{v_h},\\
\end{aligned}
\end{equation}
and
\begin{equation}\label{jumpidconsbiharm}
\begin{aligned}
\sum_{K\in\Th}\int_{\partial K}\frac{\partial(\Delta u_2)}{\partial n_{\partial K}}v_h-\Delta u_2\frac{\partial v_h}{\partial n_{\partial K}} & = \sum_{F\in\Ei}\int_F\bjump{\frac{\partial(\Delta u_2)}{\partial n_{F}}}\avg{v_h}+\sum_{F\in\Eib}\int_F\bavg{\frac{\partial(\Delta u_2)}{\partial n_{F}}}\jump{v_h}\\
&~~~~~~-\sum_{F\in\Ei}\int_F\jump{\Delta u_2}\bavg{\frac{\partial v_h}{\partial n_{F}}}-\sum_{F\in\Eib}\int_F\avg{\Delta u_2}\bjump{\frac{\partial v_h}{\partial n_{F}}}\\
& = \sum_{F\in\Eib}\int_F\bavg{\frac{\partial(\Delta u_2)}{\partial n_{F}}}\jump{v_h}-\avg{\Delta u_2}\bjump{\frac{\partial v_h}{\partial n_{F}}},
\end{aligned}
\end{equation}
where $n_F$ denotes a \emph{fixed} choice of unit normal to $F$.
Let us define
$$B_1(u_h,v_h):=-\sum_{F\in\Eib}\int_F\bavg{\frac{\partial u_h}{\partial n_F}}\jump{v_h},$$
and
$$B_2(u_h,v_h):=\sum_{F\in\Eib}\int_F\bavg{\frac{\partial(\Delta u_h)}{\partial n_{F}}}\jump{v_h}-\avg{\Delta u_h}\bjump{\frac{\partial v_h}{\partial n_{F}}}.$$
Then, defining $\tilde{A}_k:V_{k,h}\times V_{k,h}\to\mathbb R$, $k=1,2$, by
\begin{equation}
\tilde{A}_k(u_h,v_h):=a_k(u_h,v_h)+B_k(u_h,v_h)\quad\forall u_h,v_h\in V_{k,h},\,\,k=1,2,
\end{equation}
we arrive at the following DGFEMs for the approximation of the solutions $u_k$, $k=1,2$, to~(\ref{2and4}): find $u_{k,h}\in V_{k,h}$ such that
\begin{equation}\label{consfork}
\tilde{A}_k(u_{k,h},v_h) = \ell(v_h)\quad\forall v_h\in V_{k,h}.
\end{equation}
Identities~(\ref{conslap})--(\ref{jumpidconsbiharm}), imply that the FEMs given by~(\ref{consfork}) for $k=1,2$, are \emph{consistent}, that is if $u_k\in H^{2k}(\X)\cap H^k_0(\X)$ solve~(\ref{2and4}) for $k=1,2$, then,
\begin{equation}\label{consforktrue}
\tilde{A}_k(u_k,v_h) = \ell(v_h)\quad\forall v_h\in V_{k,h}.
\end{equation}
Furthermore, we see that
\begin{align}
a_1(v_h,v_h) & = \sum_{K\in\Th}|v_h|_{H^1(K)}^2=:|v_h|_{H^1(\X;\Th)}^2\quad\forall v_h\in V_{1,h},\\
a_2(v_h,v_h) & = \sum_{K\in\Th}\|\Delta v_h\|_{L^2(K)}^2=:|v_h|_{H_\Delta(\X;\Th)}^2\quad\forall v_h\in V_{2,h},\label{A2stab}
\end{align}
but, the remaining terms present in $\tilde{A}_1$ and $\tilde{A}_2$ are not bounded quite as simply.
If $F$ is a face of $K\in\Th$, trace estimates yield for $w_k\in H^{2k}(K)$
\begin{align*}
\left\|\frac{\partial w_1}{\partial n_F}\right\|_{L^2(F)}^2&\le C(\tilde{h}_F^{-1}|w_1|_{H^1(K)}^2+\tilde{h}_F|w_1|_{H^2(K)}^2),\\
\left\|\frac{\partial (\Delta w_2)}{\partial n_F}\right\|_{L^2(F)}^2&\le C(\tilde{h}_F^{-1}|\Delta w_2|_{H^1(K)}^2+\tilde{h}_F|\Delta w_2|_{H^2(K)}^2),\\
\|\Delta w_2\|_{L^2(F)}^2&\le C(\tilde{h}_F^{-1}\|\Delta w_2\|_{L^2(K)}^2+\tilde{h}_F|\Delta w_2|_{H^1(K)}^2),
\end{align*}
where $C$ depends upon the shape-regularity constant of $\Th$. Then, applying inverse estimates~\cite{MR2373954} of the form
\begin{equation}\label{22:04-1}
|w|_{H^2(K)}\le Ch_K^{-k}|w|_{H^{2-k}(K)},
\end{equation}
for $w\in\mathbb P^p(K)$,
gives us
\begin{align}
\left\|\frac{\partial w_1}{\partial n_F}\right\|_{L^2(F)}^2&\le C\tilde{h}_F^{-1}|w_1|_{H^1(K)}^2,\label{22:041}\\
\left\|\frac{\partial (\Delta w_2)}{\partial n_F}\right\|_{L^2(F)}^2&\le C\tilde{h}_F^{-3}|\Delta w_2|_{L^2(K)}^2\label{22:04}\\
\|\Delta w_2\|_{L^2(F)}^2&\le C\tilde{h}_F^{-1}\|\Delta w_2\|_{L^2(K)}^2\label{22:04:222},
\end{align}
so long as $w_k\in\mathbb P^p(K)\subset H^{2k}(K)$, $k=1,2$. Then, utilising~(\ref{22:041})--(\ref{22:04:222}), and the Cauchy--Schwarz inequality with a parameter, yields the following for any $\delta_1>0$, and any $v_h\in V_{1,h}$
\begin{equation}\label{est1:10:52}
\begin{aligned}
\tilde{A}_1(v_h,v_h)&\ge|v_h|^2_{H^1(\X;\Th)}-\frac{1}{2}\sum_{F\in\Eib}\left[\delta_1\tilde{h}_F\left\|\bavg{\frac{\partial v_h}{\partial n_F}}\right\|_{L^2(F)}^2+(\delta_1\tilde{h}_F)^{-1}\|\jump{v_h}\|_{L^2(F)}^2\right]\\
&\ge|v_h|^2_{H^1(\X;\Th)}-\frac{\delta_1C}{2}\sum_{F\in\Eib}\sum_{K\in\Th:F\subset\partial K}|v_h|_{H^1(K)}^2-\sum_{F\in\Eib}(\delta_1\tilde{h}_F)^{-1}\|\jump{v_h}\|_{L^2(F)}^2\\
&\ge\left(1-\frac{\delta_1CC(d)}{2}\right)|v_h|^2_{H^1(\X;\Th)}-\frac{1}{2\delta_1}\sum_{F\in\Eib}\tilde{h}_F^{-1}\|\jump{v_h}\|_{L^2(F)}^2,\\
\end{aligned}
\end{equation}
where the final inequality holds due to the fact that the number of elements that share a given face is bounded in terms of the dimension, $d$. Similarly, for any $\delta_2>0$, and any $v_h\in V_{2,h}$, we see that
\begin{equation}\label{est2:10:52}
\begin{aligned}
&\tilde{A}_2(v_h,v_h)\ge|v_h|_{H_\Delta(\X;\Th)}^2-\frac{1}{2}\sum_{F\in\Eib}\left[\delta_2\tilde{h}_F^3\left\|\bavg{\frac{\partial (\Delta v_h)}{\partial n_F}}\right\|_{L^2(F)}^2+\delta_2\tilde{h}_F\|\avg{\Delta v_h}\|_{L^2(F)}^2\right.\\
&~~~~~~~~~~~~~~~~~~~~~~~~~~~~~~~~~~~~~+\left.(\delta_2\tilde{h}_F^3)^{-1}\|\jump{v_h}\|_{L^2(F)}^2+(\delta_2\tilde{h}_F)^{-1}\left\|\bjump{\frac{\partial v_h}{\partial n_F}}\right\|_{L^2(F)}^2\right]\\
&~~~~~~~~~~~~\ge(1-\delta_2CC(d))|v_h|_{H_\Delta(\X;\Th)}^2-\frac{1}{2\delta_2}\sum_{F\in\Eib}\left[\tilde{h}_F^{-3}\|\jump{v_h}\|_{L^2(F)}^2+\tilde{h}_F^{-1}\left\|\bjump{\frac{\partial v_h}{\partial n_F}}\right\|_{L^2(F)}^2\right].
\end{aligned}
\end{equation}
%Thus, it seems implausible to bound the counterparts $B_k$ in the following manner with $C_k<1$ and independent of $h$, $k=1,2$,
%$$|B_1(v_h,v_h)|\le C_1|v_h|_{H^1(\Th)}^2,\quad|B_2(v_h,v_h)|\le C_2|v_h|_{H_\Delta(\Th)}^2,\quad\forall v_h\in V_h,$$
% which would yield the stability estimates
%$$\tilde{A}_1(v_h,v_h) \ge \tilde{C}_1|v_h|_{H^1(\Th)}^2,\,\,\tilde{A}_2(v_h,v_h) \ge \tilde{C}_2|v_h|_{H_\Delta(\Th)}^2,\,\,\forall v_h\in V_h,$$
%where $\tilde{C}_k := 1-C_k>0$. Furthermore, neither $|\cdot|_{H^1(\Th)}^2$, nor $|\cdot|_{H_\Delta(\Th)}^2$ constitute norms on $V_h$, as they both admit nontrivial Kernels. 
The above estimates lead one to supplement the bilinear forms $\tilde{A}_k$, $k=1,2$, with additional bilinear forms $S_k,J_k:V_h\times V_h\to\mathbb R$, where the bilinear forms $J_k$ penalise interface jumps of the inputs and their piecewise weak derivatives up to order $2k-1$ across interior faces, and up to order $k-1$ on boundary faces, and the bilinear forms $S_k$ preserve the symmetry of the scheme. Clearly, the choice of $J_k$ and $S_k$ lead to different finite element methods; in~\cite{arnold2002unified} the authors present and analyse nine DG methods from~\cite{baumann1999discontinuous,bassi2000gmres,douglas1976interior,cockburn1998local,brezzi1999discontinuous,riviere1999improved,bassi1997high,babuvska1973nonconforming,brezzi2000discontinuous} for the Poisson problem ($k=1$), and in~\cite{suli2007hp} a $hp$-finite element method is introduced for the Biharmonic problem ($k=2$) with symmetric and nonsymmetric penalties. For other examples of nonconforming methods for second- and fourth-order elliptic problems see~\cite{MR2142191,MR3077903,Kawecki:77:article:On-curved,Kawecki:77:article:On-curved-oblique-long,brenner2017posteriori,brenner2012quadratic,brenner2011c,feng2018interior}. 

Thus, we may take
\begin{align*}
S_1(u_h,v_h)&:=-\sum_{F\in\Eib}\int_F\jump{u_h}\bavg{\frac{\partial v_h}{\partial n_F}},\\
J_1(u_h,v_h)&:=\sum_{F\in\Eib}\frac{\eta_F^1}{\tilde{h}_F}\int_F\jump{u_h}\jump{v_h},\\
S_2(u_h,v_h)&:=\sum_{F\in\Eib}\int_F\jump{u_h}\bavg{\frac{\partial(\Delta v_h)}{\partial n_F}}-\bjump{\frac{\partial u_h}{\partial n_F}}\avg{\Delta v_h}\\
J_2(u_h,v_h)&:=\sum_{F\in\Eib}\int_F\frac{\eta^2_F}{\tilde{h}_F^{3}}\jump{u_h}\jump{v_h}+\frac{\eta^3_F}{\tilde{h}_F}\bjump{\frac{\partial u_h}{\partial n_F}}\bjump{\frac{\partial v_h}{\partial n_F}},
\end{align*}
where $\eta_F^j$, $j=1,2,3,$ are positive parameter choices independent of $\tilde{h}_F$, that are chosen sufficiently large, in order to compensate for the jumps across $F\in\Eib$ present in estimates~(\ref{est1:10:52}) and~(\ref{est2:10:52}), as well as the jump estimates resulting from the terms included for symmetry that are present in $S_1$ and $S_2$ (these terms are bounded in exactly the same manner as in the derivation of estimates~(\ref{est1:10:52}) and~(\ref{est2:10:52})).
  
By~(\ref{jumps0u1})--(\ref{jumps0u2}), we see that $J_k(u_k,v_h) = 0$ for all $v_h\in V_h$, and so the bilinear forms
\begin{equation}\label{Akdef}
A_k(u_h,v_h):=\tilde{A}_k(u_h,v_h)+S_k(u_h,v_h)+J_k(u_h,v_h),\quad u_h,v_h\in V_h,
\end{equation}
are also consistent, i.e., they satisfy~(\ref{consforktrue}); furthermore, they are symmetric. These particular choices of $J_k$ (and thus $A_k$) coincide with the IPDG method of~\cite{douglas1976interior} ($k=1$) and the $h$-version of the symmetric $hp$-DG method of~\cite{suli2007hp}, with the parameters $\lambda_1=\lambda_2=1$ ($k=2$).

Analogously to deriving~(\ref{est1:10:52}) and~(\ref{est2:10:52}), one can show the following stability estimates~\cite{arnold2002unified,suli2007hp}
\begin{equation}\label{truelystable}
A_1(v_h,v_h)\ge C_1\|v_h\|_{h,1}^2\quad\forall v_h\in V_{1,h}\quad\mbox{and}\quad A_2(v_h,v_h)\ge C_2\|v_h\|_{h,\Delta}^2\quad\forall v_h\in V_{2,h},
\end{equation}
where the \emph{norms}, $\|\cdot\|_{h,1}$, and $\|\cdot\|_{h,\Delta}$ are defined by
\begin{equation}\label{needsalabel14:48}
\|v_h\|_{h,1}^2:=|v|_{H^1(\X;\Th)}^2+C_{*,1}J_1(v_h,v_h),\quad\|v_h\|_{h,\Delta}^2:=|v|_{H_\Delta(\X;\Th)}^2+C_{*,\Delta}J_2(v_h,v_h),
\end{equation}
and the constants $C_{*,1}$ and $C_{*,\Delta}$ depend only on the dimension, the domain $\X$, the polynomial degree, and the shape-regularity constants. These estimates of course yield existence and uniqueness of $u_{k,h}$ satisfying
$$A_k(u_{k,h},v_h) = \ell(v_h)\quad\forall v_h\in V_{k,h},$$
for $k=1,2$.

However, in the context of curved finite elements, it does not seem to be possible to obtain the same stability estimate for $A_2$ (i.e., the second estimate of~(\ref{truelystable})).
In the polytopal case, one may see that~(\ref{22:04}) and~(\ref{22:04:222}) follow from~(\ref{22:04-1}) due to the fact that $\Delta:\mathbb P^p(K)\to\mathbb P^{\max\{p-2,0\}}(K)\subset \mathbb P^p(K)$ for each $K\in\Th$, and so we may apply the inverse estimate~(\ref{22:04-1}) to $\Delta w_2|_K$. In the case of curved finite elements, due to the chain rule, this is no longer true, in general, since a given function of the finite element space is of the form $w|_K=\rho\circ F_K^{-1}$, where $\rho$ is a polynomial, and $F_K$ is a given (sufficiently regular) nonaffine map, and so 
\begin{equation*}
\begin{aligned}\Delta w|_K=\nabla\cdot(\nabla(\rho\circ F_K^{-1}))&=\nabla\cdot((\nabla\rho\circ F_K^{-1})(DF_K^{-1})^T)\\
&=(\nabla\rho\circ F_K^{-1})\cdot(\nabla\cdot(DF_K^{-1})^T)+(DF_K^{-1}(D^2\rho\circ F_K^{-1}))\!:\!(DF_K^{-1})^T\\
&\ne\psi\circ F_K^{-1},
\end{aligned}
\end{equation*} 
for some polynomial $\psi$, unless $F_K$ (and thus $F_K^{-1}$) is affine, i.e., the mesh is polytopal.
This leads one to obtain estimates of the form
\begin{align}
\left\|\frac{\partial (\Delta w_2)}{\partial n_F}\right\|_{L^2(F)}^2&\le C\tilde{h}_F^{-3}(|w_2|_{H^2(K)}^2+|w_2|_{H^1(K)}^2),\label{invproblem}\\
\|\Delta w_2\|_{L^2(F)}^2&\le C\tilde{h}_F^{-1}(|w_2|_{H^2(K)}^2+|w_2|_{H^1(K)}^2),
\end{align}
which would not directly lead to the derivation of the stability estimate~(\ref{truelystable}) of $A_2$ (since we are no longer able to estimate in the $\|\cdot\|_{h,\Delta}$-norm, as the Laplacian structure is no longer preserved). This leads us to define a new variant of $A_2$ with the goal of replacing the inner product $$(u,v)_{\Delta,K}:=\int_K\Delta u\,\Delta v$$ with $$\langle D^2u,D^2v\rangle_K:=\int_KD^2u\!:\!D^2v,$$ leading to coercivity in the norm 
$$\|v_h\|_{h,2}^2:=|v_h|_{H^2(\X;\Th)}^2+C_{*,2}J_2(v_h,v_h).$$
In order to achieve such a stability estimate, one is required to prove a discrete Poincar\'{e}--Friedrichs' inequality, in order to bound the $H^1$-terms of the right-hand side of~(\ref{invproblem}) by $H^2$ terms, and factors that are present in $J_2(\cdot,\cdot)$.

Finally, we discuss error estimates. Since the methods are consistent, one has
$$A_k(u_{k,h},v_h) = \ell(v_h) = A_k(u_k,v_h)\quad\forall v_h\in V_{k,h},$$
and thus, for any $z_h\in V_{k,h}$, the triangle inequality yields
\begin{equation}\label{20:00:1}
\|u_k-u_{k,h}\|_{h,k} \le \|u_k-z_h\|_{h,k}+\|u_{k,h}-z_h\|_{h,k},
\end{equation} and the stability estimates~(\ref{truelystable}) give us
\begin{equation}\label{20:00:2}
\begin{aligned}
\|u_{k,h}-z_h\|_{h,k} &\le C_k^{-1}A_k(u_{k,h}-z_h,u_{k,h}-z_h)\\
& = C_k^{-1}A_k(u_{k,h},u_{k,h}-z_h)-A_k(z_h,u_{k,h}-z_h)\\
& = C_k^{-1}A_k(u_k,u_{k,h}-z_h)-A_k(z_h,u_{k,h}-z_h)\\
& = C_k^{-1}A_k(u_k-z_h,u_{k,h}-z_h).
\end{aligned}
\end{equation}
Unfortunately $u_k-z_h$ does not, in general, belong to $V_{k,h}$, and we cannot utilise the inverse estimates that lead to the stability estimates~(\ref{truelystable}) to bound $A_k(u_k-z_h,u_{k,h}-z_h)$ in the $\|\cdot\|_{h,k}$-norm for $k=1,2$. One can, however show that~\cite{arnold2002unified,suli2007hp}
\begin{equation}\label{20:00:3}
A_k(u_k-z_h,u_{k,h}-z_h)\le \widetilde{C}_k\|u_k-z_h\|_{h,k,*}\|u_{k,h}-z_h\|_{h,k},
\end{equation}
where $\|\cdot\|_{h,k,*}$ is a variant of the $\|\cdot\|_{h,k}$ norm, including piecewise derivatives of order $0\le j\le k$. Applying~(\ref{20:00:3}) to~(\ref{20:00:2}), and applying the result to~(\ref{20:00:1}), one obtains 
$$\|u_{k,h}-u\|_{h,k}\le \|u_k-z_h\|_{h,1}+\widetilde{C}_kC_k^{-1}\|u_k-z_h\|_{h,k,*}\quad\forall z_h\in V_h.$$
Choosing $z_h\in V_{k,h}$ to be a suitable interpolant, if $u_k\in H^k_0(\X)\cap H^{2k}(\X)\cap H^{\mathbf{s}_k}(\X;\Th)$, where $\mathbf{s}_k = (s_K^k)_{K\in\Th}$, and each $s_K^k\ge 2k$, one obtains
\begin{equation}\label{errorest}
\|u_{k,h}-u_k\|_{h,k}\le C_k\left(\sum_{K\in\Th}h_K^{2t_K^k-2k}\|u_k\|_{H^{s_K^k}(K)}^2\right)^{1/2},
\end{equation}
where $t_K^k := \max\{p+1,s_K^k\}$; in the case of quasiuniform meshes, the above becomes
$$\|u_{k,h}-u\|_{h,k}\le C_kh^{\max\{p+1,s_K^k\}-k}|u_k|_{H^{s_K^k}(\X;\Th)},$$
and so the estimate is optimal with respect to the mesh size. For $k=1$, the estimate is provided in~\cite{arnold2002unified} for the case that $s_K^1=2$ for all $K\in\Th$, i.e., the integer case, and for $k=2$, the estimate~(\ref{errorest}) is provided in~\cite{suli2007hp}. In the case of curved finite elements, the method for proving optimal error estimates is the same (except there are a few more terms that we must estimate), however, one still requires a \emph{suitable interpolate}. In the context of~(\ref{errorest}), this means that there is an element $z_h\in V_{k,h}$, uniquely determined by a function $w_k\in H^{\mathbf{s}_k}(\X;\Th)$, such that for each $K\in\Th$, each \emph{integer} $0\le q\le\min\{p,2k-1\}$, and each multi-index $\alpha$, with $0\le|\alpha|\le q$,
\begin{equation}\label{interpests}
\begin{aligned}
|w_k-z_h|_{H^q(K)}&\le Ch_K^{t_K^k-q}|w_k|_{H^{s_K^k}(K)},\\
\|D^\alpha(w_k-z_h)\|_{L^2(\partial K)}&\le Ch_K^{t_K^k-q-1/2}|w_k|_{H^{s_K^k}(K)},
\end{aligned}
\end{equation}
where $C$ may depend upon the polynomial degree, $\X$, and the shape regularity constant, but is independent of $h_K$. A goal of the proceeding Section will be to prove~(\ref{interpests}) in the curved case, which will yield optimal error estimates for both the schemes we propose, and, since the polytopal case can be viewed as a special case of the curved case, we will provide optimal estimates for the IPDG method of~\cite{douglas1976interior} for the Poisson problem in noninteger Sobolev norms. The first estimate of~(\ref{interpests}) is proven in~\cite{MR1014883} for the case that $s_K^k$ is integer valued, we aim to provide such estimates in $H^s$-norms, for non integer $s$. %, and to utilise the framework of Dupont and Scott~\cite{dupont1980polynomial}, which provides a computable constant in the error estimates (the Bramble--Hilbert Lemma is nonconstructive).
\end{section}
\begin{section}{Curved domain approximation and finite element estimates}\label{sec:3}
We will begin this section by providing the details of~\cite{MR1014883}, which provides us with a notion of exact domain approximation, along with essential scaling arguments that allow us to prove the desired trace and inverse estimates. Such estimates will allow us to prove that our proposed FEMs are stable, yielding existence and uniqueness of numerical solutions. This requires the following notation.
\subsection{Notation}
\begin{definition}[Face and vertex sets]\label{faceandvertexsets:def}
Given a mesh $\Th$, we denote by $\Eib$, the set of faces of $\Th$, by $\Ei$ the set of interior faces of $\Th$, and by $\Eb$, the set of boundary faces.
\end{definition}
\begin{definition}[Jump and average operators]
For each face $F\in\Eib$, we have that $F=\overline{K}\cap\overline{K'}$ for some $K,K'\in\Th$ (in the case that $F\in\Eb$ take $F=\overline{K}\cap\partial\X$), with corresponding unit normal vector $n_F$ (which, for convention, is chosen so that it is the outward normal to $K$, we define the jump operator, $\jump{\cdot}$, over $F$ by
\begin{equation}\label{jump}
\jump{\phi} = \left\{
\begin{aligned}
&(\phi|_K)|_F-(\phi|_{K'})|_F\,\,\mbox{if}\,\,F\in\Ei,\\
&(\phi|_K)|_F\,\,\mbox{if}\,\,F\in\Eb,\\
\end{aligned}
\right.
\end{equation} 
and the average operator, $\avg{\cdot}$, by
\begin{equation}\label{average}
\avg{\phi} = \left\{
\begin{aligned}
&\frac{1}{2}((\phi|_K)|_F+(\phi|_{K'})|_F)\,\,\mbox{if}\,\,F\in\Ei,\\
&(\phi|_K)|_F\,\,\mbox{if}\,\,F\in\Eb.\\
\end{aligned}
\right.
\end{equation} 
\end{definition}
\begin{definition}[Element $L^2$-inner product]
For an element $K$, we define the inner product $\langle\cdot,\cdot\rangle_K$ by
\begin{equation}\label{C4:L2inner}
\langle u,v\rangle_K :=
\begin{aligned}
\int_Ku\,v\,\,\mbox{if}\,\,u,v\in L^2(K),
\quad\int_Ku\cdot v\,\,\mbox{if}\,\,u,v\in L^2(K;\mathbb R^d),
\quad\int_Ku:v\,\,\mbox{if}\,\,u,v\in L^2(K;\mathbb R^{d\times d}).
\end{aligned}
\end{equation} 
Any ambiguity in this notation will be resolved by the arguments of the bilinear form. The bilinear forms $\langle\cdot,\cdot\rangle_{\partial K}$ and $\langle\cdot,\cdot\rangle_{F}$ for $F\in\Eib$, are defined similarly.
\end{definition}
\begin{definition}[$\lesssim$ and $\approx$ symbols]\label{leconst}
Herein we write $a\lesssim b$ for $a,b\in\mathbb R$, if there exists a constant $C>0$, such that
$$a \le Cb,$$
independent of $\mathbf{h}:=\{h_K:K\in\Th\}$, and $u$, but otherwise possibly dependent on the polynomial degree, $p$, the shape-regularity
constants of $\Th$, $C_\mathcal{T}$, and $d$. Furthermore, we write $a\approx b$ if both $a\lesssim b$ and $b\lesssim a$.
\end{definition}
\subsection{Curved simplices}
The ability to define a nonaffine approximation of a domain, $\X\subset\mathbb R^d$ relies upon the $\X$ satisfying a notion of piecewise regularity, which motivates the following definition. 
\begin{definition}[Piecewise $C^{k}$ domain]\label{pwck21:05}
A domain $\Omega\subset\mathbb R^d$ is piecewise $C^{k}$ for $k\in\mathbb N$, if we may express the boundary of $\Omega$, $\partial\Omega$, as a finite union
\begin{equation}\label{piecewiserep}
\partial\Omega = \bigcup_{n=1}^N\overline{\Gamma_n},
\end{equation}
where each $\Gamma_n\subset\mathbb R^d$ is of zero $d$-dimensional Lebesgue measure, and admits a local representation as the {}graph of a uniformly $C^{k}$ function. That is, for each $n$, and at each $x\in\Gamma_n$ there exists an open neighbourhood $V_n$ of $x$ in $\mathbb R^d$ and an orthogonal coordinate system $(y^n_1,\ldots,y^n_d)$, such that $$V_n=\{(y^n_1,\ldots,y^n_d):-a_j^n<y_j^n<a_j^n,1\le j\le d\};$$
as well as a uniformly $C^k$ function $\varphi_n$ defined on $V'_n=\{(y_1^n,\ldots,y^n_{d-1}):-a_j^n<y_j^n<a_j^n,1\le j\le d-1\}$ and such that
\begin{align*}
&|{\varphi_n}({y^{n}}')|\le a^n_d/2\mbox{ for every }{y^n}'=(y^n_1,\ldots,y^n_{d-1})\in V_{n'},\\ 
&\X\cap V = \{y^n=({y^n}',y^n_d)\in V:y^n_d<\varphi_n({y^n}')\},\\
&\Gamma_n\cap V = \{y^n=({y^n}',y^n_d)\in V:y^n_d=\varphi_n({y^n}')\}.
\end{align*}
\end{definition}
\begin{definition}[Curved $d$-simplex]\label{curveddsimp}
An open set $K\subset\mathbb R^d$ is called a curved $d$-simplex if there exists a $C^1$ mapping $F_K$ that maps a straight reference $d$-simplex $\hat K$ onto $K$, and that is of the form
\begin{equation}\label{non:aff:map}
F_K=\tilde{F}_K+\Phi_K,
\end{equation}
where 
\begin{equation}\label{def:affine-map}
\tilde{F}_K:\hat{x}\mapsto\tilde{B}_K\hat x+\tilde{b}_K
\end{equation}
is an invertible map and $\Phi_K\in C^1(\hat{K};\mathbb R^d)$ satisfies
\begin{equation}\label{CK:def}
C_K:=\sup_{\hat{x}\in\hat K}\|D\Phi_K(\hat x)\tilde{B}_K^{-1}\|<1,
\end{equation}
where $\|\cdot\|$ denotes the induced Euclidean norm on $\mathbb R^{d\times d}$.
\end{definition}
\begin{definition}[Associated straight $d$-simplex]
Given a curved $d$-simplex $K$, with the associated straight reference $d$-simplex $\hat K$, and map $F_K:\hat K\to K$, with $F_K=\tilde{F}_K+\Phi_K$, we define the \emph{associated straight $d$-simplex}:
$$\tilde{K}:=\tilde{F}_K(\hat K).$$
\end{definition}
\begin{remark}
The associated $d$-simplex, $\tilde{K}$, is a straight $d$-simplex that ``approximates" $K$.
\end{remark}
\begin{lemma}[Affine invariance of $C_K$]\label{affineinv:17:07}
Given a $d$-simplex triple $(K,\hat K,\tilde{K})$, another reference $d$-simplex $\hat K'$, and a map $\tilde{F}_{K'}\in GL(\mathbb R^d)$ that maps $\hat K'$ onto $\hat K$, there is a map $F_{K'}:\hat K'\to K$ that also satisfies~(\ref{CK:def}). Moreover, $C_{K'}=C_K$.
\end{lemma}
\emph{Proof:} See Remark 2.3 of~\cite{MR1014883}.$\quad\quad\square$%$$F_K'(\cdot) =(\tilde{B}_K\hat A(\cdot)+\tilde{B}_K\hat a+\tilde{b}_K)+\Phi_K'(\cdot),\quad\mbox{with}\quad \Phi_K'(\cdot)=\Phi_K(\hat A(\cdot)+\hat a),$$
%which gives us
%$$C_{K'}:=\sup_{\hat x'\in\hat K'}\|D\Phi_K'(\hat x')\cdot(\tilde{B}_K\hat A)^{-1}\|=\sup_{\hat x\in\hat K}\|D\Phi_K(\hat x)\hat A\cdot\hat{A}^{-1}\tilde{B}_K^{-1}\|=C_K.\quad\quad\square$$
%\begin{remark}
%Lemma~\ref{affineinv:17:07} is rather significant, as it means that given $K\in\Th$ with an approximating straight simplex $\tilde K$, we can always redefine the map $F_K$ (without relabelling it), so that $F_K:\hat K\to K$, and $\hat K:=\tilde{K}/\diam(\tilde{K})$, and the reference simplex enjoys similar shape regularity to $\tilde K$. So, in general, we assume a given reference simplex is of this form. This fact will be used in the proof of Corollary~\ref{shaperegcorr}.
%\end{remark}
\begin{remark}[Affine mesh]
In the case that the domain has a flat boundary, one employs an affine approximation of the domain, in which case, the corresponding functions $\Phi_K$ in~(\ref{non:aff:map}) are all zero.
\end{remark}
%The (possibly curved) open  $K$ are the images of an open reference element $\hat K$ under a collection of mappings
\begin{definition}[Mesh size]\label{meshsize:def} For each $K\in\Th$, let $h_K:=\operatorname{diam}(\tilde{K})\ge C(d)\|\tilde{B}_K\|$ (where $\tilde{K}=\tilde{B}_K(\hat K)$). It is assumed that $h = \max_{K\in\Th} h_K$ for each mesh $\Th$.
\end{definition}
\begin{definition}[Face-mesh size]
For each face $F\in\mathscr{E}^{i,b}_h$, we define
\begin{equation}\label{hF:def}
\tilde{h}_F:=\left\{\begin{array}{l l}
\min(h_K,h_{K'}) & \mbox{if}~F\in\mathscr{E}^i_h,\\
h_K & \mbox{if}~F\in\mathscr{E}^b_h.
\end{array}\right.
%\quad
%\tilde{p}_F:=\left\{\begin{array}{l l}
%\min(p,p_{K'}) & \mbox{if}~F\in\mathscr{E}^i_h,\\
%p & \mbox{if}~F\in\mathscr{E}^b_h,
%\end{array}\right.
\end{equation}
where $K$ and $K'$ are such that $F=\partial K\cap\partial K'$ if $F\in\mathscr{E}^i_h$, or $F\subset\partial K\cap\partial\X$ if $F\in\mathscr{E}^b_h$. 
\end{definition}
%\emph{Mesh conditions.} We shall adopt the following assumptions on the meshes.
%\begin{assumption}\label{Meshconds}
%We assume that there is a uniform upper bound on the number of faces composing the boundary of any given element; in other words, there is a $C_\mathcal{F}>0$, independent of $h$, such that
%\begin{equation}\label{meshcond1}
%\max_{K\in\mathscr{T}_h}\operatorname{card}\{F\in\mathscr{E}^{i,b}_h:F\subset\partial K\}\le C_\mathcal{F}\quad\forall K\in\mathscr{T}_h,~\forall h>0.
%\end{equation}
%It is also assumed that any two elements sharing a face have commensurate diameters, i.e., there is a $C_{\mathcal{T}}\ge 1$, independent of $h$, such that
%\begin{equation}\label{meshcond2}
%\max(h_K,h_{K'})\le C_\mathcal{T}\min(h_K,h_{K'}),
%\end{equation}
%for any $K$ and $K'$ in $\mathscr{T}_h$ that share a face. 
%
%Finally, we assume that each $F\in\Eb$ satisfies
%\begin{equation}\label{Fcontained}
%F = F\cap\Gamma_n,
%\end{equation}
%for some $n\in\{1,\ldots,N\}$, with $\Gamma_n$ given as in~(\ref{piecewiserep}). This implies that each boundary face is completely contained in a boundary portion $\Gamma_n$.
%\end{assumption}
%\begin{remark}
%The assumptions on the mesh given by Assumption~\ref{Meshconds}, 
%%and the polynomial degrees, 
%in particular~(\ref{meshcond2}), show that if $F$ is a face of $K$, then
%\begin{equation}\label{meshcondcons}
%h_K\le C_\mathcal{T}\tilde{h}_F.
%%\quad\mbox{and}\quad\tilde{p}_F\le c_\mathcal{P}p.
%\end{equation}
%\end{remark}
\begin{definition}[Class $m$ curved $d$-simplex]
A curved $d$-simplex $K$ is of class $C^m$, $m\ge1$, if the mapping $F_K$ is of class $C^m$ on $\hat K$.
\end{definition}
The proofs of the next four lemmas can be found in~\cite{MR1014883} (i.e., Lemmas 2.1, 2.2, 2.3 and 2.4).
\begin{lemma}
The mapping $F_K$ is a $C^1$-diffeomorphism from $\hat K$ onto $K$ and satisfies
\begin{align}
\sup_{\hat x\in\hat K}\|DF_K(\hat x)\|&\le(1+C_K)\|\tilde B_K\|,\label{indinvproof1}\\
\sup_{x\in K}\|DF_K^{-1}(x)\|&\le(1-C_K)^{-1}\|\tilde B_K^{-1}\|,\label{indinvproof2}\\
\forall\hat x\in\hat K,\quad(1-C_K)^d|\det\tilde{B}_K|&\le|\det DF_K(\hat x)|\le(1+C_K)^d|\det\tilde{B}_K|.\label{indinvproof3}
\end{align}
\end{lemma}
\begin{lemma}\label{scalingassumption}
Let us denote by $c_\ell$, $2\le\ell\le m$, $m\in\mathbb N$, the constants
\begin{equation}\label{cell:def}
c_\ell(K):=\sup_{\hat x\in\hat K}\|D^\ell F_K(\hat x)\|\|\tilde{B}_K\|^{-\ell}.
\end{equation}
There exist constants $c_{-\ell}$, $2\le\ell\le m$, depending continuously on $c_K$, $c_2(K),\ldots,c_m(K)$, such that
\begin{equation}\label{indinvproof4}
\sup_{x\in K}\|D^\ell F_K^{-1}(x)\|\le c_{-\ell}\|\tilde{B}_K\|^{2(\ell-1)}\|\tilde{B}_K^{-1}\|^\ell.
\end{equation}
\end{lemma}
\subsection{Scaling arguments}
\begin{lemma}\label{inverselemma}
Assume that $K$ is a curved $d$-simplex of class $C^m$. Let $l$ be an integer, $0\le l\le m$, and $q\in\{2,\infty\}$. A function $v$ belongs to $W^{m,q}(K)$ if and only if the function $\hat v:=v\circ F_K$ belongs to $W^{m,q}(\hat K)$. We also have for any $v\in W^{m,q}(K)$
\begin{align}
|v|_{W^{l,q}(K)}&\le C|\operatorname{det}\tilde{B}_K|^{1/q}\|\tilde{B}_K^{-1}\|^l\left(\sum_{r=\min\{l,1\}}^l\|\tilde{B}_K\|^{2(l-r)}|\hat v|_{W^{r,q}(\hat K)}\right),\label{nonaff:scaling}\\
|\hat v|_{W^{l,q}(\hat K)}&\le C|\operatorname{det}\tilde{B}_K|^{-1/q}\|\tilde{B}_K\|^l\left(\sum_{r=\min\{l,1\}}^l|v|_{W^{r,q}(K)}\right),\label{nonaffine:scaling}
\end{align}
where the constants $C$ depend continuously on $c_K$, $c_2(K),\ldots,c_m(K)$.
\end{lemma}
\begin{lemma}\label{tracelemma}
Assume that $K$ is a curved $d$-simplex of class $C^m$, and that $F$ is a face of $K$; we denote by $\tilde{B}_F$ the restriction of $\tilde{B}_K$ to $\hat F:=F_K^{-1}(F)$. Let $l$ be an integer, $1\le l\le m$, $s\in[0,l-1/2)$. Then, for any $v\in H^l(K)$, the function $\tau_F(v)$ belongs to $H^{s}(F)$, and we have
\begin{equation}\label{tracelikebound}
\|v\|_{H^s(F)}\le C|\operatorname{det}\tilde{B}_F|^{1/2}|\operatorname{det}\tilde{B}_K|^{-1/2}\|\tilde{B}_K^{-1}\|^s(\|v\|_{L^2(K)}+\|\tilde{B}_K\|^l|v|_{H^{l}(K)}),
\end{equation}
where the constant $C$ depends continuously on $c_K$, $c_2(K),\ldots c_m(K)$.
\end{lemma}
%\begin{remark}[Interpretation of the Lemmas~\ref{scalingassumption}--\ref{tracelemma}]
%One may view the seemingly rather abstract Lemmas~\ref{scalingassumption},~\ref{inverselemma}, and~\ref{tracelemma} as a necessary prerequisite for the standard scaling argument used to prove optimal interpolation estimates on affine meshes, inverse inequalities in $W^{m,p}$-norms, and trace inequalities in $H^s$-norms, respectively. In particular, the following example should demonstrate the importance of having a uniform upper bound on the value $c_\ell$.
%\end{remark}
A key tool in the derivation of optimal interpolation estimates on affine meshes is the following scaling argument (see Theorem 3.1.2 of~\cite{MR1930132}): for $l\in\mathbb N_0$, $p\in[1,\infty]$, assuming $v\in W^{l,p}(\tilde K)$, and $\hat v := v\circ F_K\in W^{l,p}(\hat K)$, we have
\begin{equation}\label{affine:scaling}
|\hat v|_{W^{l,p}(\hat K)}\le C\|\tilde B_K\|^l|\operatorname{det}\tilde{B}_K|^{-1/p}|v|_{W^{l,p}(\tilde{K})}.
\end{equation}
Here, we are considering the \emph{affine equivalent} straight $d$-simplices $\hat K$ and $\tilde{K}$, and an invertible affine map $F_K$. That is, $\tilde{K} = F_K(\hat K)$, where $F_K$ is of the form~(\ref{non:aff:map}) with $\Phi_K\equiv0$. 

One can see that~(\ref{affine:scaling}) and~(\ref{nonaffine:scaling}) are similar. The main difference is the presence of the lower order seminorms on the right-hand side of~(\ref{nonaffine:scaling}).

Let us look at the particular example of the $H^2$-seminorm when $F_K$ is not affine. The chain rule, and the multivariable change of variables formula yields
%\begin{equation}\label{abuseofnot}
%\begin{aligned}
%|\hat v|_{H^2(\hat K)}^2 &= \int_{\hat K}|D^2\hat v|^2\,\\
%%& = \int_{\hat K}|\nabla(D(v\circ F_K))|^2\,\\
%%& = \int_{\hat K}|\nabla((Dv\circ F_K)DF_K)|^2\\
%%& = \int_{\hat K}\left|[D((Dv\circ F_K)DF_K)]^T\right|^2\\
%%& = \int_{\hat K}\left|[D^2v\circ F_K(DF_K)^2+Dv\circ F_KD^2F_K]^T\right|^2\\
%%& = \int_{\hat K}\left|((DF_K)^2)^T(D^2v)^T\circ F_K+(D^2F_K)^T(Dv\circ F_K)^T\right|^2\\
%& = \int_{\hat K}\left|((DF_K)^2)^T(D^2v\circ F_K)+(D^2F_K)^T(\nabla v\circ F_K)\right|^2\\
%&\le C(d)\int_{\hat K}\|DF_K\|^4|D^2v\circ F_K|^2+\|D^2F_K\|^2|\nabla v\circ F_K|^2\\
%&\le C(d)\left(\sup_{\hat x\in\hat K}\|DF_K(\hat x)\|^4\int_{\hat K}|(D^2v)\circ F_K|^2\right.\\
%&\left.~~~~~~~~~~~~~~~~~~+\sup_{\hat x\in\hat K}\|D^2F_K(\hat x)\|^2\int_{\hat K}|(\nabla v)\circ F_K|^2\right)\\
%&\le C(d)\sup_{x\in K}|\operatorname{det}DF_K^{-1}(x)|\left(\sup_{\hat x\in\hat K}\|D^2F_K(\hat x)\|^2|v|^2_{H^1(K)}\right.\\
%&\left.~~~~~~~~~~~~~~~~~~~~~~~~~~~~~~~~~~~~~~~~~~~~~+\sup_{\hat x\in\hat K}\|DF_K(\hat x)\|^4|v|^2_{H^2(K)}\right).
%\end{aligned}
%\end{equation}
%Thus, taking square roots, we obtain
\begin{equation}\label{non:aff:bound:12:59}
\begin{aligned}
|\hat v|_{H^2(\hat K)} &\le C(d)\sup_{x\in K}|\operatorname{det}DF_K^{-1}(x)|^{1/2}\left(\sup_{\hat x\in\hat K}\|D^2F_K(\hat x)\||v|_{H^1(K)}+\sup_{\hat x\in\hat K}\|DF_K(\hat x)\|^2|v|_{H^2(K)}\right).
\end{aligned}
\end{equation}
Note that if $F_K$ were affine, then $DF_K=\tilde{B}_K$, $DF_K^{-1} = \tilde{B}_K^{-1}$, and $D^2F_K\equiv0$, thus from the above, we immediately obtain~(\ref{affine:scaling}) with $l=p=2$.

A sufficient assumption that yields an estimate of the same order as~(\ref{affine:scaling}) with $l=p=2$ (in terms of $\|\tilde{B}_K\|$), is to assume that $c_\ell$, given by~(\ref{cell:def}), is uniformly bounded for $\ell=2$. This, coupled with the fact that $C_K<1$ gives us
\begin{align*}
\sup_{\hat x\in\hat K}\|DF_K(\hat x)\|&\le(1+C_K)\|\tilde{B}_K\|,\\
\sup_{\hat x\in\hat K}\|D^2F_K(\hat x)\|&=(\sup_{\hat x\in\hat K}\|D^2F_K(\hat x)\|\tilde{B}_K\|^{-2})\|\tilde{B}_K\|^2=c_2\|\tilde{B}_K\|^2.
\end{align*}
Applying the above to~(\ref{non:aff:bound:12:59}) yields
$$|\hat v|_{H^2(\hat K)}\le C\sup_{x\in K}|\operatorname{det}DF_K^{-1}(x)|^{1/2}\|\tilde{B}_K\|^2(|v|_{H^1(K)}+|v|_{H^2(K)}).$$
In order to appropriately bound the determinant term, one must note that
$DF_K^{-1} = (DF_K)^{-1},$
and so
$$|\operatorname{det}DF_K^{-1}|=|\operatorname{det}DF_K|^{-1}\le|\operatorname{det}\tilde{B}_K|^{-1}(1-C_K)^d.$$
Ultimately, this gives us
\begin{equation}\label{*normmotiv}
|\hat v|_{H^2(\hat K)}\le C|\operatorname{det}\tilde{B}_K|^{-1/2}\|\tilde{B}_K\|^2(|v|_{H^1(K)}+|v|_{H^2(K)}).
\end{equation}
This motivates the two following definitions, generalising the prerequisite assumptions, allowing one to obtain analogous estimates in higher order seminorms.
%\begin{remark}[Minor abuse of notation]
%One may notice that from the seventh equality of~(\ref{abuseofnot}) the ``$D^2$" applied to $v$ is different to the ``$D^2$" applied to $F_K$. In particular, the first ``$D^2$" is the Hessian, and the latter ``$D^2$" is the actual second derivative. This stems from the fact that we define the Hessian of $v$, $D^2v$, in the following way:
%$$D^2v = \nabla(Dv) = [D(Dv)]^T = [D^2v]^T,$$
%where the final $D^2$ above is the actual second derivative.
%\end{remark}
\begin{definition}
The family $(\Th)_h$ of meshes is said to be regular if there exist two constants, $\sigma$ and $c$, independent of $h$, such that, for each $h$, any $K\in\Th$ satisfies
\begin{equation}\label{shape-reg}
h_K/\rho_K\le\sigma,
\end{equation}
where $\rho_K$ is the diameter of the sphere inscribed in $\tilde{K}$.
Furthermore, we have
\begin{equation}\label{regularpart}
\sup_h\sup_{K\in\Th} C_K\le c<1.
\end{equation}
\end{definition}
\begin{remark}
Condition~(\ref{shape-reg}) is referred to as \emph{nondegeneracy} (for example in~\cite{MR2373954}).
\end{remark}
\begin{definition}\label{regoforderm}
The family $(\Th)_h$ of meshes is said to be regular of order $m$ if it is regular and if, for each $h$, any $K\in\Th$ is of class $C^{m+1}$, with
\begin{equation}\label{regoforderms}
\sup_h\sup_{K\in\Th}\sup_{\hat x\in\hat K}\|D^lF_K(\hat x)\|\|\tilde{B}_K\|^{-l}<\infty,\quad 2\le l\le m+1.
\end{equation}
\end{definition}
%\begin{definition}[Curved and flat partition of $\Th$]\label{candfpart}
%Given a mesh $\Th$, we denote by $\Thf$ the set of all $K\in\Th$ that satisfy $D^2F_K\equiv0$, and define $\Thc:=\Th\setminus\Thf$.
%\end{definition}
%\begin{definition}[Curved and flat partition of $\Eb$]\label{candfpartEb}
%Given a mesh $\Th$, we denote by $\Ebc$ the set of all faces $F\in\Eb$ that satisfy $F\subset\partial K$, where $K\in\Thc$, and define $\Ebf:=\Eb\setminus\Ebc$.
%\end{definition}
\begin{assumption}\label{Meshconds}
%We assume that there is a uniform upper bound on the number of faces composing the boundary of any given element; in other words, there is a $C_\mathcal{F}>0$, independent of $h$, such that
%\begin{equation}\label{meshcond1}
%\max_{K\in\mathscr{T}_h}\operatorname{card}\{F\in\mathscr{E}^{i,b}_h:F\subset\partial K\}\le C_\mathcal{F}\quad\forall K\in\mathscr{T}_h,~\forall h>0.
%\end{equation}
We assume that any two elements sharing a face have commensurate diameters, i.e., there is a $C_{\mathcal{T}}\ge 1$, independent of $h$, such that
\begin{equation}\label{meshcond2}
\max(h_K,h_{K'})\le C_\mathcal{T}\min(h_K,h_{K'}),
\end{equation}
for any $K$ and $K'$ in $\mathscr{T}_h$ that share a face. 

Finally, we assume that each $F\in\Eb$ satisfies
\begin{equation}\label{Fcontained}
F = F\cap\Gamma_n,
\end{equation}
for some $n\in\{1,\ldots,N\}$, with $\Gamma_n$ given as in~(\ref{piecewiserep}). This implies that each boundary face is completely contained in a boundary portion $\Gamma_n$, as well as ensuring that our approximation of the domain $\X$ is exact.
\end{assumption}
\begin{remark}
The assumptions on the mesh given by Assumption~\ref{Meshconds}, 
%and the polynomial degrees, 
in particular~(\ref{meshcond2}), show that if $F$ is a face of $K$, then
\begin{equation}\label{meshcondcons}
h_K\le C_\mathcal{T}\tilde{h}_F.
%\quad\mbox{and}\quad\tilde{p}_F\le c_\mathcal{P}p.
\end{equation}
\end{remark}
A final, necessary step, before providing optimal interpolation estimates and inverse estimates for (continuous and discontinuous) curved Lagrange finite element spaces, is to relate the estimates of this section to the local mesh size, $h_K$. The general rule of thumb in this context is that $\|\tilde{B}_K\|$ is of order $h_K$, and $\|\tilde{B}_K^{-1}\|$ is of order $h_K^{-1}$. This notion is made more concrete by the following theorem from~\cite{MR1930132}.
\begin{theorem}\label{ciarlet:thm}
Let $\hat K$ and $\tilde{K}=\tilde{F}_K(\hat K)$ be two affine-equivalent open subsets of $\mathbb R^d$, where $\tilde{F}_K:\hat x\to\tilde{B}_K\hat x+\tilde{b}_K$ is an invertible affine mapping. Then we have the upper bounds
\begin{equation}\label{preB_Khbehaviour}
\|\tilde{B}_K\|\le\frac{h(\tilde K)}{\rho(\hat K)}\quad\mbox{and}\quad\|\tilde{B}_K^{-1}\|\le\frac{h(\hat K)}{\rho(\tilde K)},
\end{equation}
where, for a given open subset $E$ of $\mathbb R^d$, we define
\begin{equation}\label{ciarletessential:thm}
\begin{aligned}
&h(E) = \operatorname{diam}(E),\\
&\rho(E) = \sup\{\operatorname{diam}(S):S\mbox{ is a ball contained in E}\}.
\end{aligned}
\end{equation}
\end{theorem}
\begin{corollary}\label{shaperegcorr}
Assume that the family $(\Th)_h$ of meshes satisfies~(\ref{shape-reg}). % and let $K^*$ be a fixed reference $d$-simplex. 
Then, there exists a positive constant $C$ depending only on $\sigma$, such that for any $K\in\Th$ with an associated straight element $\tilde K$, that
\begin{equation}\label{B_Khbehaviour}
\|\tilde{B}_K\|\le Ch_K\quad\mbox{and}\quad\|\tilde{B}_K^{-1}\|\le Ch_K^{-1}.
\end{equation}
\end{corollary}
\emph{Proof:} Firstly, by Lemma~\ref{affineinv:17:07}, we may, without loss of generality, assume %that $0\in\tilde{K}$, 
%and
that the reference simplex $\hat K\subset B_1(0)\subset\mathbb R^d$, and that $\rho(\hat K)\ge C^{-1}$, where the constant $C>0$ depends upon $\sigma$, but is independent of $\tilde K$ (see~\cite{MR1014883}). Since $\tilde{K}$ and $\hat K$ are affine equivalent, we may apply~(\ref{preB_Khbehaviour}), which gives us
%Recall that for a given $K\in\Th$, we associate the map $F_K=\tilde{F}_K+\Phi_K:\hat K\to K$, with $\tilde{F}_K(\hat x) = \tilde{B}_K\hat x+\tilde{b}_K$, $\hat x\in\hat K$, where $\tilde{B}_K$ is an invertible matrix, and $\tilde{b}_K\in\mathbb R^d$. 
\begin{equation*}
\|\tilde{B}_K\|\le\frac{h(\tilde K)}{\rho(\hat K)}\quad\mbox{and}\quad\|\tilde{B}_K^{-1}\|\le\frac{h(\hat K)}{\rho(\tilde K)}.
\end{equation*}
Recall that we define $h_K:=h(\tilde{K})$, and $\rho_K:=\rho(\tilde{K})$ and so we have $\|\tilde{B}_K\|\le h_K/\rho(\hat K)\le Ch_K$, which gives us the first estimate of~(\ref{B_Khbehaviour}). By~(\ref{shape-reg}), we also have $\|\tilde{B}_K^{-1}\|\le h(\hat K)/\rho(\tilde{K})\le\sigma h(\hat K)h_K^{-1}\le\sigma h_K^{-1}$, which is the second estimate of~(\ref{B_Khbehaviour}).$\quad\quad\square$
%. To obtain the first estimate of~(\ref{B_Khbehaviour}), we note that 
%\begin{equation*}
%\begin{aligned}
%\rho(\hat K) &= \sup\{\diam(S):S\mbox{ is a ball contained in }\hat K\}\\
%&=\sup\left\{\diam(S):S\mbox{ is a ball contained in }\frac{1}{\diam\tilde{K}}\tilde K\right\}\\
%&=\frac{1}{\diam\tilde{K}}\rho(\tilde K)\ge\frac{h(\tilde K)}{\sigma h(\tilde K)}=\frac{1}{\sigma}\quad(\mbox{by }\ref{shape-reg}).
%\end{aligned}
%\end{equation*}
%Thus $\|\tilde{B}_K\|\le \sigma h_K$.
%
\begin{definition}[$v$, $\hat v$, and $v^*$]\label{vvv:def}
Given a triple $(K^*,\hat K,K)$ (fixed reference simplex, reference simplex, and curved simplex), a pair of invertible maps $(G_K:K^*\to\hat K,F_K:\hat K\to K)$, and a function $v:K\to\mathbb R^N$, for some $N\in\mathbb N$, we define the functions $\hat v:\hat K\to\mathbb R^N$, $v^*:K^*\to\mathbb R^N$, as follows:
\begin{equation}\label{starhatetc}
\hat v:=v\circ F_K,\quad v^*:=\hat v\circ G_K=v\circ F_K\circ G_K.
\end{equation}
Furthermore, given $v^*:K^*\to\mathbb R^N$, we also define
\begin{equation}\label{starhatetotherway}
\hat v:=v^*\circ G_K^{-1},\quad v:=\hat v\circ F_K^{-1}=v^*\circ G_K^{-1}\circ F_K^{-1}.
\end{equation}
\end{definition}
\subsection{Lagrange finite element spaces}
The finite element spaces we consider in this paper consist of \emph{discontinuous} piecewise polynomial functions, which fall into the class of discontinuous (curved) Lagrange finite element spaces. In general, a finite element is a triple $(K,P_K,\Sigma_K)$ where $K$ is a subset of $\mathbb R^d$, $P_K$ is a finite dimensional space on $K$, and $\Sigma_K$ is a set of continuous linear forms on $P_K$, which we will call the degrees of freedom. In the context of Lagrange finite element spaces, the continuous linear forms are given by (local) point evaluations. In the simplicial case, the placement of these points is naturally described using the barycentric coordinates of the simplex.
\begin{definition}[Barycentric coordinates]\label{barydef}
Given a straight $d$-simplex $\hat K$, with vertices $\hat a_1,\ldots,$\\$\hat a_{d+1}\in\mathbb R^d$, we define the barycentric coordinates of $\hat K$, $\lambda_1,\ldots,\lambda_d,\lambda_{d+1}$ via the following (invertible) system
\begin{equation}
\begin{bmatrix}
1 & 1 & \ldots & 1\\
(\hat a_1)_1 & (\hat a_2)_1 & \ldots & (\hat a_{d+1})_1\\
\vdots & \vdots  & \ddots & \vdots \\
(\hat a_{1})_d & (\hat a_{2})_d & \ldots & (\hat a_{d+1})_d\\
\end{bmatrix}
\begin{bmatrix}
\lambda_1\\
\lambda_2\\
\vdots\\
\lambda_{d+1}
\end{bmatrix}
=
\begin{bmatrix}
1\\
\hat x_1\\
\vdots\\
\hat x_d
\end{bmatrix}
,
\end{equation}
where $\hat x = (x_1,\ldots,x_d)^T\in\hat K$.
\end{definition}
\begin{definition}[Straight Lagrange finite element]\label{straightlagrdef}
For a straight $d$-simplex $\hat K$ with vertices $\hat a_1,\ldots,\hat a_{d+1}\in\mathbb R^d$, with barycentric coordinates $\lambda_1,\ldots,\lambda_{d+1}$, we set 
\begin{equation}
J(p) = \{\alpha\in\mathbb N_0^{d+1}:|\alpha|=p\},
\end{equation}
and for any $\alpha\in J(p)$, we associate the point $\hat a_\alpha\in\hat K$ with barycentric coordinates $\lambda_i=\alpha_i/p$, $i=1,\ldots d+1$. Then, we call $(\hat K,\hat P_K,\hat \Sigma_K)$ a straight Lagrange finite element of type $p$, where
\begin{equation}\label{lineardofs}
\hat P_K = \mathbb P^p(\hat K),\quad\hat\Sigma_K = \{\hat\mu_\alpha,\alpha\in J(p)\},
\end{equation}
with $\hat\mu_\alpha(\hat f):=\hat f(\hat a_\alpha)$, for $f\in\hat P_K$, and we recall that $\mathbb P^{p}(K)$ is the space of all polynomials with total degree less than or equal to $p$.
\end{definition}
\begin{definition}[Curved Lagrange finite element]\label{curvedlagrdef}
The triple $(K,P_K,\Sigma_K)$ is a curved Lagrange finite element of type $(m,p)$ if $K$ is a curved $d$-simplex of class $C^{m+1}$, and
\begin{align}
P_K&=\{\rho=\hat\rho\circ F_K^{-1},\hat\rho\in\hat P_K=\mathbb P^{p}(\hat K)\},\\
\Sigma_K&=\{\mu:\forall v\in C^0(\overline{K}),\mu(v) = \hat\mu(v\circ F_K),\,\hat\mu\in\hat\Sigma_K\},\label{nonlineardofs}
\end{align}
where $(\hat K,\hat P_K,\hat\Sigma_K)$ is a straight Lagrange finite element of type $p$.
\end{definition}
\begin{definition}[Discontinuous Galerkin finite element space]\label{dgspaces:def} The discontinuous Galerkin finite element space $\dg$ is defined by
\begin{equation}\label{FEspace:def}
\dg:=\{v\in L^2(\X):~v|_K=\hat\rho\circ F_K^{-1},\,\hat\rho\in\mathbb P^{p}(\hat K),~\forall K\in\Th\},
\end{equation}
where $p\in\mathbb N_0$.
\end{definition}
\begin{remark}
One could equivalently define $\dg:=\cup_{K\in\Th}P_K$, where $P_K$ is a curved Lagrange finite element of type $(m,p)$.
\end{remark}
%\begin{remark}[Degrees of freedom]
%As mentioned before, we call the collection $\Sigma_K$ the ``degrees of freedom" of a finite element. One reason for this is the notion of unisolvency, that is, the collection $\Sigma_K$ should uniquely determine the basis functions of $P_K$, and thus determines all of the ``degrees of freedom" that a function in $P_K$ may have. Furthermore, in the case of DG vs CG, one should not call the vertices of $K$ the degrees of freedom, particularly when we move from a local to a global setting. This may go unnoticed in the continuous case, where the set of all vertices of all $K\in\Th$ is in a one-to-one correspondence with all of the linear forms in $\cup_{K\in\Th}\Sigma_K$. However, in the discontinuous case this is no longer true, in the sense that the degrees of freedom associated with face vertices are no longer ``shared" between neighbouring elements (indeed given a vertex shared between two elements, the point evaluation of a function $v\in\dg$ at such a vertex is not always well defined).
%\end{remark}
Piecewise polynomial functions naturally satisfy a property of piecewise regularity. This is accurately captured by considering the notion of broken Sobolev spaces.

\begin{definition}[Broken Sobolev spaces] Let $\mathbf s=(s_K:K\in\Th)$ denote a vector of nonnegative real numbers and let $r\in[1,\infty]$.
The broken Sobolev space $W^{\mathbf s,r}(\X;\Th)$ is defined by
\begin{equation}\label{brokensob}
W^{\mathbf s,r}(\X;\Th):=\{v\in L^2(\X):~v|_K\in W^{s_K,r}(K)~\forall K\in\Th\}.
\end{equation}
We denote $H^\mathbf{s}(\X;\Th):=W^{\mathbf s,2}(\X;\Th)$, and set $W^{s,r}(\X;\Th):=W^{\mathbf s,r}(\X;\Th)$, in the case that $s_K=s,~s\ge 0$, for all $K\in\Th$. For $v\in W^{1,r}(\X;\Th)$, let 
$\nabla_hv\in L^r(\X;\mathbb R^d)$ denote the discrete (also known as broken) gradient of $v$, i.e., $(\nabla_hv)|_K=\nabla(v|_K)$ for all $K\in\Th$. Higher order discrete derivatives are defined in a similar way. We define a norm on $W^{s,r}(\X;\Th)$ by
\begin{equation}\label{brokensobnorm}
\|v\|^r_{W^{s,r}(\X;\Th)}:=\sum_{K\in\Th}\|v\|^r_{W^{s,r}(K)}
\end{equation}
with the usual modification when $r=\infty$.
\end{definition}
\begin{definition}
We define the following for $K\in\Th$, $s\in\mathbb N_0$, $r\in[1,\infty)$:
\begin{align}
|v|^r_{W_*^{s,r}(K)}&:=\sum_{j=\min\{1,s\}}^s|v|^r_{W^{j,r}(K)},\label{brokensob*normK}\\
|v|^r_{W_*^{s,r}(\X;\Th)}&:=\sum_{K\in\Th}|v|^r_{W_*^{s,r}(K)},\label{brokensob*norm}
\end{align}
with the usual modification when $r=\infty$. Note that $|\cdot|_{W_*^{s,r}(\X;\Th)}$ is a norm when $s=0$, and a semi-norm when $s\in\mathbb N$. We also define $|\cdot|_{H_*^s(K)}$ and $|\cdot|_{H_*^s(\X;\Th)}$ in the usual way.
\end{definition}
\begin{remark}
We can use these semi-norms to equivalently phrase estimates such as~(\ref{*normmotiv}), which can now be written as
$$
|\hat v|_{H^2(\hat K)}\le C|\operatorname{det}\tilde{B}_K|^{-1/2}\|\tilde{B}_K\|^2|v|_{H_*^2(K)}.
$$
\end{remark}
We now state and prove trace and inverse estimates that we will be utilised frequently. In particular, the noninteger order trace estimate will be utilised in proving the second estimate of~(\ref{interpests}).
\subsection{Trace and inverse estimates}
\begin{lemma}
Assume that $\Th$ is a regular mesh on $\overline{\X}$. Then, for any $K\in\Th$, we have that
\begin{equation}\label{1simptrace}
\|v\|^2_{2,\partial K}\le C_{\operatorname{Tr}}(h_K^{-1}\|v\|_{2,K}^2+h_K\|\nabla v\|_{2,K}^2)\quad\forall v\in H^1(K),
\end{equation}
where $C_{\operatorname{Tr}}$ is independent of $K$ and $h_K$.
\end{lemma}
\emph{Proof:} Applying~(\ref{tracelikebound}) of Lemma~\ref{tracelemma} with $l=m=1$ and $s=0$, for any $K\in\Th$ and any face $F$ of $K$, we obtain
$$\|v\|_{L^2(F)}^2\le C|\operatorname{det}\tilde{B}_F||\operatorname{det}\tilde{B}_K|^{-1}(\|v\|_{L^2(K)}^2+\|\tilde{B}_K\|^{2}|v|_{H^1(K)}^2),$$
where we recall that $\tilde{B}_F$ is the restriction of $\tilde{B}_K$ to $\hat F:=F_K^{-1}(F)$ (and thus acts as a map on $\mathbb R^{d-1}$). Now, applying~(\ref{B_Khbehaviour}) yields
$$\|\tilde{B}_K\|\le Ch_K\quad\mbox{and}\quad\|\tilde{B}_K^{-1}\|\le Ch_K^{-1},$$
where the constant $C$ is independent of $K$ and $h_K$.
Thus, as the determinant is a continuous $d$-linear ($(d-1)$-linear in the case of $\tilde{B}_F$) map, we obtain
\begin{equation}\label{trace:proof}
\begin{aligned}
\|v\|_{L^2(F)}^2&\le C\|\tilde{B}_K\|^{d-1}\|\tilde{B}_K^{-1}\|^d(\|v\|_{L^2(K)}^2+\|\tilde{B}_K\|^{2}|v|_{H^1(K)})\\
&\le Ch_K^{-1}(\|v\|_{L^2(K)}^2+h_K^{2}|v|_{H^1(K)})=C(h_K^{-1}\|v\|_{L^2(K)}^2+h_K|v|_{H^1(K)}).
\end{aligned}
\end{equation}
Since the number of faces of an element $K\in\Th$ is uniformly bounded with respect to the dimension, we obtain~(\ref{1simptrace}) by summing~(\ref{trace:proof}) over all faces $F\subset\partial K$. $\quad\quad\square$
\begin{lemma}[Noninteger order trace estimate]\label{fractracelem}
Assume that $\{\Th\}_h$ is a regular family of triangulations on $\overline{\X}$. Then, for any $K\in\Th$, and any $(d-1)$ face $F$ of $K$, we have that
\begin{equation}\label{tracestfrac}
\|v\|_{L^2(F)}\le  Ch_K^{-1/2}(\|v\|_{L^2(K)}+h_K^r|v|_{H^r(K)}),
\end{equation}
for all $v\in H^r(K)$, $1/2<r<1$. Furthermore, the constant $C$ is independent of $h_K$ and the choice of $K\in\Th$.
\end{lemma}
\emph{Proof:} From the multivariable change of variables formula, we obtain 
$$\|v\|_{L^2(F)}\le C|\det\tilde{B}_F|^{1/2}\|\hat v\|_{L^2(\hat F)},$$
where $\tilde{B}_F$ is the restriction of $\tilde{B}_K$ to $\hat F=F_K^{-1}(F)$.
Under a second change of variables, we obtain
$$\|\hat v\|_{L^2(\hat F)}=|\det{\tilde{A}_{\hat F}}|^{1/2}\|v^*\|_{L^2(F^*)},$$
where $F^*$ is a $(d-1)$-face of a \emph{fixed} reference $d$-simplex, $K^*$, and $G_K:K^*\to\hat K$, $G_K(x^*):=\tilde{A}_{\hat K}x^*+\tilde{a}_{\hat K}$, with $\tilde{A}_{\hat K}\in GL(\mathbb R^d)$, $\tilde{a}_{\hat K}\in\mathbb R^d$,
and $\tilde{A}_{\hat F}$ is the restriction of $\tilde{A}_{\hat K}$ to $F^*=G_K^{-1}(\hat F)$.

Since the trace operator is continuous from $H^r(K^*)\to L^2(K^*)$ for $r>1/2$~\cite{MR3396210}, we see that
\begin{equation}\label{16:33:wed}
\begin{aligned}
\|v^*\|_{L^2(F^*)}&\le C(K^*,d)(\|v^*\|_{L^2(K^*)}+|v^*|_{H^r(K^*)})\\
&\le C(K^*,d)(\chi_1(A_K)\|\hat v\|_{L^2(\hat K)}+\chi_2(A_K)|\hat v|_{H^r(\hat K)}),
\end{aligned}
\end{equation}
where $\chi_1$, and $\chi_2$ are positive, continuous functions that we will soon provide.

Recall the definition of the $H^r$-semi norm:
\begin{equation}\label{hrseminorm}
|\hat v|_{H^r(\hat K)}^2:=\int_{\hat K}\int_{\hat K}\frac{|\hat v(\hat x_1)-\hat v(\hat x_2)|^2}{|\hat x_1-\hat x_2|^{d+2r}}.
\end{equation}
We note that since $\hat x_1,\hat x_2\in\hat K$,
$$|F_K(\hat x_1)-F_K(\hat x_2)|\le C(d)\sup_{\hat x\in\hat K}\|DF_K(\hat x)\||\hat x_1-\hat x_2|,$$
which, when applied to~(\ref{hrseminorm}), gives us
\begin{equation*}
\begin{aligned}
\int_{\hat K}\int_{\hat K}\frac{|\hat v(\hat x_1)-\hat v(\hat x_2)|^2}{|\hat x_1-\hat x_2|^{d+2r}} & = \int_{\hat K}\int_{\hat K}\frac{(C(d)\sup_{\hat x\in\hat K}\|DF_K(\hat x)\|)^{d+2r}|\hat v(\hat x_1)-\hat v(\hat x_2)|^2}{(C(d)\sup_{\hat x\in\hat K}\|DF_K(\hat x)\||\hat x_1-\hat x_2|)^{d+2r}}\\
\end{aligned}
\end{equation*}
\begin{equation}\label{scalingarg1}
\begin{aligned}
&\le\int_{\hat K}\int_{\hat K}\frac{(C(d)\sup_{\hat x\in\hat K}\|DF_K(\hat x)\|)^{d+2r}|\hat v(\hat x_1)-\hat v(\hat x_2)|^2}{|F_K(\hat x_1)-F_K(\hat x_2)|^{d+2r}}\\
&\le C\|\tilde B_K\|^{d+2r}\int_{\hat K}\int_{\hat K}\frac{|\hat v(\hat x_1)-\hat v(\hat x_2)|^2}{|F_K(\hat x_1)-F_K(\hat x_2)|^{d+2r}}.
\end{aligned}
\end{equation}
We apply the multivariable change of variables formula once more, obtaining
\begin{equation}\label{scalingarg2}
\begin{aligned}
|\hat v|_{H^r(\hat K)}^2&\le C\|\tilde B_K\|^{d+2r}\int_K\int_K\frac{|v(x_1)-v(x_2)|^2}{|x_1-x_2|^{d+2r}}|\det(DF_K^{-1}(x_1))||\det(DF_K^{-1}(x_2))|\\
&\le C\|\tilde B_K\|^{d+2r}\|\tilde B_K^{-1}\|^{2d}|v|_{H^r(K)}^2.
\end{aligned}
\end{equation}
Of course, we also have
$$\|\hat v\|_{L^2(\hat K)}\le C\|\tilde B_K^{-1}\|^d\|v\|_{L^2(K)}.$$
We obtain the functions $\chi_1$ and $\chi_2$ in a similar manner, except since $G_K$ is affine, the scaling argument is simpler, and we have that
$$\chi_1(A)=|\det A^{-1}|,\quad\mbox{and}\quad\chi_2(A)=|\det A^{-1}|^2\|A\|^{d+2r}.$$
From the nondegeneracy condition~(\ref{shape-reg}), it follows (from the proof of Theorem 4.4.20 in~\cite{MR2373954}) that the collection of the invertible matrices given by the affine maps from $K^*$ to $\hat K$  is contained in a \emph{compact} subset $BL:=\{B\in GL(\mathbb R^d):|\det B|\ge\varepsilon,|B_{ij}|\le r\}$ of $GL(\mathbb R^d)$, where $\varepsilon=\varepsilon(\sigma,d,K^*)$, and $r=r(K^*)$. That is, if
\begin{equation*}
\begin{aligned}
\tilde{G}_{\hat K}:&\,K^*\to\hat K,\quad K^*\ni x^*\mapsto \tilde{A}_{\hat K}x^*+\tilde{a}_{\hat K}\in\hat K,
\end{aligned}
\end{equation*}
then $\tilde{A}_{\hat K}\in BL$. Thus we have
$$\chi_i(\tilde{A}_{\hat K})^2\le\sup_{A\in BL}\chi_i(A)^2\le C(K^*,\sigma),\quad i=1,2.$$
Overall, we have obtained
\begin{equation*}
\begin{aligned}
\|v\|_{L^2(F)}&\le C(d,\sigma,K^*)|\det\tilde{B}_F|^\frac{1}{2}\|\tilde{B}_K^{-1}\|^\frac{d}{2}(\|v\|_{L^2(K)}+\|\tilde B_K^{-1}\|^\frac{d}{2}\|\tilde{B}_K\|^\frac{d}{2}\|\tilde{B}_K\|^r|v|_{H^r(K)})\\
&\le Ch_K^{-1/2}(\|v\|_{L^2(K)}+h_K^r|v|_{H^r(K)}),
\end{aligned}
\end{equation*}
where the final inequality follows from~(\ref{B_Khbehaviour}). Furthermore, the estimate is independent of $h_K$, and the choice of $K$. Thus, we have obtained the desired estimate.$\quad\quad\square$
\begin{lemma}\label{inverselemma19:00}
Assume that $(\Th)_h$ is a family of meshes on $\overline{\X}$ that is regular of order $m\in\mathbb N$. For any $v\in\dg$, the following inverse estimate holds for any $K\in\Th$, with $0\le s\le m$, and $q\in[2,\infty]$:
\begin{equation}\label{local:invineq}
|v|_{W_*^{m,q}(K)}\le C_Ih_K^{s-m}|v|_{W_*^{s,q}(K)},
\end{equation}
where the positive constant $C_I$ is independent of $K$ and $h_K$.
\end{lemma}
\emph{Proof:} We first note that ~(\ref{local:invineq}) is trivial when $m=0$, since then $s=0$, and $|\cdot|_{W_*^{m,q}}=|\cdot|_{W_*^{s,q}}=\|\cdot\|_{L^{q}}$, so we will assume that $m\ge1$. We will first prove~(\ref{local:invineq}) when $s=0$. By~(\ref{nonaff:scaling}), for $j\in\mathbb N$,  $1\le j\le m$, $q\in\{2,\infty\}$, and any $K\in\Th$, we have 
\begin{equation}\label{inv1stbound}
|v|_{W^{j,q}(K)}\le C|\det\tilde{B}_K|^{1/q}\|\tilde{B}_K^{-1}\|^j\left(\sum_{r=\ms{j}}^j\|\tilde{B}_K\|^{2(j-r)}|\hat v|_{W^{r,q}(\hat K)}\right).
\end{equation}
Now, for $0\le r\le j$, 
\begin{equation*}
\|\tilde{B}_K\|^{2(j-r)}|\hat v|_{W^{r,q}(\hat K)}\le C(\sigma) h_K^{2(j-r)}|\hat v|_{W^{r,q}(\hat K)},
\end{equation*}
where the inequality is due to~(\ref{B_Khbehaviour}). Now, let $K^*$ be a \emph{fixed} reference element, and take $\tilde{G}_{\hat K}:K^*\to\hat K$, with $\tilde{G}_{\hat K}(x^*)=\tilde{A}_{\hat K}x^*+\tilde{a}_{\hat K}$, with $\tilde{A}_{\hat K}\in GL(\mathbb R^d)$ and $\tilde{a}_{\hat K}\in\mathbb R^d$. As in the proof of Lemma~\ref{fractracelem}, it follows that $\tilde{A}_{\hat K}$ belongs to a compact subset $BL$ of $GL(\mathbb R^d)$.

Now, defining $v^*(x^*) = \hat v(\tilde{G}_{\hat K}(x^*))$, it follows that $v^*\in\mathbb P^p(K^*)$, where $\mathbb P^p(K^*)$ is of finite dimension, depending only on $K^*$, $d$ and $p$, thus by the equivalence of norms on finite dimensional spaces, we see that
\begin{equation}\label{equivofnormsused19:08}
\begin{aligned}
|\hat v|_{W^{r,q}(\hat K)} & \le\|\tilde{A}_{\hat K}\|^r|\det\tilde{A}_{\hat K}|^{1/q}|v^*|_{W^{r,q}(K^*)}\\
&\le\|\tilde{A}_{\hat K}\|^r|\det\tilde{A}_{\hat K}|^{1/q}\|v^*\|_{W^{r,q}(K^*)}\\
&\le C(d,p,K^*)\|\tilde{A}_{\hat K}\|^r|\det\tilde{A}_{\hat K}|^{1/q}\|v^*\|_{L^q(K^*)}\\
&\le C(d,p,K^*)\|\tilde{A}_{\hat K}\|^r|\det\tilde{A}_{\hat K}|^{1/q}|\det\tilde{A}_{\hat K}^{-1}|^{1/q}\|\hat v\|_{L^q(\hat K)}\\
&= C(d,p,K^*)\|\tilde{A}_{\hat K}\|^r\|\hat v\|_{L^q(\hat K)}\\
&\le C(d,p,K^*)(\max_{B\in BL}\|\tilde{A}_{\hat K}\|^r)\|\hat v\|_{L^q(\hat K)}\\
&\le C(d,p,\sigma,K^*,r)\|\hat v\|_{L^q(\hat K)}.
\end{aligned}
\end{equation}
Thus, applying the above inequality,~(\ref{nonaffine:scaling}) with $l=0$, and~(\ref{B_Khbehaviour}), to~(\ref{inv1stbound}), we obtain
\begin{equation*}
\begin{aligned}
|v|_{W^{j,q}(K)}&\le C(d,p,\sigma,K^*)|\det\tilde{B}_K|^{1/q}\|\tilde{B}_K^{-1}\|^j\|\hat v\|_{L^{q}(\hat K)}\\
&\le C(d,p,\sigma,K^*)|\det\tilde{B}_K|^{1/q}\|\tilde{B}_K^{-1}\|^j|\det\tilde{B}_K|^{-1/q}\|v\|_{L^{q}(K)}\\
&\le C(d,p,\sigma,K^*,j) h_K^{-j}\|v\|_{L^{q}(K)}\\
&\le C(d,p,\sigma,K^*,m) h_K^{-j}\|v\|_{L^{q}(K)}.
\end{aligned}
\end{equation*}
Since our choice of $1\le j\le m$ was arbitrary, we may take $1\le k\le m$, and sum the above over $1\le j\le k$, obtaining
\begin{equation}\label{star:17:14}
|v|_{W_*^{k,q}(K)}\le C(d,p,\sigma,K^*,m)h_K^{-k}\|v\|_{L^{q}(K)}\quad 1\le k\le m.
\end{equation}
We obtain~(\ref{local:invineq}) with $s=0$, by setting $k=m$ above.
We will now prove~(\ref{local:invineq}) for $1\le s\le m$. 

In this case we will argue by induction, and as our base case, we shall prove the result for $s=1$. Take $1\le j\le m$, and let $|\alpha|=j$. Then we may write $D^\alpha v=D^\beta(D^\gamma v)$ for some $|\beta|=j-1$, $|\gamma|=1$. One must note that by the chain rule, $Dv|_K = D(\hat v\circ F_K^{-1})|_K=(D\hat v\circ F_K^{-1})DF_K^{-1}$, where the components of $(D\hat v\circ F_K^{-1})DF_K^{-1}$ do not necessarily belong to $\mathbb P^p(\hat K)$. It is the case, however, that $D^{\delta}\hat v\in\mathbb P^p(\hat K)$ for any $|\delta|=1$.

One can see that
\begin{equation}\label{0:18:03}
\begin{aligned}
\|D^\alpha v\|_{L^q(K)}&\le|D^\gamma v|_{W_*^{j-1,q}(K)}\\
&\le|Dv|_{W_*^{j-1,q}(K)}\\
&=\left|(D\hat v\circ F_K^{-1})DF_K^{-1}\right|_{W_*^{j-1,q}(K)}\\
&\lesssim\sum_{r=\ms{j-1}}^{j-1}\sup_{x\in K}\|D^r(DF_K^{-1}(x))\||D\hat v\circ F_K^{-1}|_{W_*^{j-1-r,q}(K)}\\
&\lesssim\max_{\min\{2,j\}\le r\le j}\sup_{x\in K}\|D^r F_K^{-1}(x)\||D\hat v\circ F_K^{-1}|_{W_*^{j-1,q}(K)}.
\end{aligned}
\end{equation}
By~(\ref{indinvproof2}) and~(\ref{indinvproof4}), we have that
\begin{equation}\label{1:18:03}
\max_{\min\{2,j\}\le r\le j}\sup_{x\in K}\|D^r F_K^{-1}(x)\|\le\max_{\min\{2,j\}\le r\le j}c_{-r}\|\tilde{B}_K\|^{2(r-1)}\|\tilde{B}_K^{-1}\|^r,\end{equation}
where we are denoting $c_{-1}:=1/(1-C_K)$. Furthermore, since $D \hat v\in[\mathbb P^{p-1}(\hat K)]^d\subset$\\$[\mathbb P^{p}(\hat K)]^d$, we can apply~(\ref{star:17:14}) with $k=j-1$, obtaining
\begin{equation}\label{2:18:03}
|D\hat v\circ F_K^{-1}|_{W_*^{j-1,q}(K)}\lesssim h_K^{1-j}\|D\hat v\circ F_K^{-1}\|_{L^{q}(K)}.
\end{equation}
We also have that
\begin{equation}\label{3:18:03}
\begin{aligned}
\|D\hat v\circ F_K^{-1}\|_{L^{q}(K)}&=\|(D\hat v\circ F_K^{-1}DF_K^{-1})(DF_K^{-1})^{-1}\|_{L^{q}(K)}\\
&\le\sup_{\hat x\in \hat K}\|DF_K\|\|Dv\|_{L^q(K)}\\
&=\sup_{\hat x\in \hat K}\|DF_K\||v|_{W_*^{1,q}(K)}.
\end{aligned}
\end{equation}
Applying~(\ref{1:18:03}),~(\ref{2:18:03}), and~(\ref{3:18:03}) to~(\ref{0:18:03}), and summing over all $|\alpha|=j$, we obtain\\
$$|v|_{W^{j,q}(K)}\lesssim\max_{\min\{2,j\}\le r\le j}c_{-r}\|\tilde{B}_K\|^{2(r-1)}\|\tilde{B}_K^{-1}\|^r\sup_{\hat x\in \hat K}\|DF_K\|h_K^{1-j}|v|_{W_*^{1,q}(K)}.$$
Lastly, applying~(\ref{indinvproof1}) and~(\ref{B_Khbehaviour}) to the above estimate, we obtain (noting that $\Th$ is regular of order $m$)
\begin{equation}\label{10:11}|v|_{W^{j,q}(K)}\lesssim\max_{\min\{2,j\}\le r\le j}h_K^{r-1}h_K^{1-j}|v|_{W_*^{1,q}(K)}\le h_K^{1-j}|v|_{W_*^{1,q}(K)}.
\end{equation}
Again, our choice of $1\le j\le m$ was arbitrary, and so we can sum~(\ref{10:11}) over $1\le j\le k$ for any $1\le k\le m$, obtaining
$$|v|_{W_*^{k,q}(K)}\lesssim\max_{\min\{2,k\}\le r\le k}h_K^{r-1}h_K^{1-k}|v|_{W_*^{1,q}(K)}\le h_K^{1-k}|v|_{W_*^{1,q}(K)}\quad1\le k\le m.$$
To proceed to argue by induction, we will assume that for $1\le s\le k\le m-1$ we have
\begin{equation}\label{indhyp:11:46}
|v|_{W_*^{k,q}(K)}\lesssim h_K^{s-k}|v|_{W_*^{s,q}(K)},
\end{equation}
and we will use this to show that
$$|v|_{W_*^{k,q}(K)}\lesssim h_K^{s+1-k}|v|_{W_*^{s+1,q}(K)},$$
for $1\le s+1\le k\le m$.

To this end, let us take $s+1\le j\le m$ and let $|\alpha|=j$. Again we write $D^\alpha v = D^\beta(D^\gamma v)$ for some $|\beta|=j-1$, and $|\gamma|=1$, and so, analogous to our previous argument, we obtain
\begin{equation*}
\begin{aligned}
\|D^\alpha v\|_{L^q(K)}\le |D^\gamma v|_{W_*^{j-1,q}(K)}&\le|Dv|_{W_*^{j-1,q}(K)}\\
&\lesssim h_K^{-1}|D\hat v\circ F_K^{-1}|_{W_*^{j-1,q}(K)}.
\end{aligned}
\end{equation*}
Applying our inductive hypothesis~(\ref{indhyp:11:46}) with $k=j-1\ge s$, we obtain
\begin{equation}\label{12:00}
\|D^\alpha v\|_{L^q(K)}\lesssim h_{K}^{-1}h_K^{s-(j-1)}|D\hat v\circ F_K^{-1}|_{W_*^{s,q}(K)}.
\end{equation}
Now, 
\begin{equation*}
\begin{aligned}
|D\hat v\circ F_K^{-1}|_{W_*^{s,q}(K)} & = |(D\hat v\circ F_K^{-1}DF_K^{-1})DF_K\circ F_K^{-1}|_{W_*^{s,q}(K)}\\
&=|(Dv)DF_K\circ F_K^{-1}|_{W_*^{s,q}(K)}\\
&\le\sum_{r=\ms{s}}^s\sup_{x\in K}\|D^r(DF_K\circ F_K^{-1})(x)\||Dv|_{W_*^{s-r,q}(K)}\\
&\le\sum_{r=\ms{s}}^s\sup_{x\in K}\|D^r(DF_K\circ F_K^{-1})(x)\||Dv|_{W_*^{s,q}(K)}\\
&\le\sum_{r=\ms{s}}^s\sup_{x\in K}\|D^r(DF_K\circ F_K^{-1})(x)\||v|_{W_*^{s+1,q}(K)}.
\end{aligned}
\end{equation*}
Applying the above to~(\ref{12:00}), we obtain
$$\|D^\alpha v\|_{L^q(K)}\lesssim\left(h_K^{-1}\sum_{r=\ms{s}}^s\sup_{x\in K}\|D^r(DF_K\circ F_K^{-1})(x)\|\right)h_K^{s+1-j}|v|_{W_*^{s+1,q}(K)}.$$
Let us momentarily assume that, for any $1\le s\le m-1$,
\begin{equation}\label{it:remains}
h_K^{-1}\sum_{r=\ms{s}}^s\sup_{x\in K}\|D^r(DF_K\circ F_K^{-1})(x)\|\lesssim 1.
\end{equation}
Then we obtain
$$\|D^\alpha v\|_{L^q(K)}\lesssim h_K^{s+1-j}|v|_{W_*^{s+1,q}(K)},$$
where $|\alpha|=j$, and $s+1\le j\le m$ was arbitrary. Summing over all $|\alpha|=j$, and then all $s+1\le j\le k\le m$, we obtain
$$\sum_{s+1\le|\alpha|\le k}\|D^\alpha v\|_{L^q(K)}\lesssim h_K^{s+1-k}|v|_{W_*^{s+1,q}(K)}\quad s+1\le k\le m.$$
It is also clear that
$$\sum_{\ms{s}\le|\alpha|\le s}\|D^\alpha v\|_{L^q(K)}\lesssim |v|_{W_*^{s,q}(K)}\le |v|_{W_*^{s+1,q}(K)}\lesssim h_K^{s+1-k}|v|_{W_*^{s+1,q}(K)},$$
and so we obtain
$$|v|_{W_*^{k,q}(K)}\lesssim \sum_{\ms{s}\le|\alpha|\le s}\|D^\alpha v\|_{L^q(K)}+\sum_{s+1\le|\alpha|\le k}\|D^\alpha v\|_{L^q(K)}\lesssim h_K^{s+1-k}|v|_{W_*^{s+1,q}(K)},$$
$s+1\le k\le m$, which concludes our inductive argument, and yields~(\ref{local:invineq}) for $1\le s\le m$, by taking $k=m$. It remains to show that~(\ref{it:remains}) is in fact true. Let us recall the formula~\cite{MR1014883}
 \begin{equation*}
 D^r(f\circ g) = \sum_{i=1}^r(D^if\circ g)\left(\sum_{\alpha\in E(r,i)}c_\alpha\prod_{l=1}^r(D^lg)^{\alpha_l}\right),
 \end{equation*}
where $E(r,i)$ is the set given by:
\begin{equation}\label{1Eri:def}
E(r,i):=\left\{\alpha\in\mathbb N_0^r:|\alpha|=i\quad\mbox{and}\quad\sum_{l=1}^rl\alpha_l=r\right\},
\end{equation}
and the $c_\alpha$'s, $\alpha\in E(m,r)$ are some given constants, bounded independently of $h_K$. From this, we obtain
\begin{equation*}
\begin{aligned}
&h_K^{-1}\sum_{r=\ms{s}}^s\sup_{x\in K}\|D^r(DF_K\circ F_K^{-1})(x)\|\\
&~~~=h_K^{-1}\sum_{r=\ms{s}}^s\sup_{x\in K}\left\|\sum_{i=1}^r(D^{i+1}F_K)\circ F_K^{-1}(x)\sum_{\alpha\in E(r,i)}c_\alpha\prod_{l=1}^r(D^lF_K^{-1})^{\alpha_l}(x)\right\|\\
&~~~\lesssim h_K^{-1}\sum_{r=\ms{s}}^s\sum_{i=1}^r\sup_{x\in K}\|(D^{i+1}F_K)\circ F_K^{-1}(x)\|\sum_{\alpha\in E(r,i)}\prod_{l=1}^r\sup_{x\in K}\|(D^lF_K^{-1})(x)\|^{\alpha_l}\\
&~~~\lesssim h_K^{-1}\sum_{r=\ms{s}}^s\sum_{i=1}^rc_{i+1}\|\tilde{B}_K\|^{i+1}\sum_{\alpha\in E(r,i)}\prod_{l=1}^r\|\tilde{B}_K\|^{2(l-1)\alpha_l}\|\tilde{B}_K^{-1}\|^{l\alpha_l},
\end{aligned}
\end{equation*}
where the final inequality follows from~(\ref{indinvproof4}), and the fact that the mesh is regular of order $m\ge s+1$. Applying~(\ref{B_Khbehaviour}), and noting that by definition, if $\alpha\in E(r,i)$, then $|\alpha| = i$ and $\sum_{l=1}^rl\alpha_l = r$, we obtain
\begin{equation*}
\begin{aligned}
h_K^{-1}\sum_{r=\ms{s}}^s&\sup_{x\in K}\|D^r(DF_K\circ F_K^{-1})(x)\|\lesssim~~~~~~~~~~~~~~~~~~~~~~~~~\\
\end{aligned}
\end{equation*}
\begin{equation*}
\begin{aligned}
&~~~~~~~~~~~~~~\lesssim h_K^{-1}\sum_{r=\ms{s}}^s\sum_{i=1}^rh_K^{i+1}\sum_{\alpha\in E(r,i)}h_K^{\sum_{l=1}^r(2l\alpha_l-2\alpha_l)}h_K^{\sum_{l=1}^r\alpha_l}\\
&~~~~~~~~~~~~~~=h_K^{-1}\sum_{r=\ms{s}}^s\sum_{i=1}^rh_K^{i+1}\sum_{\alpha\in E(r,i)}h_K^{2(r-i)}h_K^r\\
&~~~~~~~~~~~~~~\lesssim h_K^{-1}\sum_{r=\ms{s}}^sh_K^r\sum_{i=1}^rh_K^{1-i}\\
&~~~~~~~~~~~~~~= \sum_{r=\ms{s}}^sh_K^r\sum_{i=1}^rh_K^{-i}\lesssim\sum_{r=0}^sh_K^rh_K^{-r}\lesssim 1,
\end{aligned}
\end{equation*}
as desired.
Note that the estimates we have derived are \emph{independent} of the choice of $K\in\Th$.$\quad\quad\square$
\subsection{Interpolation estimates}
The proofs of the following lemmas can both be found in~\cite{MR1014883}, i.e., Theorem 4.1 and Corollary 4.1; one must note that they are both given in a more general context. However, we are considering Lagrange finite element spaces, which satisfy the hypotheses of Theorem 4.1 and Corollary 4.1 (see examples 1 and 2 on page 1221 of~\cite{MR1014883}).
\begin{lemma}[Optimal local interpolation in $\dg$]\label{thm:star:16:27}
Assume that the family\\$(\Th)_h$ is regular of order $m$. Let $\ell,s,p\in\mathbb N_0$, $p\ge2$, with $\ell\le s\le\min\{p,m\}+1$. Then for any $K\in\Th$, and any $u\in H^s(\X;\Th)$, there exists a $z_h\in\dg$ such that
\begin{equation}\label{local:opt:iterp}
\|u-z_h\|_{H^\ell(K)}\le Ch_K^{s-\ell}|u|_{H_*^s(K)},
\end{equation}
where the constant $C$ is independent of $h_K$, $u$, and $K$.
\end{lemma}
\begin{remark}
We note that the classical Lagrange interpolation operator can only be applied functions that have well defined point values. Even in two dimensions, it is not in general the case that functions in $H^1$ have well defined point values. This leads one to define other interpolation operators that require less regularity, in particular, we define a local interpolation operator that is well defined on $L^2$ functions (one of the first examples is due to P. Clem\'{e}nt~\cite{MR0400739}, using local averaging; however the one we will define is provided in~\cite{MR1014883} and is slightly different).
\end{remark}
\begin{definition}[Local $L^2$ projection]
For $v\in L^2(\X)$, and $K\in\Th$, we define $\hat\rho_v$ to be the unique element of $\mathbb P^p(\hat K)$ that satisfies
\begin{equation}\label{localL2proj}
\int_{\hat K}(\hat v-\hat\rho_v)\hat\rho\quad\forall\hat\rho\in\mathbb P^p(\hat K).
\end{equation}
\end{definition}
\begin{definition}[Local Lagrange interpolation operator]
For $K\in\Th$, we define the Lagrange interpolation operator $\Pi_h:L^2(K)\to P_K$, where $(K,P_K,\Sigma_K)$ is a curved Lagrange finite element of type $(m,p)$, by
\begin{equation}\label{nonclassinterp}
\Pi_h(v)=\sum_{\mu\in\Sigma_K}\mu(\rho_v)\rho_\mu,
\end{equation}
where $\rho_v=\hat\rho_v\circ F_K^{-1}$, with $\hat\rho_v$ satisfying~(\ref{localL2proj}), and  $\{\rho_\mu\}_{\mu\in\Sigma_K}$ forms a basis of $P_K$.
\end{definition}
%Note that the following definition is given in the context of a general finite element space, that could be given by any of the above spaces. As such, we shall denote this space by $V_h$. When using the interpolation operator, it will be clear which finite element space we are considering.
%\begin{definition}[Global $L^2$ projection]\label{projoper}
%For $v\in[L^2(\X)]^{m\times n}$, $m,n\in\mathbb N$, given a finite element space $V_h\subset[L^2(\X)]^{m\times n}$, we define the global $V_h$-projection operator $\proj_{V_h}:[L^2(\X)]^{m\times n}\to V_h$, by
%\begin{equation}\label{proj:oper}
%\int_\X\proj_{V_h}(v)\!:\!\Phi = \int_\X v\!:\!\Phi\quad\forall\Phi\in V_h.
%\end{equation}
%\end{definition}
\begin{lemma}[$H^r$-multipliers]
Assume that $u\in H^r(\X)$, $0<r<1$ and $\psi\in C^{0,1}(\overline{\X})$. Then, there exists a constant $C$ depending only on $d$ and $r$, such that
\begin{equation}\label{www15:48}
|u\,\psi|_{H^r(\X)}\le\sqrt{2}\|\psi\|_{L^\infty(\X)}|u|_{H^r(\X)}+
\sqrt{2}C(d,r)\sqrt{1+\diam(\X)^2}|\psi|_{C^{0,1}(\overline{\X})}\|u\|_{L^2(\X)}.
\end{equation}
\end{lemma}
\emph{Proof:} We see that
\begin{equation}\label{0:15L45}
\begin{aligned}
&|u\,\psi|_{H^r(\X)}^2 = \int_\X\int_\X\frac{|u(x_1)\psi(x_1)-u(x_2)\psi(x_2)|^2}{|x_1-x_2|^{d+2r}}\\
& ~~~~\le 2\int_\X\int_{\X}\frac{|u(x_1)\psi(x_1)-u(x_2)\psi(x_1)|^2}{|x_1-x_2|^{d+2r}}+2\int_\X\int_{\X}\frac{|u(x_2)\psi(x_1)-u(x_2)\psi(x_2)|^2}{|x_1-x_2|^{d+2r}}.
\end{aligned}
\end{equation}
It then follows that
\begin{equation}\label{1:15L45}
\begin{aligned}
\int_\X\int_{\X}\frac{|u(x_1)\psi(x_1)-u(x_2)\psi(x_1)|^2}{|x_1-x_2|^{d+2r}}&\le\|\psi\|_{L^\infty(\X)}^2\int_\X\int_{\X}\frac{|u(x_1)-u(x_2)|^2}{|x_1-x_2|^{d+2r}}\\
&=\|\psi\|_{L^\infty(\X)}^2|u|_{H^r(\X)}^2.
\end{aligned}
\end{equation}
Furthermore, we note that for $x_1,x_2\in\X\cap\{|x_1-x_2|\ge1\}$, $|x_1-x_2|\le\operatorname{diam}(\X)$, and so
\begin{equation}\label{14:14:14}
\begin{aligned}
\int_\X\int_{\X}\frac{|u(x_2)\psi(x_1)-u(x_2)\psi(x_2)|^2}{|x_1-x_2|^{d+2r}}&\le|\psi|_{C^{0,1}(\overline{\X})}^2\int_\X\int_\X\frac{|u(x_2)|^2|x_1-x_2|^2}{|x_1-x_2|^{d+2r}}\\
&=|\psi|_{C^{0,1}(\overline{\X})}^2\int_\X\int_{\X\cap\{|x_1-x_2|\le1\}}\frac{|u(x_2)|^2}{|x_1-x_2|^{d+2(r-1)}}\\
&~~~~~+|\psi|_{C^{0,1}(\overline{\X})}^2\operatorname{diam}(\X)^2\int_\X\int_{\X\cap\{|x_1-x_2|\ge1\}}\frac{|u(x_2)|^2|x_1-x_2|^2}{|x_1-x_2|^{d+2r}}\\
& = |\psi||_{C^{0,1}(\overline{\X})}^2\left(\|\tilde u^2*\tilde{R}_1\|_{L^1(\mathbb R^d)}+\operatorname{diam}(\X)^2\|\tilde u^2*\tilde{R}_2\|_{L^1(\mathbb R^d)}\right),
\end{aligned}
\end{equation}
where $*$ denotes convolution over $\mathbb R^d$, and the functions $\tilde u,\tilde R_1,\tilde R_2\in L^2(\mathbb R^d)$ are defined as follows:
\begin{equation*}
\begin{aligned}
\tilde u|_\X &= u,\quad&\tilde u|_{\mathbb R^d\setminus\X} = 0,\\
\tilde R_1|_{\X\cap\{|x|\le 1\}} & = |x|^{-(d+2(r-1))},\quad&\tilde R_1|_{\mathbb R^d\setminus(\X\cap\{|x|\le 1\})} = 0,\\
\tilde R_2|_{\X\cap\{|x|\ge 1\}} & = |x|^{-(d+2r)},\quad&\tilde R_2|_{\mathbb R^d\setminus(\X\cap\{|x|\ge 1\})} = 0.
\end{aligned}
\end{equation*}
Then, applying Young's inequality for convolutions, we find that 
\begin{align*}
\|\tilde u^2*\tilde{R}_1\|_{L^1(\mathbb R^d)}&\le \|\tilde u^2\|_{L^1(\mathbb R^d)}\|\tilde{R}_1\|_{L^1(\mathbb R^d)}\\
& \le \|\tilde{u}\|_{L^2(\mathbb R^d)}^2\int_{\mathbb R^d\setminus(\X\cap\{|x|\le 1\})}|x|^{-(d+2(r-1))}\\
& \le \|u\|_{L^2(\X)}^2C(d)\int_{0}^1z^{-(d+2(r-1))}z^{d-1}\\
& = \|u\|_{L^2(\X)}^2C(d)\int_{0}^1z^{1-2r}\\
& = \|u\|_{L^2(\X)}^2C(d)\left.\frac{z^{2(1-r)}}{2(1-r)}\right|^1_0\\
& = \|u\|_{L^2(\X)}^2C(d)/(2(1-r)) = C(d,r)\|u\|_{L^2(\X)}^2.
\end{align*}
Similarly, we obtain
\begin{align*}
\|\tilde u^2*\tilde{R}_2\|_{L^1(\mathbb R^d)}&\le C(d)\|u\|_{L^2(\X)}^2\int_1^\infty z^{-(1+2r)}\\
& = \frac{C(d)}{2r}\|u\|_{L^2(\X)}^2 = C(d,r)\|u\|_{L^2(\X)}^2.
\end{align*}
Applying these two estimates to~(\ref{14:14:14}), we obtain
\begin{equation}\label{2:15L45}
\begin{aligned}
\int_\X\int_{\X}\frac{|u(x_2)\psi(x_1)-u(x_2)\psi(x_2)|^2}{|x_1-x_2|^{d+2r}}\le C(d,r)(1+\diam(\X)^2)|\psi|_{C^{0,1}(\overline{\X})}^2\|u\|_{L^2(\X)}^2.
\end{aligned}
\end{equation}
Applying~(\ref{1:15L45}) and~(\ref{2:15L45}) to~(\ref{0:15L45}), we obtain
$$|u\,\psi|_{H^r(\X)}^2\le2\|\psi\|_{L^\infty(\X)}^2|u|_{H^r(\X)}^2+2C(d,r)(1+\diam(\X)^2)|\psi|_{C^{0,1}(\overline{\X})}\|u\|_{L^2(\X)}^2.$$
Taking square roots above, we obtain~(\ref{www15:48}).$\quad\quad\square$
\begin{lemma}[Integer and non integer regularity interpolation estimates]
Assume that $\X$ is piecewise $C^{m+1}$, with $m\in\mathbb N$, $m\ge k+1$, $k\in\{1,2\}$. Let $\{\Th\}_h$ be a family of triangulations on $\overline{\X}$ that is regular of order $m$, satisfying Assumption~\ref{Meshconds}. Let $u\in H^{\mathbf{s}_k}(\X;\Th)$, $\mathbf{s}_k = (s_K^k)_{K\in\Th}$, with $s_K^k>2k-1/2$ for all $K\in\Th$. Then, there exists a $z_{k,h}\in\dg$, $p\ge 2k-1$, and a constant $C$, independent of $u_k$,  and $h_K$, but dependent on $\max_Ks^k_K$, such that for each $K\in\Th$, each nonnegative integer $q\le 2k-1$, and each multi-index $\beta$ with $|\beta|=q$, we have
\begin{equation}\label{opt:interp:2}
\begin{aligned}
\|u-z_{k,h}\|_{H^q(K)}&\le Ch_K^{t_K^k-q}\|u_k\|_{H^{s_K^k}(K)},\\
\|D^\beta(u-z_{k,h})\|_{L^2(\partial K)}&\le Ch_K^{t_K-q-1/2}\|u\|_{H^{s^k_K}(K)},
\end{aligned}
\end{equation}
where $t^k_K=\min\{p+1,m+1,s^k_K\}$. Furthermore, under the domain and mesh hypotheses for $k=1$, if $u\in H^{\mathbf{s}_1}(\X;\Th)$, $\mathbf{s}_1 = (s_K^1)_{K\in\Th}$, with $s_K^1>5/2$, for all $K\in\Th$, then~(\ref{opt:interp:2}) holds for $q\le2\le p$.
\end{lemma}
\emph{Proof:} We will first discuss how we will obtain the second bound of~(\ref{opt:interp:2}). Let $k\in\{1,2\}$. We either have that $q\le 2k-1$ is a nonnegative integer, and $t_K^k-q>2k-1/2-q\ge1/2$, or we have that $q\le2\le p$, and $t_K^k-q = \min\{p+1,m+1,s_K^1\}-q>5/2-2=1/2$. Thus, under the hypotheses of the lemma, for $k=1,2$ we have that $t_K^k-q>1/2$.

 Since the family of triangulations is regular of order $m$, it follows that for any $\beta$ such that $|\beta|=q$, and any $v\in\dg$, that $D^\beta(u-v)\in H^{t_K^k-q}(K)$. In particular, $t_K^k-q>1/2$. Thus, we may apply the trace estimate~(\ref{tracestfrac}) with $r_K=t_K^k-q>1/2$, obtaining
$$\|D^\beta(u-v)\|_{L^2(\partial K)}\le Ch_K^{-1/2}(\|D^\beta(u-v)\|_{L^2(K)}+h_K^{r_K}|D^\beta(u-v)|_{H^{r_K}(K)}).$$
Let us assume that there exists a $z_h\in\dg$ satisfying the first estimate of~(\ref{opt:interp:2}). Then, setting $v=z_h$ above we obtain
\begin{equation}\label{want"10:59}
\begin{aligned}
\|D^\beta(u-z_h)\|_{L^2(\partial K)}&\le Ch_K^{-1/2}(\|D^\beta(u-z_h)\|_{L^2(K)}+h_K^{r_K}|D^\beta(u-z_h)|_{H^{r_K}(K)})\\
&\le Ch_K^{-1/2}(h_K^{t_K-q}\|u\|_{H^{s_K}(K)}+h_K^{r_K}|D^\beta(u-z_h)|_{H^{r_K}(K)}).
& 
\end{aligned}
\end{equation}
Thus, to obtain both estimates of~(\ref{opt:interp:2}), it suffices to prove that the there exists a $z_h\in\dg$ such that the first estimate of~(\ref{opt:interp:2}) holds, as well as the following:
\begin{equation}\label{otherestotherway}
|u-z_h|_{H^{q+r_K}(K)}=|u-z_h|_{H^{t_K}(K)}\le Ch_K^{t_K-q-r_K}\|u\|_{H^{s_K}(K)}.
\end{equation}
Since, applying the above estimate to~(\ref{want"10:59}), and noting the factor $h_K^{r_K}$ in the second inequality of~(\ref{want"10:59}), we obtain the second estimate of~(\ref{opt:interp:2}). Note that we already have such bounds in the case that $s_K$ is an integer, and as such, we shall assume from this point on that $s_K\notin\mathbb N$.

We will now prove the first estimate of~(\ref{opt:interp:2}). Let $\beta$ satisfy $|\beta|=q$, and let $K^*$ be a fixed reference simplex. Then, from~(\ref{nonaff:scaling}) we obtain
\begin{equation}\label{11759}
\begin{aligned}
|D^\beta(u-z_h)|_{L^2(K)}&\le|u-z_h|_{H^q(K)}\\
&\le C|\det\tilde B_K|^{1/2}\|\tilde{B}_K^{-1}\|^q\sum_{j=\min\{1,q\}}^q\|\tilde B_K\|^{2(q-j)}|\hat u-\hat z_h|_{H^j(\hat K)}\\
&\le C(K^*,\sigma)|\det\tilde B_K|^{1/2}\|\tilde{B}_K^{-1}\|^q\sum_{j=\min\{1,q\}}^q\|\tilde B_K\|^{2(q-j)}|u^*-z_h^*|_{H^j(K^*)}.\\
\end{aligned}
\end{equation}
We take the function $z_h\in\dg$, defined as follows: $z_h|_K = \Pi_hu|_K$ where $\Pi_h$ is the local interpolation operator, given by~(\ref{nonclassinterp}). Due to~(\ref{localL2proj}), this operator reproduces polynomials in $\mathbb P^p(\hat K)$, and so we may apply Theorem 5 of~\cite{MR0336957} in conjunction with Theorem 1.8 of~\cite{MR3136501} (applying Theorem 1.8 of~\cite{MR3136501} allows us to consider noninteger Sobolev spaces when applying the Bramble--Hilbert Lemma), obtaining for $\min\{1,q\}\le j\le q$
\begin{equation}\label{21759}
\begin{aligned}
|u^*-z_h^*|_{H^j(K^*)} & \le \|u^*-z_h^*\|_{H^{t^k_K}(K^*)}\\
&\le C(K^*)|u^*|_{H^{t^k_K}(K^*)}\\
&\le C(K^*,\sigma)|\hat u|_{H^{t^k_K}(\hat K)},
\end{aligned}
\end{equation}
%\begin{equation}\label{21759}
%\begin{aligned}
%|u^*-z_h^*|_{H^j(K^*)} & \le \|u^*-z_h^*\|_{H^{j+r}(K^*)}\\
%&\le C(K^*)|u^*|_{H^{j+r}(K^*)}\\
%&\le C(K^*,\sigma)|\hat u|_{H^{j+r}(\hat K)},
%\end{aligned}
%\end{equation}
where, by assumption $t_K^k>q$ (note that the final inequality follows from a scaling argument similar to the one used in estimate~(\ref{16:33:wed}), noting that $K^*$ and $\hat K$ are affine equivalent, and the mesh is shape regular). We now decompose $t_K=\ell_K+r_K$, where $\ell_K\ge k$ is an integer, and $r_K\in(0,1)$. We see that
$$|\hat u|^2_{H^{t_K}(\hat K)}=|\hat u|^2_{H^{\ell_K+r_K}(\hat K)} = \int_{\hat K}\int_{\hat K}\frac{|D^{\ell_K}\hat u(\hat x_1)-D^{\ell_K}\hat u(\hat x_2)|^2}{|\hat x_1-\hat x_2|^{d+2r_K}}=|D^{\ell_K}\hat u|^2_{H^{r_K}(\hat K)}.$$
Let us recall the formula
 \begin{equation}\label{1compdiff:rule}
 D^r(f\circ g) = \sum_{i=1}^r(D^if\circ g)\left(\sum_{\alpha\in E(r,i)}c_\alpha\prod_{l=1}^r(D^lg)^{\alpha_l}\right),
 \end{equation}
where $E(r,i)$ is the set
\begin{equation}\label{1Eri:def}
E(r,i):=\left\{\alpha\in\mathbb N_0^r:|\alpha|=i\quad\mbox{and}\quad\sum_{l=1}^rl\alpha_l=r\right\},
\end{equation}
and the $c_\alpha$'s, $\alpha\in E(m,r)$ are some given constants, bounded independently of $h_K$. By the triangle inequality, we obtain
$$|D^{\ell_K}\hat u|_{H^{r_K}(\hat K)}\le\sum_{i=1}^{\ell_K}\left|(D^iu\circ F_K)\sum_{\alpha\in E(\ell_K,i)}c_\alpha\prod_{l=1}^{\ell_K}(D^lF_K)^{\alpha_l}\right|_{H^{r_K}(\hat K)}.$$
We now apply~(\ref{www15:48}) to the above estimate, obtaining (noting that $\hat K$ is contained in the unit ball, and thus $\diam(\hat K)\le 2$)
\begin{equation}\label{19:54;;}
\begin{aligned}
|D^{\ell_K}\hat u|_{H^{r_K}(\hat K)}&\le C(d,r)\sum_{i=1}^{\ell_K}\left\|\sum_{\alpha\in E(\ell_K,i)}c_\alpha\prod_{l=1}^{\ell_K}(D^lF_K)^{\alpha_l}\right\|_{L^\infty(\hat K)}|D^iu\circ F_K|_{H^{r_K}(\hat K)}\\
&~~~~+C(d,r)\sum_{i=1}^{\ell_K}\left|\sum_{\alpha\in E(\ell_K,i)}c_\alpha\prod_{l=1}^{\ell_K}(D^lF_K)^{\alpha_l}\right|_{C^{0,1}(\overline{\hat K})}\|D^iu\circ F_K\|_{L^2(\hat K)}.
\end{aligned}
\end{equation}
By~(\ref{indinvproof1}), and the fact that the triangulation is regular of order $m\ge k+1$ (and that $\mathbb N\ni\ell_K<t_K\le m+1$, so $\ell_K\le m$), we estimate the first term on the right-hand side of~(\ref{19:54;;}) as follows
\begin{equation}\label{firstterm2120}
\begin{aligned}
&C(d,r)\sum_{i=1}^{\ell_K}\left\|\sum_{\alpha\in E(\ell_K,i)}c_\alpha\prod_{l=1}^{\ell_K}(D^lF_K)^{\alpha_l}\right\|_{L^\infty(\hat K)}|D^iu\circ F_K|_{H^{r_K}(\hat K)}\\
&~~~~~~\le C\sum_{i=1}^{\ell_K}\sum_{\alpha\in E(\ell_K,i)}c_\alpha\prod_{l=1}^{\ell_K}\|\tilde{B}_K\|^{l\alpha_l}|D^iu\circ F_K|_{H^{r_K}(\hat K)}\\
&~~~~~~\le C\sum_{i=1}^{\ell_K}\|\tilde{B}_K\|^{\ell_K}|D^iu\circ F_K|_{H^{r_K}(\hat K)}.
\end{aligned}
\end{equation}
For the second term, we see that
\begin{equation}\label{second:term:20:06}
\begin{aligned}
&C(d,r)\sum_{i=1}^{\ell_K}\left|\sum_{\alpha\in E(\ell_K,i)}c_\alpha\prod_{l=1}^{\ell_K}(D^lF_K)^{\alpha_l}\right|_{C^{0,1}(\overline{\hat K})}\|D^iu\circ F_K\|_{L^2(\hat K)}\\
&~~~~~\le C\sum_{i=1}^{\ell_K}\left\|\sum_{\alpha\in E(\ell_K,i)}c_\alpha D\left(\prod_{l=1}^{\ell_K}(D^lF_K)^{\alpha_l}\right)\right\|_{L^\infty(\hat K)}\|D^iu\circ F_K\|_{L^2(\hat K)}\\
&~~~~~\le C\sum_{i=1}^{\ell_K}\sum_{\alpha\in E(\ell_K,i)}c_\alpha\left\|D\left(\prod_{l=1}^{\ell_K}(D^lF_K)^{\alpha_l}\right)\right\|_{L^\infty(\hat K)}\|D^iu\circ F_K\|_{L^2(\hat K)}\\
\end{aligned}
\end{equation}
Furthermore, for $\alpha\in E(\ell_K,i)$
\begin{equation*}
\begin{aligned}
&\left\|D\left(\prod_{l=1}^{\ell_K}(D^lF_K)^{\alpha_l}\right)\right\|_{L^\infty(\hat K)} = \left\|\sum_{l=1}^{\ell_K}D((D^lF_K)^{\alpha_l})\prod_{j=1,j\ne l}^{\ell_K}(D^jF_K)^{\alpha_j}\right\|_{L^\infty(\hat K)}\\
&~~~\le C(d,\ell_K) \sum_{l=1}^{\ell_K}\alpha_l\sup_{\hat x\in\hat K}\|D^{l+1}F_K(\hat x)\|\sup_{\hat x\in\hat K}\|D^{\alpha_l}F_K(\hat x)\|^{\alpha_l-1}\prod_{j=1,j\ne l}^{\ell_K}\sup_{\hat x\in\hat K}\|D^jF_K(\hat x)\|^{\alpha_j},
\end{aligned}
\end{equation*}
since the triangulation is regular of order $m\ge k+1$ (and $\ell_K+1\le m+1$), applying~(\ref{indinvproof1}) to the above yields
\begin{equation}\label{21:13}
\begin{aligned}
\left\|D\left(\prod_{l=1}^{\ell_K}(D^lF_K)^{\alpha_l}\right)\right\|_{L^\infty(\hat K)}&\le C\sum_{l=1}^{\ell_K}\alpha_l\|\tilde{B}_K\|^{l+1+l(\alpha_l-1)}\prod_{j=1,j\ne l}^{\ell_K}\|\tilde{B}_K\|^{j\alpha_j}\\
&\le C\sum_{l=1}^{\ell_K}\alpha_l\|\tilde{B}_K\|^{l+1+l(\alpha_l-1)}\|\tilde{B}_K\|^{\sum_{j=1}^{\ell_K}j\alpha_j-l\alpha_l}\\
&= C\sum_{l=1}^{\ell_K}\alpha_l\|\tilde{B}_K\|^{l+1+l(\alpha_l-1)-\ell_K-l\alpha_l}\\
& = C\|\tilde B_K\|^{1+\ell_K}\sum_{l=1}^{\ell_K}\alpha_l=Ci\|\tilde B_K\|^{1+\ell_K}.
\end{aligned}
\end{equation}
Applying~(\ref{21:13}) to~(\ref{second:term:20:06}) gives us
\begin{equation}\label{theiradsas}
\begin{aligned}
&C(d,r)\sum_{i=1}^{\ell_K}\left|\sum_{\alpha\in E(\ell_K,i)}c_\alpha\prod_{l=1}^{\ell_K}(D^lF_K)^{\alpha_l}\right|_{C^{0,1}(\overline{\hat K})}\|D^iu\circ F_K\|_{L^2(\hat K)}\\
&~~~~~\le C\sum_{i=1}^{\ell_K}i\sum_{\alpha\in E(\ell_K,i)}c_\alpha \|\tilde{B}_K\|^{1+\ell_K}\|D^iu\circ F_K\|_{L^2(\hat K)}\\
&~~~~~\le C\|\tilde{B}_K\|^{1+\ell_K}\sum_{i=1}^{\ell_K}\|D^iu\circ F_K\|_{L^2(\hat K)}.
\end{aligned}
\end{equation}
We now apply~(\ref{firstterm2120}) and~(\ref{theiradsas}) to~(\ref{19:54;;}), obtaining
\begin{equation}\label{§7:51}
|D^{\ell_K}\hat u|_{H^{r_K}(\hat K)} \le C\|\tilde{B}_K\|^{\ell_K}\sum_{i=1}^{\ell_K}|D^iu\circ F_K|_{H^{r_K}(\hat K)}+\|\tilde{B}_K\|\|D^iu\circ F_K\|_{L^2(\hat K)}.
\end{equation}
Applying the change of variables formula in the $L^2$-norms in~(\ref{§7:51}), and the scaling argument~(\ref{scalingarg1})--(\ref{scalingarg2}) to the $|D^iu\circ F_K|_{H^{r_K}(\hat K)}$ term for $i=\ell_K$ (noting that this argument is valid for any $r_K\in(0,1)$, as long as the function has $H^{r_K}$-regularity) in~(\ref{§7:51}), in conjunction with~(\ref{B_Khbehaviour}), we obtain
\begin{equation*}
\begin{aligned}
|D^{\ell_K}\hat u|_{H^{r_K}(\hat K)} &\le Ch_K^{\ell_K-\frac{d}{2}}(h_K^{r_K}|D^{\ell_K}u|_{H^{r_K}(K)}+\sum_{i=1}^{\ell_K}h_K|u|_{H^i(K)})+Ch_K^{\ell_K}\sum_{i=1}^{\ell_K-1}|D^{i}u\circ F_K|_{H^{r_K}(\hat K)}\\
&\le Ch_K^{\ell_K+r_K-d/2}\|u\|_{H^{\ell_K+r_K}(K)}+Ch_K^{\ell_K}\sum_{i=1}^{\ell_K-1}|D^{i}u\circ F_K|_{H^{r_K}(\hat K)}\\
\end{aligned}
\end{equation*}
\begin{equation}\label{31759}
\begin{aligned}
&\le Ch_K^{t_K-d/2}\|u\|_{H^{t_K}(K)}+Ch_K^{\ell_K}\sum_{i=1}^{\ell_K-1}|D^{i}u\circ F_K|_{H^{r_K}(\hat K)}~~~~~~~~~~~~~~~\\
&\le Ch_K^{t_K-d/2}\|u\|_{H^{s_K}(K)}+Ch_K^{\ell_K}\sum_{i=1}^{\ell_K-1}|D^{i}u\circ F_K|_{H^{r_K}(\hat K)},
\end{aligned}
\end{equation}
where the constant $C$ is independent of $h_K$ and the choice of $K\in\Th$ (note that we have utilised the continuous embedding $H^{s_K}(K)\subseteq H^{t_K}(K)$, where the constant in the embedding only depends upon $d$ and $r_K$, due to Proposition 2.1 of~\cite{MR2944369}).
We note, however, that the terms of the sum on the right-hand side of the final inequality of~(\ref{31759}) are not present in the $H^{t_K}$-norm.
Furthermore, for $1\le i<\ell_K$, we note the following:
\begin{equation}\label{fixed:14:39}
\begin{aligned}
|D^iu\circ F_K|_{H^{r_K}(\hat K)}&\le C(\sigma,K^*)|(D^iu)^*|_{H^{r_K}(K^*)}\\
&=C(\sigma,K^*)|(D^iu)^*-M|_{H^{r_K}(K^*)},\\
\end{aligned}
\end{equation}
for any $M\in[\mathbb P^0(K^*)]^{\operatorname{dim}(D^iu)}$, 
where the first inequality follows from a scaling argument, and the fact that the mesh is regular, and the final equality holds due to the fact that constant functions are in the kernel of $|\cdot|_{H^r}$. 
%We now apply an interpolation (between $L^2$ and $H^1$) estimate given by Corollary 5.1 of~\cite{BREZIS2017}, yeilding
We now use the fact that the embedding $H^1(K^*)\subseteq H^{r_K}(K^*)$ is continuous, obtaining
\begin{equation}\label{really?fixed:14:39}
\begin{aligned}
|D^iu\circ F_K|_{H^{r_K}(\hat K)}&\le C(K^*,\sigma,d,r_K)\inf_{M\in[\mathbb P^0(K^*)]^{\operatorname{dim}(D^iu)}}\|(D^iu)^*-M\|_{H^{r_K}(K^*)}\\
&\le C(K^*,\sigma,d,r_K)|(D^iu)^*|_{H^1(K^*)}\\
&\le C(K^*,\sigma,d,r_K)|D^iu\circ F_K|_{H^1(\hat K)},
\end{aligned}
\end{equation}
where the penultimate inequality follows from an application of Theorem 1.8 of~\cite{MR3136501}, and the final inequality follows from the fact that the mesh is regular.

Thus, we obtain
$$\sum_{i=1}^{\ell_K-1}|D^iu\circ F_K|_{H^{r_K}(\hat K)}\le C\sum_{i=1}^{\ell_K-1}|D^iu\circ F_K|_{H^1(\hat K)}\le Ch_K^{-d/2+1}\sum_{i=1}^{\ell_K-1}|D^iu|_{H^1(K)}.$$
Applying the above to~(\ref{31759}) gives us
\begin{equation}\label{§7:52}
|D^{\ell_K}\hat u|_{H^{r_K}(\hat K)} \le Ch_K^{t_K-\frac{d}{2}}\|u\|_{H^{s_K}(K)}+Ch_K^{\ell_K+1-\frac{d}{2}}\sum_{i=1}^{\ell_K-1}|D^iu|_{H^1(K)}\le Ch_K^{t_K-\frac{d}{2}}\|u\|_{H^{s_K}(K)}.
\end{equation}
Finally, applying~(\ref{§7:52}),~(\ref{21759}), and~(\ref{B_Khbehaviour}) to~(\ref{11759}), we obtain
$$|D^\beta(u-z_h)|_{L^2(K)}\le Ch_K^{-q}\sum_{j=\min\{1,q\}}^qh_K^{2(q-j)}h_K^{t_K}\|u\|_{H^{s_K}(K)}\le Ch_K^{t_K-q}\|u\|_{H^{s_K}(K)},$$
which is the first estimate of~(\ref{want"10:59}). Estimate~(\ref{otherestotherway}) is obtained in a similar manner, utilising~(\ref{indinvproof4}). $\quad\quad\square$

\subsection{Discrete Poincar\'{e}--Friedrichs' inequalities}
\begin{lemma}[Discrete Poincar\'{e}--Friedrichs' inequality]\label{disc:PF}
Assume that $\{\Th\}_h$ is regular of order $2$ family of triangulations, and let $v\in\dg$. Then, the following inequality holds
\begin{equation}\label{est:disc:PF}
\|v\|_{L^2(\X)}^2\le C\left(|v|_{H^1(\X;\Th)}^2+\sum_{F\in\Eb}\|u_h\|_{L^2(F)}^2+\sum_{F\in\Ei}\tilde{h}_F^{-1}\|\jump{u_h}\|_{L^2(F)}^2\right),
\end{equation}
where the positive constant, $C$, depends only on the shape-regularity constants of the mesh, $d$, and $\X$.
\end{lemma}
\emph{Proof:} Let $K\in\Th$, and take $v\in\dg$. We see that
\begin{equation*}
\begin{aligned}
\int_K|v|^2 &= \frac{1}{d}\int_K\nabla\cdot(xv^2)-2v\sum_{i=1}^dx_iD_{i}v\\
&\le\frac{1}{d}\left(\int_{\partial K}(xv^2)\cdot n_{\partial K}+\int_K\frac{d}{2}|v|^2+2\sum_{i=1}^dx_i^2|D_{i}v|^2\right),
\end{aligned}
\end{equation*}
subtracting $(1/2)\int_Kv^2$ from each side and multiplying by $2$ yields
$$\int_K|v|^2\le\frac{2}{d}\left(\int_{\partial K}(xv^2)\cdot n_{\partial K}+2\sum_{i=1}^dx_i^2|D_{i}v|^2\right).$$
Summing the above over all $K\in\Th$, and denoting $n_F$ to be a \emph{fixed} choice of unit normal to $F\in\Eib$, we obtain
\begin{equation}\label{star:20:47:81}
\begin{aligned}
\sum_{K\in\Th}\|v\|_{L^2(K)}^2&\le\frac{2}{d}\left(\sum_{F\in\Eib}\int_F\jump{xv^2}\cdot n_F+\sum_{K\in\Th}2\sum_{i=1}^d\int_Kx_i^2|D_iv|^2\right)\\
&\le\frac{2}{d}\left(\sum_{F\in\Ei}\int_F\jump{xv^2}\cdot n_F+C(\X)\sum_{F\in\Eb}\|v\|_{L^2(F)}^2+2C(\X)^2|v|_{H^1(\X;\Th)}^2\right),
\end{aligned}
\end{equation}
where $C(\X):=\max_{x\in\overline{\X}}\max_{i=1,\ldots,d}|x_i|$. Furthermore, we have that
\begin{equation*}
\begin{aligned}
\sum_{F\in\Ei}\int_F\jump{xv^2}\cdot n_F & = \sum_{F\in\Ei}\int_F(\jump{x}\avg{v^2}+\avg{x}\jump{v^2})\cdot n_F\\
& = \sum_{F\in\Ei}\int_F(2x\avg{v}\jump{v})\cdot n_F\\
& \le 2C(\X)\sum_{F\in\Ei}\|\avg{v}\|_{L^2(F)}\|\jump{v}\|_{L^2(F)}\\
& \le \sum_{F\in\Ei}C(\X)^2(\delta\tilde{h}_F)^{-1}\|\jump{v}\|_{L^2(F)}^2+\delta\tilde{h}_F\|\avg{v}\|_{L^2(F)}^2,
\end{aligned}
\end{equation*}
for any $\delta>0$. We then apply the trace inequality~(\ref{1simptrace}), obtaining
\begin{equation*}
\begin{aligned}
\sum_{F\in\Ei}\int_F\jump{xv^2}\cdot n_F & \le C(\X)^2\sum_{F\in\Ei}(\delta\tilde{h}_F)^{-1}\|\jump{v}\|_{L^2(F)}^2\\
&~~~~~+\delta C(d)\sum_{K\in\Th:F\subset\partial K}\tilde{h}_F(h_K^{-1}\|v\|_{L^2(K)}^2+h_K\|\nabla v\|_{L^2(K)}^2)\\
&\le C(\X)^2\sum_{F\in\Ei}(\delta\tilde{h}_F)^{-1}\|\jump{v}\|_{L^2(F)}^2\\
&~~~~~+\delta C(d)\sum_{K\in\Th}\|v\|_{L^2(K)}^2+h_K\tilde{h}_F\|\nabla v\|_{L^2(K)}^2.
\end{aligned}
\end{equation*}
Applying the above estimate to~(\ref{star:20:47:81}), we obtain, for any $\delta>0$,
\begin{equation*}
\begin{aligned}
\sum_{K\in\Th}\|v\|_{L^2(K)}^2&\le\frac{2}{d}\left(2C(\X)^2|v|_{H^1(\X;\Th)}^2+\sum_{F\in\Ei}C(\X)^2(\delta\tilde{h}_F)^{-1}\|\jump{v}\|_{L^2(F)}^2\right.\\
\end{aligned}
\end{equation*}
\begin{equation}\label{21:35:35:21}
\begin{aligned}
&\left.~~~~~~~~~~~~+\delta C(d)\sum_{K\in\Th}\|v\|_{L^2(K)}^2+\|\nabla v\|_{L^2(K)}^2+C(\X)\sum_{F\in\Eb}\|v\|_{L^2(F)}^2\right).
\end{aligned}
\end{equation}
Choosing $\delta$ sufficiently small, so that $2\delta C(d)/d\le1/2$, subtracting $(1/2)\|v\|_{L^2(\X)}^2$ from each side of~(\ref{21:35:35:21}) and multiplying by $2$ we obtain the desired estimate.$\quad\quad\square$
\begin{lemma}[Gradient Poincar\'{e}--Friedrichs' inequality]\label{gradpoin16:54}
Assume that $\{\Th\}_h$ is regular of order $2$ family of triangulations, and let $v\in\dg$. Then, the following inequality holds
\begin{equation}\label{GPF:est:disc}
|v|_{H^1(\X;\Th)}^2\le C_P\left(|v|_{H^2(\X;\Th)}^2+\sum_{F\in\Ei}\tilde{h}_F^{-1}\|\jump{\nabla v\cdot n_F}\|_{L^2(F)}^2+\sum_{F\in\Eib}\tilde{h}_F^{-1}\|\jump{v}\|_{L^2(F)}^2\right),
\end{equation}
where the positive constant, $C_P$, depends only on the shape regularity constants of the mesh, $d$, and $\X$.
\end{lemma}
\emph{Proof:} Let $v\in\dg$, and take any $K\in\Th$. An application of the divergence theorem gives us
$$\int_K|\nabla v|^2 = -\int_K(\Delta v)v+\int_{\partial K}(\nabla v\cdot n_F)v.$$
Summing this equality over all $K\in\Th$, gives us
\begin{equation}\label{poin:13:06}
\begin{aligned}
\sum_{K\in\Th}|v|_{H^1(K)}^2 & = -\sum_{K\in\Th}\langle\Delta v,v\rangle_K+\sum_{F\in\Eib}\langle\jump{\nabla v\cdot n_F},\avg{v}\rangle_F+\sum_{F\in\Ei}\langle\avg{\nabla v\cdot n_F},\jump{v}\rangle_F\\
&\le\sum_{K\in\Th}\frac{\delta}{2}\|v\|_{L^2(K)}^2+\frac{1}{2\delta}C(d)|v|_{H^2(K)}^2+\sum_{F\in\Eb}\frac{\delta\tilde{h}_F}{2}\|\nabla v\|_{L^2(F)}^2+\frac{1}{2\delta\tilde{h}_F}\|v\|_{L^2(F)}^2\\
&~~~~~+\sum_{F\in\Ei}\frac{1}{2\delta\tilde{h}_F}\|\jump{\nabla v\cdot n_F}\|_{L^2(F)}^2+\frac{\delta\tilde{h}_F}{2}\|\avg{v}\|_{L^2(F)}^2\\
&~~~~~+\sum_{F\in\Ei}\frac{1}{2\delta\tilde{h}_F}\|\jump{v}\|_{L^2(F)}^2+\frac{\delta\tilde{h}_F}{2}\|\avg{\nabla v\cdot n_F}\|_{L^2(F)}^2,
\end{aligned}
\end{equation}
for any $\delta>0$. Applying the trace estimate~(\ref{1simptrace}), we obtain (noting that $\tilde{h}_F\le h_K$)
\begin{equation*}
\begin{aligned}
&\sum_{F\in\Ei}\frac{\delta\tilde{h}_F}{2}(\|\avg{v}\|_{L^2(F)}^2+\|\avg{\nabla v\cdot n_F}\|_{L^2(F)}^2)+\sum_{F\in\Eb}\frac{\delta\tilde{h}_F}{2}\|\nabla v\|_{L^2(F)}^2\\
&\le C\frac{\delta}{2}\sum_{F\in\Eib}\sum_{K\in\Th:F\subset\partial K}\tilde{h}_F(\|v\|_{L^2(\partial K)}^2+\|\nabla v\|_{L^2(\partial K)}^2)\\
&\le C\frac{\delta}{2}\sum_{F\in\Eib}\sum_{K\in\Th:F\subset\partial K}\tilde{h}_F(h_K^{-1}\|v\|_{L^2(K)}^2+(h_K+h_K^{-1})\|\nabla v\|_{L^2(K)}^2+h_K^{-1}|\nabla v|_{H^1(K)}^2)\\
&\le C\frac{\delta}{2}\|v\|_{H^2(\X;\Th)}^2;
\end{aligned}
\end{equation*}
applying this to~(\ref{poin:13:06}) gives us
\begin{equation*}
\begin{aligned}
|v|_{H^1(\X;\Th)}^2&\le C\delta\|v\|_{L^2(\X)}^2+C\delta|v|_{H^1(\X;\Th)}^2+C(\delta^{-1}+\delta)|v|_{H^2(\X;\Th)}^2\\
&~~~~+\delta^{-1}\sum_{F\in\Eb}\tilde{h}_F^{-1}\|v\|_{L^2(F)}^2+\sum_{F\in\Ei}\tilde{h}_F^{-1}(\|\jump{\nabla v\cdot n_F}\|_{L^2(F)}^2+\|\jump{v}\|_{L^2(F)}^2).
\end{aligned}
\end{equation*}
We now apply~(\ref{est:disc:PF}) to the estimate above, which yields (noting that $1\le\tilde{h}_F^{-1}$)
\begin{equation*}
\begin{aligned}
|v|_{H^1(\X;\Th)}^2&\le 2C\delta|v|_{H^1(\X;\Th)}^2+C(\delta^{-1}+\delta)|v|_{H^2(\X;\Th)}^2\\
&~~~~+\delta^{-1}\sum_{F\in\Eb}\tilde{h}_F^{-1}\|v\|_{L^2(F)}^2+\sum_{F\in\Ei}\tilde{h}_F^{-1}(\|\jump{\nabla v\cdot n_F}\|_{L^2(F)}^2+\|\jump{v}\|_{L^2(F)}^2).
\end{aligned}
\end{equation*}
We now choose $\delta$ sufficiently small, so that $2C\delta\le1/2$, which gives us
\begin{equation*}
\begin{aligned}
|v|_{H^1(\X;\Th)}^2&\le\frac{1}{2}|v|_{H^1(\X;\Th)}^2+C_P\left[|v|_{H^2(\X;\Th)}^2\!+\!\sum_{F\in\Eib}\tilde{h}_F^{-1}\|\jump{v}\|_{L^2(F)}^2+\sum_{F\in\Ei}\tilde{h}_F^{-1}\|\jump{\nabla v\!\cdot\!n_F}\|_{L^2(F)}^2\right].
\end{aligned}
\end{equation*}
Subtracting $(1/2)|v|_{H^1(\X;\Th)}^2$ from both sides of the inequality and multiplying by $2$ yields the desired estimate.$\quad\quad\square$
\subsection{Tangential operators and curved simplex curvature bounds}
In order to appropriately define and bound the bilinear forms that define our method, we need to be able to define tangential differential operators (i.e., operators that involve derivatives that are tangential to the faces of the curved simplices $K$ of the mesh), and bound the curvature terms arising in the bilinear form (these curvature terms appear both on boundary faces, and on interior faces if the dimension $d\ge3$).

\emph{Tangential differential operators.} For $F\in\mathscr{E}^{i,b}$, denote for $s>1/2$ the space of $H^s$-regular tangential vector fields on $F$ by $H^s_\Ta(F):=\{v\in H^s(F)^d:v\cdot n_F=0~\mbox{on}~F\}$. Below we define the tangential gradient 
$\nabla_\Ta:H^s(F)\to H^{s-1}_T(F)$ and the tangential divergence $\operatorname{div}_\Ta:H^s_\Ta(F)\to H^{s-1}(F)$, where $1\le s\le 2$ (note that in the case that $\partial\X$ is piecewise $C^{m}$, with $m\ge 2$, we are able to consider $1\le s\le m$).

We see that $F\subset\partial K$, for some $K\in\Th$. Since $K$ is piecewise $C^2$ (see the proof of Lemma~\ref{tange:oper:bounds}), for a.e. $x\in\partial K$, there exists a neighbourhood $W_x$ of $x$ in $\partial K$, sufficiently small to allow the existence of a family of $C^2$ curves that satisfy the following: a curve of each family passes through every point of $W_x$, and the unit tangent vectors to these curves form an orthonormal system (assumed to be oriented with respect to $\overline{n}$, where $\overline{n}$ is the unit outward normal to $\partial K$) at every point of $W_x$. We take the lengths $s_1,\ldots,s_{d-1}$ along each of these curves, respectively, to be the local coordinate system, and denote $t_1,\ldots,t_{d-1}$ to be the unit tangent vectors along each curve, respectively. In this notation, we have the following for $v:\partial K\to\mathbb R^d$:
$$v = v_\Ta+(v\cdot\overline{n})\overline{n},\quad v_\Ta:=\sum_{j=1}^{d-1}(v\cdot t_j)t_j.$$
For $\phi\in C^1(\overline{K})$, and $\psi\in C^1(\overline{K})^d$, with $\psi|_{\partial K} = \sum_{j=1}^{d-1}\psi_jt_j$, we obtain 
\begin{equation}\label{tangeopers}
\nabla\phi|_{\partial K} = \nabla_\Ta\phi+\frac{\partial\phi}{\partial\overline{n}}\overline{n},\quad\nabla_\Ta\phi = \sum_{j=1}^{d-1}\frac{\partial\phi}{\partial s_j}t_j,
\end{equation}
and
\begin{equation}\label{tangeopers2}
\operatorname{div}_\Ta\psi =\nabla_\Ta\cdot\psi= \sum_{j=1}^{d-1}\frac{\partial\psi_j}{\partial s_j},
\end{equation}
which extend to $\phi\in H^s(K)$, $s>3/2$, by density and the construction of the trace operator. Furthermore, one can see that by rearranging the first identity of~(\ref{tangeopers}), that $\nabla_\Ta = \nabla-\overline{n}\frac{\partial}{\partial\overline{n}}$ (and thus $\operatorname{div}_\Ta$) is well defined $a.e.$ on $\partial K$, and is independent of the choice of normal $\overline{n}$.

We approach~(\ref{tangeopers}) and~(\ref{tangeopers2}) in the context of traces and Sobolev spaces, in the following lemma. In particular we are able to decompose the Laplacian, $\Delta$, in terms of the tangential Laplacian $\Delta_\Ta :=\operatorname{div}_\Ta\nabla_\Ta$, the mean curvature of the face, and first and second order normal derivatives.
\begin{lemma}
Let $\X$ be a piecewise $C^{2}$ domain, and let $\{\mathcal{T}_h\}_{h>0}$ be a family of meshes on $\overline{\X}$ that is regular of order $1$ and satisfies Assumption~\ref{Meshconds}.  Then, for any $h>0$, for each $K\in\mathscr{T}_h$ {} and each face $F\subset\partial K$, the following identities hold:
\begin{align}
\tau_F(\nabla v)  &= \nabla_\Ta(\tau_Fv)+\left(\tau_F\frac{\partial v}{\partial n_F}\right)n_F\quad\forall v\in H^s(K),\,s>3/2,\label{identitiesgrad}\\
\tau_F(\Delta v)  &= \operatorname{div}_\Ta\nabla_\Ta(\tau_Fv)+\mc_F\left(\tau_F\frac{\partial v}{\partial n_F}\right)+\tau_F\frac{\partial}{\partial n_F}(\nabla v\cdot n_F),\label{identitieslap}
\end{align}
for all $v\in H^s(K)$, $s>5/2,$ where $n_F$ is a fixed choice of unit normal to $F$, $\mc_F:=\nabla_\Ta\cdot n_F$ is the mean curvature of the face $F$, and $\tau_F$ is the trace operator from $K$ to $F$. 
%Let $s_1,\ldots,s_{d-1}$ be an orthonormal coordinate system on $F$, with tangent v. Then, for $u\in H^s(F)$ and $v=\sum_{i=1}^{d-1}v_it_i$, with $v_i\in H^s(F)$ for $i=1,\ldots,d-1$, we define
%\begin{equation}\label{tangeopers}
%\nabla_\Ta u:=\sum_{i=1}^{d-1}t_i\frac{\partial u}{\partial t_i},\quad\operatorname{div}_\Ta v:=\sum_{i=1}^{d-1}\frac{\partial v}{\partial t_i}.
%\end{equation}
\end{lemma}

\emph{Proof:} Let us take $U\in C^3(\overline{K})$, and for $F\in\Eib$, let $u = U|_F$. Then, as the family of meshes $\{\mathcal{T}_h\}_{h>0}$ is regular of order $1$, it follows that $F\subset\partial K$ for some $K\in\Th$, where $K$ is piecewise $C^2$ (see the proof of Theorem~\ref{tange:oper:bounds}). Thus, we may extend (without relabelling) the unit normal to $F$, $n_F$ (note that this choice of unit normal is fixed, and that~(\ref{identitiesgrad}) and~(\ref{identitieslap}) are independent of this choice), by $n_F\in C^1(\overline{K})$ (note that the extension may not be normal to the other faces of $\partial K$, when restricted there), and so also define an extension of the tangential gradient, $\nabla_\Ta: C^3(\overline{K})\to C^1(\overline{K})^d$, by
$$\nabla_\Ta U = \nabla U -\frac{\partial U}{\partial n_F}n_F.$$
This can be rearranged to yield
$$\nabla U = \nabla_\Ta U + \frac{\partial U}{\partial n_F}n_F.$$
Upon restricting to $F$, we obtain
\begin{equation*}
\begin{aligned}
\nabla U|_F & = \nabla_\Ta U|_F + \left.\left(\left(\frac{\partial U}{\partial n_F}\right)n_F\right)\right|_F\\
& = \nabla_\Ta|_F(U|_F) + \left.\left(\frac{\partial U}{\partial n_F}\right)\right|_Fn_F|_F\\
& = \nabla_\Ta(U|_F) + \left.\left(\frac{\partial U}{\partial n_F}\right)\right|_Fn_F.
\end{aligned}
\end{equation*}
Thus, by density and the construction of the trace operator, this extends to $u\in H^s(K)$, $s>3/2$, giving us
\begin{equation}\label{identities15:03}
\tau_F(\nabla u) = \nabla_\Ta(\tau_Fu)+\left(\tau_F\frac{\partial u}{\partial n_F}\right)n_F,
\end{equation}
which is~(\ref{identitiesgrad}).

For the identity~(\ref{identitieslap}), we follow a similar approach to~\cite{MR645791}, in which the statement is essentially proven for $d=2,3$.
Now, for $x\in F$ let us take a local coordinate system $s_1,\ldots,s_{d-1}$, on a neighbourhood $W_x$ of $x$ in $F$. Expressing $F$ locally as the graph of a $C^2$ function $\phi$, we see that
$$u(s_1,\ldots,s_{d-1}) = U(s_1,\ldots,s_{d-1},\phi(s_1,\ldots,s_{d-1})).$$
Furthermore, let us assume that the coordinates have been chosen so that $\nabla_{s'}\phi(0) = 0$ (denoting $s'=(s_1,\ldots,s_{d-1})$), so that the local coordinates $\{s',s_d\}=\{s',\phi(s')\}$ are tangent to the hyperplane $\{s_d=0\}$ at $x = (0,\phi(0))$. Then, in $W_x$, we have that
$$\operatorname{div}_\Ta\nabla_\Ta u = \sum_{j=1}^{d-1}\frac{\partial^2 u}{\partial s_j^2},$$
where, for $j=1,\ldots,d-1$,
\begin{equation*}
\begin{aligned}
\frac{\partial^2 u}{\partial s_j^2} & = \frac{\partial}{\partial s_j}\left(\frac{\partial}{\partial s_j}(U(s',\phi(s')))\right)\\
& = \frac{\partial}{\partial s_j}\left(U_j(s',\phi(s'))+\frac{\partial\phi}{\partial s_j}U_d(s',\phi(s'))\right)\\
& = U_{jj}(s',\phi(s'))+2\frac{\partial\phi}{\partial s_j}U_{dj}(s',\phi(s'))+\frac{\partial^2\phi}{\partial s_j^2}U_d(s',\phi(s'))+\left(\frac{\partial\phi}{\partial s_j}\right)^2U_{dd}(s',\phi(s')),
\end{aligned}
\end{equation*}
where $U_{j}$, $U_{jk}$ denote the first and second order partial derivatives in the $j$ and $j,k$ components of $U$, respectively.
Thus, at $x$, i.e., at $s'=0$, we have
\begin{equation*}
\operatorname{div}_\Ta\nabla_\Ta u = \sum_{j-1}^{d-1}U_{jj}(0,\phi(0))+U_d(0,\phi(0))\sum_{j=1}^{d-1}\frac{\partial^2\phi(0)}{\partial s_j^2}.
\end{equation*}
Moreover, at $x$, $U_{dd} = \frac{\partial^2 U}{\partial n_F^2}$, $U_d=\frac{\partial U}{\partial n_F}$, and $\sum_{j=1}^{d-1}\frac{\partial^2\phi}{\partial s_j^2} = -\mc_F$. Thus, at $x$,
$$\Delta U = \operatorname{div}_\Ta\nabla_\Ta u  + \mc_F\frac{\partial U}{\partial n_F}+\frac{\partial^2 U}{\partial n_F^2}.$$
This decomposition is valid at any $x\in F$, and so we obtain
\begin{equation*}
\begin{aligned}
\Delta U|_F & = \operatorname{div}_\Ta\nabla_\Ta u+\left.\mc_F\frac{\partial U}{\partial n_F}\right|_F+\left.\frac{\partial^2 U}{\partial n_F^2}\right|_F\\
&=\operatorname{div}_\Ta\nabla_\Ta (U|_F)+\left.\left.\mc_F\frac{\partial U}{\partial n_F}\right|_F+\frac{\partial^2 U}{\partial n_F^2}\right|_F.
\end{aligned}
\end{equation*}
Thus, by density, applying~(\ref{identities15:03}), for $u\in H^s(K)$, $s>5/2$, we obtain
$$\tau_F(\Delta u) = \operatorname{div}_\Ta\nabla_\Ta(\tau_Fu)+\mc_F\left(\tau_F\frac{\partial u}{\partial n_F}\right)+\tau_F\frac{\partial}{\partial n_F}(\nabla u\cdot n_F),$$
which is~(\ref{identitieslap}).$\quad\quad\square$
\begin{lemma}\label{tange:oper:bounds}
Let $\X$ be a piecewise $C^{2}$ domain, and let $\{\mathcal{T}_h\}_{h>0}$ be a family of meshes on $\overline{\X}$ that is regular of order $1$, and satisfies Assumption~\ref{Meshconds}.  %Then, for any $h>0$, for each $K\in\mathscr{T}_h$ {} and each face $F\subset\partial K$, the following identities hold:
%\begin{equation}\label{identities}
%\begin{aligned}
%\tau_F(\nabla v)  &= \nabla_\Ta(\tau_Fv)+\left(\tau_F\frac{\partial v}{\partial n_F}\right)n_F\quad\forall v\in H^s(K),\,s>3/2,\\
%\tau_F(\Delta v)  &= \operatorname{div}_\Ta\nabla_\Ta(\tau_Fv)+\mc_F\left(\tau_F\frac{\partial v}{\partial n_F}\right)+\tau_F\frac{\partial}{\partial n_F}(\nabla v\cdot n_F),\quad\forall v\in H^s(K),
%\end{aligned}
%\end{equation}
%$s>5/2,$ where $\mc_F$ is the mean curvature of the face $F$, and $\tau_F$ is the trace operator from $K$ to $F$. 
Then, there exists a constant $C$, depending on $\X,d$ and the family of triangulations $\{\mathcal{T}_h\}_{h>0}$, such that for $\Eib\ni F\subset\partial K$, on ${F}$ we have that
\begin{equation}\label{bounded2ndff}
(\nabla_\Ta v)^T\nabla_\Ta n_F^T\nabla_\Ta w\le C|\nabla_\Ta v||\nabla_\Ta w|\quad\forall v,w\in H^s(K),\,\,s>3/2.
\end{equation}
%If, in addition, $\X$ is a piecewise $C^2$ piecewise convex domain, then, for $\Eb\ni F\subset\partial K$, on $F$ we have that
%\begin{equation}\label{coercive2ndff}
%\mc_F\ge 0.
%\end{equation}

%\begin{equation}\label{coercive2ndff}
%(\nabla_\Ta v)^T\nabla_\Ta n_F^T\nabla_\Ta v\ge 0\quad\forall v\in H^s(K),\,\,s>3/2.
%\end{equation}
{}
\end{lemma}
%\emph{Proof:} First, if $F\in\Ei$, then $F$ is flat, and both identities in~(\ref{identities}) follow from Lemma 4 in~\cite{MR3077903}. 
%
%If $F\in\Eb$, then the identities follow similarly to the proof of Lemma 4 in~\cite{MR3077903}, taking into account the fact that for $K\in\Th$ such that $\Eb\ni F\subset\partial K$, the Laplacian of a smooth function $u\in C^\infty(\overline{K})$ can be decomposed as follows
%$$\Delta u|_F =\operatorname{div}_\Ta\nabla_\Ta u+\left.\mc_F\frac{\partial u}{\partial n}\right|_F+\left.\frac{\partial^2u}{\partial n^2}\right|_F.$$
%Noting that the trace operator, $\tau$, commutes with tangential partial derivatives, by the density of smooth functions in $H^s(K)$, we obtain
%$$\tau_F(\Delta v)  = \operatorname{div}_\Ta\nabla_T(\tau_Fv)+\mc_F\left(\tau_F\frac{\partial v}{\partial n_F}\right)+\tau_F\frac{\partial}{\partial n_F}(\nabla v\cdot n_F).\quad\quad\square$$
%{}
{}
\emph{Proof:} First, let us assume that $F\in\Eb$.
Then, since $\X$ is piecewise $C^2$, $F\subset\Gamma_n\subset\partial\X$, where $\Gamma_n$ is a $C^2$ portion of $\partial\X$. It then follows that for a given $F\in\Eb$, $n_F$  is of class $C^{1}(F)$.
%(note however that its restriction to $\partial K$ may only remain normal to $F$). 
For any two vector-valued functions $\xi^1,\xi^2:F\to\mathbb R^d$ tangent to $F$, it then follows that
$$(\xi^1)^T\nabla_\Ta n_F^T\,\xi^2\le\sup_{x\in F}|\nabla_\Ta n_F^T(x)||\xi_1||\xi_2|\le\sup_{x\in\Gamma_n}|\nabla_\Ta n_{\Gamma_n}^T(x)||\xi_1||\xi_2|.$$
%In order to obtain a uniform bound (i.e., independent of the choice of $F$), we must first distinguish between the cases where $F\in\Eb$ and $F\in\Ei$. If $F\in\Eb$, then  $F\subset\Gamma_n\subset\partial\X$, where $\Gamma_n$ is a $C^2$ portion of $\partial\X$, and thus the extension of $n_F$ can be chosen to coincide with a $C^1$ extension of the normal to $\Gamma_n$, and thus, $\sup_{\overline{K}}|\nabla_\Ta n_F^T|\le\max_{F\subset\Gamma_n}\sup_{\overline{K}}|\nabla_\Ta n_{\Gamma_n}^T|=:C(n)$.
Thus, for an arbitrary $F\in\Eb$,
$$(\xi^1)^T\nabla_\Ta n_F^T\,\xi^2\le\max_{i=1,\ldots,N}\sup_{x\in\Gamma_n}|\nabla_\Ta n_{\Gamma_n}^T(x)||\xi_1||\xi_2|=C(\X)|\xi_1||\xi_2|,$$
where the constant above depends on $\X$, as the portions $\Gamma_n$ are determined by $\X$.
If $F\in\Ei$, then we may express $F$ locally as the graph of a function determined by one of the maps $F_K$ that make up the mesh $\Th$; we also have that $F_K\in C^2$, as the family of meshes is regular of order $1$. That is, since $F\subset\partial K$ for some $K\in\Th$, there exists a (straight) reference face $\hat F$, such that $F = F_K(\hat F)$. Furthermore, there exists a straight approximating face $\tilde F = \tilde{F}_K(\hat F)$ ($\tilde F_K$ is the affine part of $F_K$), which provides us with a local coordinate system. As $\tilde F$ is flat, after a suitable change of coordinates, one has that $\tilde{F}\subset\{(x',0):x'\in\mathbb R^{d-1}\}$. Furthermore, without loss of generality, we may assume that $F$ does not intersect $\tilde{F}$, since in such a case, we may define another flat face $\tilde{F}_a$, and an invertible affine map $\tilde{A}:\tilde{F}_a\to\tilde{F}$ that consists only of rotation and translation, which does not effect the bounds that we are about to obtain (i.e., the Euclidean norm of the matrix $D\tilde{A}$ is equal to $1$). Let us denote $\tilde{F}' = \{x'\in\mathbb R^{d-1}:(x',0)\in\tilde{F}\}$. Now, defining $\varphi_{F_K}:\tilde{F}'\to\mathbb R$ by 
$$\varphi_{F_K}(x') = [F_K]^d\circ\tilde{F}_K^{-1}(x',0),\quad x'\in\tilde{F}',$$
we see that
$F = \{(x',\varphi_{F_K}(x')):x'\in\tilde{F}'\}$. Note that we have now shown that all $F\in\Eib$ are of class $C^2$, and furthermore, for any $K\in\Th$, $\partial K$ may be expressed as the finite union of the closures of $F\in\Eib$, and thus for all $K\in\Th$, $\partial K$ is piecewise $C^2$.

Furthermore, expressing $F$ as the zero level set of the function $h_{F_K}(x',x_d) = x_d-\varphi_{F_K}(x')$, we see that 
$$n_F = \frac{\nabla h_{F_K}}{|\nabla h_{F_K}|}= -\frac{(\nabla_{x'}\varphi_{F_K},-1)^T}{|(\nabla_{x'}\varphi_{F_K},-1)|}=-\frac{(\nabla_{x'}\varphi_{F_K},-1)^T}{\sqrt{1+|\nabla_{x'}\varphi_{F_K}|^2}}.$$
Then, since $\nabla_\Ta n_F^T=\nabla n_F^T-n_F\frac{\partial n_F^T}{\partial n_F}$, for any two vectors $\xi^1,\xi^2:F\to\mathbb R^d$ tangent to $F$ (with components $\xi^k_1,\ldots,\xi^k_d$, $k=0,1$), and hence orthogonal to $n_F$, it follows that $(\xi^1)^T\nabla_\Ta n_F^T\,\xi^2 = (\xi^1)^T\nabla n_F^T\,\xi^2$.
%\begin{equation*}
%\begin{aligned}
%(\xi^1)^T\nabla_\Ta n_F\xi^2 & = \sum_{i,j=1}^d\frac{\partial[n_F]^j}{\partial x_i}\xi^1_i\xi^2_j-\sum_{i,j=1}^d[n_F]^i\sum_{k=1}^d\frac{\partial[n_F]^j}{\partial x_k}[n_F]^k\xi^1_i\xi^2_j\\
%& = \sum_{i,j=1}^d\frac{\partial[n_F]^j}{\partial x_i}\xi^1_i\xi^2_j-\sum_{j=1}^d\sum_{k=1}^d\frac{\partial[n_F]^j}{\partial x_k}[n_F]^k\xi^2_j\sum_{i=1}^d\xi^1_i[n_F]^i\\
%& = \sum_{i,j=1}^d\frac{\partial[n_F]^j}{\partial x_i}\xi^1_i\xi^2_j-\sum_{j=1}^d\sum_{k=1}^d\frac{\partial[n_F]^j}{\partial x_k}[n_F]^k\xi^2_j(\xi^1\cdot n_F) = \sum_{i,j=1}^d\frac{\partial[n_F]^j}{\partial x_i}\xi^1_i\xi^2_j.
%\end{aligned}
%\end{equation*}
%Thus
%$$(\xi^1)^T\nabla_\Ta n_F\xi^2=(\xi^1)^T\nabla n_F\xi^2.$$
Furthermore, denoting $\delta^{ij}:=1-\delta_{ij}$, where $\delta_{ij}$ is the Kronecker-delta symbol, we see that
\begin{equation*}
\begin{aligned}
\frac{\partial [n_F]^j}{\partial x_i} & = -\delta^{id}\frac{\partial}{\partial x_i}\left(\frac{(\delta^{jd}\frac{\partial \varphi_{F_K}}{\partial x_j}-\delta_{jd})}{\sqrt{|\nabla_{x'}\varphi_{F_K}|^2+1}}\right)\\
& = -\delta^{id}\frac{\delta^{jd}(|\nabla_{x'}\varphi_{F_K}|^2+1)\frac{\partial^2\varphi_{F_K}}{\partial x_i\partial x_j}-(\delta^{jd}\frac{\partial \varphi_{F_K}}{\partial x_j}-\delta_{jd})\sum_{k=1}^{d-1}\frac{\partial^2\varphi_{F_K}}{\partial  x_i\partial x_k}\frac{\partial \varphi_{F_K}}{\partial x_k}}{(|\nabla_{x'}\varphi_{F_K}|^2+1)^{3/2}},
\end{aligned}
\end{equation*}
and so
\begin{equation*}
\begin{aligned}
(\xi^1)^T\nabla n_F^T\,\xi^2 & = -\frac{\sum_{i,j=1}^{d-1}\frac{\partial^2\varphi_{F_K}}{\partial x_i\partial x_j}\xi^1_i\xi^2_j}{\sqrt{|\nabla_{x'}\varphi_{F_K}|^2+1}}+\frac{\sum_{i,k=1}^{d-1}\xi_i^1\frac{\partial^2\varphi_{F_K}}{\partial  x_i\partial x_k}\frac{\partial \varphi_{F_K}}{\partial x_k}\sum_{j=1}^d(\delta^{jd}\frac{\partial \varphi_{F_K}}{\partial x_j}-\delta_{jd})\xi^2_j}{(|\nabla_{x'}\varphi_{F_K}|^2+1)^{3/2}}\\
& = -\frac{\sum_{i,j=1}^{d-1}\frac{\partial^2\varphi_{F_K}}{\partial x_i\partial x_j}\xi^1_i\xi^2_j}{\sqrt{|\nabla_{x'}\varphi_{F_K}|^2+1}}-\frac{\sum_{i,k=1}^{d-1}\xi_i^1\frac{\partial^2\varphi_{F_K}}{\partial  x_i\partial x_k}\frac{\partial \varphi_{F_K}}{\partial x_k}(\xi^2\cdot n_F)}{(|\nabla_{x'}\varphi_{F_K}|^2+1)}\\
& = -\frac{\sum_{i,j=1}^{d-1}\frac{\partial^2\varphi_{F_K}}{\partial x_i\partial x_j}\xi^1_i\xi^2_j}{\sqrt{|\nabla_{x'}\varphi_{F_K}|^2+1}}\le\|D^2_{x'}\varphi_{F_K}\||\xi^1||\xi^2|.
\end{aligned}
\end{equation*}
One also has that
\begin{equation*}
\begin{aligned}
\sup_{x'\in\tilde{F}'}\|D^2_{x'}\varphi_{F_K}(x')\| &\le \sup_{x'\in\tilde{F}'}\|D^2_{x'}(F_K\circ \tilde{F}_K^{-1})(x',0)\|\\
&\le \sup_{x\in\tilde{F}}\|D^2_{x}(F_K\circ \tilde{F}_K^{-1})(x)\|\\
&= \sup_{x\in\tilde{F}}\|D^2F_K\circ\tilde{F}_K^{-1}(x)(\tilde{B}_K^{-1})^2\|\\
& = \sup_{x\in\hat F}\|D^2F_K(\hat x)(\tilde{B}_K^{-1})^2\|\\
&\le \sup_{x\in\hat F}\|D^2F_K(\hat x)\|\|\tilde{B}_K^{-1}\|^2=c_2\le C_{\operatorname{int}},
\end{aligned}
\end{equation*}
where the final inequality follows from~(\ref{regoforderms}), and $C_{\operatorname{int}}$ is independent of both $h$, and the choice of $F$, since the family of meshes is regular of order $1$. Thus, defining $C:=\max\{C(\X),C_{\operatorname{int}}\}$. For all $F\in\Eib$, we have
\begin{equation}\label{bound2ff18:26}
(\xi^1)^T\nabla_\Ta n_F^T\,\xi^2 \le C|\xi^1||\xi^2|,
\end{equation}
on $F$, for any tangent vectors to $F$.
Upon noting that $\nabla_\Ta u$, and $\nabla_\Ta v$ are tangent vectors to $F$, we obtain~(\ref{bounded2ndff}).$\quad\quad\square$

\begin{lemma}
Assume that $\X$ is piecewise $C^2$, and let $\{\Th\}_{h>0}$ be a family of meshes on $\overline{\X}$  that satisfies Assumption~\ref{Meshconds}. Then, there exists a constant $C$ depending on the family $\{\Th\}_{h>0}$, $d$, and $\X$ such that for any $F\in\Eib$, the following estimates hold on $F$:
\begin{align}
\sup_{x\in F}|\mc_F(x)|\le C(d)\sup_{x\in F}|\nabla_\Ta n_F^T(x)|&\le C,\label{1.1in}\\
\left|\nablaTa\left(\tau_F\frac{\partial v}{\partial n_F}\right)\right| & \le C(|\tau_F(D^2v)|+|\tau_F(\nabla v)|),\label{pre5.5.5}\\
|\operatorname{div}_\Ta\nabla_\Ta\tau_F(v)| & \le C(|\tau_F(D^2v)|+|\tau_F(\nabla v)|),\label{18:43ddd}
\end{align}
for all $v\in H^s(K)$, $s>5/2$, where $F\subset\partial K$, and $\tau_F$ is the trace operator from $K$ to $F$.
%\end{equation}
%\begin{equation}\label{pre5.5.5}
%\nablaTa\left(\tau_F\frac{\partial v}{\partial n_F}\right) = \sum_{k=1}^{d-1}((t_k)^T\tau_F(D^2v)\,n_F+(t_k)^T\nabla n_F^T\nablaTa (\tau_F v))t_k,\,\,\forall v\in H^s(K), s>5/2,
%\end{equation}
%where $\{t_k\}_{i=1}^{d-1}$ is an orthonormal coordinate system on $F$.
%\begin{equation}\label{pre5.5.5}
%\nablaTa\left(\tau_F\frac{\partial v}{\partial n_F}\right) = \sum_{k=1}^{d-1}((t_k)^T\tau_F(D^2v)\,n_F+(t_k)^T\nabla n_F^T\nablaTa (\tau_F v))t_k,\,\,\forall v\in H^s(K), s>5/2,
%\end{equation}
%where $\{t_k\}_{i=1}^{d-1}$ is an orthonormal coordinate system on $F$.
\end{lemma}
\emph{Proof:} Let $F\in\Eib$. Then, by definition, we see that
\begin{equation}\label{1.1}
\sup_{x\in F}|\mc_F(x)| = \sup_{x\in F}|\nabla_\Ta\cdot n_F(x)|\le C(d)|\sup_{x\in F}\nabla_\Ta n_F^T(x)|.
\end{equation}
Furthermore, let us take $\xi^1,\xi^2\in\mathbb R^d$, and decompose them in terms of their tangential and normal components, i.e., $\xi^k = (\xi^k)_\Ta+(\xi^k_{n_F})n_F$, $k=1,2$. Then, we see that on $F$
\begin{equation}\label{bball18:38}
\begin{aligned}
(\xi^1)^T\nabla_\Ta n_F^T\,\xi^2 & =  (\xi^1_\Ta)^T\nabla_\Ta n_F^T\,\xi^2_\Ta\\
&\le C|\xi^1_\Ta||\xi^2_\Ta|\le C|\xi^1||\xi^2|,
\end{aligned}
\end{equation}
where the penultimate inequality is due to~(\ref{bound2ff18:26}), as $(\xi^k)_\Ta$ are tangent vectors. Since this holds for all $\xi^1,\xi^2\in\mathbb R^d$, we deduce that 
\begin{equation}\label{1.2}
\sup_{x\in F}|\nabla_\Ta n_F^T(x)|\le C,
\end{equation}
which, combined with~(\ref{1.1}) yields
\begin{equation}\label{1.3}
\sup_{x\in F}|\mc_F(x)|\le C(d)\sup_{x\in F}|\nabla_\Ta n_F^T(x)|\le C,
\end{equation}
which is~(\ref{1.1in}).

Then, by~(\ref{identitieslap}) and~(\ref{1.3}) we see that on $F$
\begin{equation*}
\begin{aligned}
|\operatorname{div}_\Ta\nabla_\Ta\tau_F(v)| & = \left|\tau_F(\Delta v)+\mc_F\tau_F\frac{\partial v}{\partial n_F}+\tau_F\frac{\partial^2 v}{\partial n_F^2}\right|\\
&\le C(d)|\tau_F(D^2v)|+\sup_{x\in F}|\mc_F(x)||\tau_F(\nabla u)|\\
&\le C(|\tau_F(D^2v)|+|\tau_F(\nabla v)|),\\
\end{aligned}
\end{equation*}
which is~(\ref{18:43ddd}). Finally, from~(\ref{1.3}) we obtain the following
\begin{equation*}
\begin{aligned}
\left|\nabla_\Ta\tau_F\left(\frac{\partial v}{\partial n_F}\right)\right| & = \left|\nabla_\Ta(\tau_F(Dv))\cdot n_F+(\nabla_\Ta n_F^T)\cdot\nabla v\right|\\
&\le |\nabla_\Ta\tau_F(Dv)|+|\nabla_\Ta n_F^T||\nabla v|\\
&\le C(|\tau_F(D^2v)|+|\tau_F(\nabla v)|),
\end{aligned}
\end{equation*}
which is~(\ref{pre5.5.5}).$\quad\quad\square${}
\end{section}
\begin{section}{Finite element schemes and stability estimates}\label{sec:4}
We now provide DGFEMs for the approximation of solutions to~(\ref{2and4}) for $k=1,2$. For the case $k=1$, we will not need to alter the bilinear form $A_1:V_h\times V_h\to\mathbb R$, given by~(\ref{Akdef}). However, as mentioned in Section~\ref{sec:2}, obtaining an estimate in a $H_\Delta$-type norm (given by~(\ref{needsalabel14:48})) in the case $k=2$ does not seem possible when considering curved finite elements, due to the form that the inverse inequality takes. This means that we must define a different bilinear form, which relies on a discrete analogue of the following identity
\begin{equation}\label{calderon-1}
\int_\X\Delta u\,\Delta v = \int_\X D^2u\!:\!D^2v,\quad\forall u,v\in H^2_0(\X).
\end{equation}
Indeed, assuming that $\partial\X$ is Lipschitz continuous, the above estimate follows from an application of integration by parts (twice), for $i,j=1,\ldots,d$, we see that for $u,v\in C^\infty_c(\X)$
\begin{equation}\label{calderon}
\begin{aligned}
\int_\X\partial^2_{ij}u\,\partial^2_{ij}v & = \int_{\partial\X}\partial^2_{ij}u\,\partial_jv\,[n_{\partial\X}]^i-\int_K\partial^3_{iji}u\,\partial_jv\\
& = \int_\X\partial^2_{ii}u\,\partial^2_{jj}v-\int_{\partial\X}[\partial^2_{ii}u\,\partial_jv\,[n_{\partial\X}]^j-\partial^2_{ij}w\,\partial_jv\,[n_{\partial\X}]^i]\\
& = \int_\X\partial^2_{ii}u\,\partial^2_{jj}v,
\end{aligned}
\end{equation}
where the last equality is due to the fact that $v|_{\partial\X}=\partial_jv|_{\partial\X}=0$ for $j=1,\ldots,d$. Summing~(\ref{calderon}) over all $i,j=1,\ldots,d$, we obtain~(\ref{calderon-1}). Furthermore~(\ref{calderon-1}) extends to $u,v\in H^2_0(\X)$ by density, and, coupled with the Poincar\'{e} inequality, allows one to prove that the biharmonic problem~(\ref{2and4}) (for $k=2$) is well posed in $H^2_0(\X)$.

Let us define the bilinear form $C:V_{2,h}\times V_{2,h}\to\mathbb R$ as follows:
\begin{equation}\label{Cdef}
\begin{aligned}
C(u_h,v_h)&:=\sum_{F\in\Ei}\left\langle\Delta_\Ta\avg{u_h},\bjump{\frac{\partial v_h}{\partial n_F}}\right\rangle_F+\left\langle\mc_F\bavg{\frac{\partial u_h}{\partial n_F}},\bjump{\frac{\partial v_h}{\partial n_F}}\right\rangle_F-\left\langle\nabla_\Ta\bavg{\frac{\partial u_h}{\partial n_F}},\jump{\nabla_\Ta v_h}\right\rangle_F\\
&~~~~~~~~+\sum_{F\in\Ei}\langle\mathcal{Q}_F(\nabla_\Ta\avg{u_h},\jump{\nabla_\Ta v_h})+\mathcal{Q}_F(\avg{\nabla u_h\cdot n_F}n_F,\jump{\nabla_\Ta v_h})\rangle_F,\quad u_h,v_h\in V_h,
\end{aligned}
\end{equation}
where $n_F$ is a fixed choice of unit normal to $F$, $\mc_F:=\nabla_\Ta\cdot n_F$, and $\mathcal{Q}:\mathbb R^d\times\mathbb R^d\to\mathbb R$ is defined by
\begin{equation}\label{Q:def}
\mathcal{Q}_F(\xi_1,\xi_2) := \xi_1^T\nabla n_F^T\,\xi_2,\quad\xi_1,\xi_2\in\mathbb R^d.
\end{equation}
We now show that $C$ satisfies the following consistency identity.
\begin{lemma}
Assume that $\X$ is piecewise $C^2$ and that $(\Th)_{h>0}$ is a regular of order $1$ meshes on $\overline{\X}$ satisfying Assumption~\ref{Meshconds}. Then, the bilinear form $C:V_{2,h}\times V_{2,h}\to\mathbb R$ satisfies the following consistency identity:
\begin{equation*}
\begin{aligned}
&\sum_{K\in\Th}\langle\Delta w,\Delta v_h\rangle_K = \sum_{K\in\Th}\langle D^2w,D^2v_h\rangle_K+C(w,v_h),\quad\forall w\in H^2_0(\X)\cap H^4(\X),\,\forall v_h\in V_{2,h}.
\end{aligned}
\end{equation*}
\end{lemma}
\emph{Proof:} Under the hypotheses of the Lemma, it follows that an arbitrary $K\in\Th$ is Lipschitz continuous, and piecewise $C^2$, and, thus for $w\in H^4(\X)\cap H^2_0(\X)$, $v_h\in V_{2,h}$, we also have that for $i,j=1,\ldots,d$,
\begin{equation*}
\int_K\partial^2_{ij}u\,\partial^2_{ij}v_h = \int_K\partial^2_{ii}u\,\partial^2_{jj}v_h-\int_{\partial K}[\partial^2_{ii}u\,\partial_jv_h-\partial^2_{ij}w\,\partial_jv_h[n_{\partial K}]^i],
\end{equation*}
summing this expression over all $i,j=1,\ldots,d$ yields
$$\int_KD^2w\!:\!D^2v+\int_{\partial K}\Delta w\frac{\partial v_h}{\partial n_{\partial K}}-\nabla(\nabla w\cdot n_{\partial K})\cdot\nabla v_h+(\nabla w)^T\nabla n_{\partial K}^T\nabla v_h = \int_K\Delta w\,\Delta v_h.$$
Summing the above over all $K\in\Th$, we obtain
\begin{equation*}
\begin{aligned}
\sum_{K\in\Th}\langle\Delta w,\Delta v_h\rangle_K & = \sum_{K\in\Th}\langle D^2w,D^2v_h\rangle_K\\
&+\sum_{F\in\Eib}\langle \avg{\Delta w},\jump{\nabla v_h\cdot n_F}\rangle_F-\langle\avg{\nabla(\nabla w\cdot n_F)},\jump{\nabla v_h}\rangle_F+\langle\nabla n_F^T\jump{\nabla v_h},\avg{\nabla w}\rangle_F\\
&+\sum_{F\in\Ei}\langle \jump{\Delta w},\avg{\nabla v_h\cdot n_F}\rangle_F-\langle\jump{\nabla(\nabla w\cdot n_F)},\avg{\nabla v_h}\rangle_F+\langle\nabla n_F^T\avg{\nabla v_h},\jump{\nabla w}\rangle_F.
\end{aligned}
\end{equation*}
Since $w\in H^4(\X)$, it follows that
$$\sum_{F\in\Ei}\langle \jump{\Delta w},\avg{\nabla v_h\cdot n_F}\rangle_F-\langle\jump{\nabla(\nabla w\cdot n_F)},\avg{\nabla v_h}\rangle_F+\langle\nabla n_F^T\avg{\nabla v_h},\jump{\nabla w}\rangle_F=0.$$
Furthermore, since $w\in H^2_0(\X)$, for $\Eb\ni F\subset\partial\X$, we have that $w|_F=\partial_jw|_F=0$, for $j=1,\ldots,d$, and thus
\begin{equation*}
\begin{aligned}
&\sum_{F\in\Eb}\langle \avg{\Delta w},\jump{\nabla v_h\cdot n_F}\rangle_F-\langle\avg{\nabla(\nabla w\cdot n_F)},\jump{\nabla v_h}\rangle_F+\langle\nabla n_F^T\jump{\nabla v_h},\avg{\nabla w}\rangle_F\\
&=\sum_{F\in\Eb}\left\langle\frac{\partial^2w}{\partial n_F^2},\frac{\partial v_h}{\partial n_F}\right\rangle_F-\left\langle\frac{\partial^2w}{\partial n_F^2}n_F,\frac{\partial v_h}{\partial n_F}n_F+{\nabla_\Ta}v_h\right\rangle_F = 0.
\end{aligned}
\end{equation*}
It then follows that
\begin{equation*}
\begin{aligned}
\sum_{K\in\Th}\langle\Delta w,\Delta v_h\rangle_K & = \sum_{K\in\Th}\langle D^2w,D^2v_h\rangle_K\\
&+\sum_{F\in\Ei}\langle \avg{\Delta w},\jump{\nabla v_h\cdot n_F}\rangle_F-\langle\avg{\nabla(\nabla w\cdot n_F)},\jump{\nabla v_h}\rangle_F+\langle\nabla n_F^T\jump{\nabla v_h},\avg{\nabla w}\rangle_F.
\end{aligned}
\end{equation*}
We also see that
\begin{equation*}
\begin{aligned}
&\sum_{F\in\Ei}\langle\avg{\Delta w},\jump{\nabla v_h\cdot n_F}\rangle_F-\langle\avg{\nabla(\nabla w\cdot n_F)},\jump{\nabla v_h}\rangle_F+\langle\nabla n_F^T\jump{\nabla v_h},\avg{\nabla w}\rangle_F\\
& = \sum_{F\in\Ei}\left[\left\langle{\Delta_\Ta}\avg{w}+\mc_F\bavg{\frac{\partial w}{\partial n_F}}+\bavg{\frac{\partial^2 w}{\partial n_F^2}},\bjump{\frac{\partial v_h}{\partial n_F}}\right\rangle_F\right.\\
&~~~~~~~~-\left\langle\bavg{\frac{\partial^2w}{\partial n_F^2}}n_F+\nabla_\Ta\bavg{\frac{\partial w}{\partial n_F}},\bjump{\frac{\partial v_h}{\partial n_F}}n_F+\jump{\nabla_\Ta v_h}\right\rangle_F\\
&~~~~~~~~+\left.\left\langle\nabla n_F^T\left(\bjump{\frac{\partial v_h}{\partial n_F}}n_F+\jump{\nabla_\Ta v_h}\right),\nabla_\Ta\avg{w}+\bavg{\frac{\partial w}{\partial n_F}}n_F\right\rangle_F\right]\\
&=\sum_{F\in\Ei}\left[\langle\Delta_\Ta\avg{w},\jump{\nabla v_h\cdot n_F}\rangle_F+\langle\mc_F\avg{\nabla w\cdot n_F},\jump{\nabla w\cdot n_F}\rangle_F-\langle\nabla_\Ta\avg{\nabla w\cdot n_F},\jump{\nabla_\Ta v_h}\rangle_F\right.\\
&~~~~~~~~+\left.\langle\mathcal{Q}_F(\nabla_\Ta\avg{w},\jump{\nabla_\Ta v_h})+\mathcal{Q}_F(\avg{\nabla w\cdot n_F},\jump{\nabla_\Ta v_h})\rangle_F\right],
\end{aligned}
\end{equation*}
where we recall that $\mathcal{Q}$ is defined by~(\ref{Q:def}).
Finally, we see that
\begin{equation*}
\begin{aligned}
&\sum_{K\in\Th}\langle\Delta w,\Delta v_h\rangle_K = \sum_{K\in\Th}\langle D^2w,D^2v_h\rangle_K\\
&~~~~~~~~+\sum_{F\in\Ei}\left[\langle\Delta_\Ta\avg{w},\jump{\nabla v_h\cdot n_F}\rangle_F+\langle\mc_F\avg{\nabla w\cdot n_F},\jump{\nabla w\cdot n_F}\rangle_F-\langle\nabla_\Ta\avg{\nabla w\cdot n_F},\jump{\nabla_\Ta v_h}\rangle_F\right.\\
&~~~~~~~~+\left.\langle\mathcal{Q}_F(\nabla_\Ta\avg{w},\jump{\nabla_\Ta v_h})+\mathcal{Q}_F(\avg{\nabla w\cdot n_F}n_F,\jump{\nabla_\Ta v_h})\rangle_F\right]=\sum_{K\in\Th}\langle D^2w,D^2v_h\rangle_K+C(w,v_h),
\end{aligned}
\end{equation*}
for any $w\in H^4(\X)\cap H^2_0(\X)$, $v_h\in V_{2,h}$.$\quad\quad\square$
\subsection{Numerical methods}
We are now ready to provide the finite element schemes for the approximation of solutions to the Poisson, and biharmonic problem~(\ref{2and4}) for $k=1$ and $k=2$, respectively. Note that in the sequel we set $V_{k,h}:=V_{h,p}$, where we assume that $p\ge k$.
\subsubsection{Poisson problem} One seeks $u_{1,h}\in V_{1,h}$ such that
\begin{equation}\label{A1def:14:47}
\begin{aligned}
&A_1(u_{1,h},v_h):=\sum_{K\in\Th}\langle\nabla u_{1,h},\nabla v_h\rangle_K+B_1(u_{1,h},v_h)+B_1(v_h,u_{1,h})+J_1(u_{1,h},v_h)=\ell(v_h),
\end{aligned}
\end{equation}
for all $v_h\in V_{1,h}$, and we recall that $B_1,J_1:V_{1,h}\times V_{1,h}\to\mathbb R$ are defined as follows:
\begin{align*}
B_1(u_h,v_h)&:=-\sum_{F\in\Eib}\int_F\bavg{\frac{\partial u_h}{\partial n_F}}\jump{v_h},\\
J_1(u_h,v_h)&:=\sum_{F\in\Eib}\frac{\eta^1_F}{\tilde{h}_F}\langle\jump{u_h},\jump{v_h}\rangle_F
\end{align*}
where $n_F$ is a \emph{fixed} choice of unit normal to $F$, and $\eta_F^1$ is a positive face dependent parameter.
\subsubsection{Biharmonic problem}
One seeks $u_{2,h}\in V_{2,h}$ such that
 \begin{equation}\label{A2def:14:47}
\begin{aligned}
&A_2(u_{2,h},v_h):=\sum_{K\in\Th}\langle D^2u_{2,h},D^2v_h\rangle_K+C(u_{2,h},v_h)+C(v_h,u_{2,h})+B_2(u_{2,h},v_h)+B_2(v_h,u_{2,h})\\
&~~~~~~~~~~~~~~~~~~~~~~~~~~~~~~~~~+J_2(u_{2,h},v_h)=\ell(v_h),\quad\forall v_h\in V_{2,h},
\end{aligned}
\end{equation}
 where $B_2,J_2:V_{2,h}\times V_{2,h}\to\mathbb R$ are defined as follows:
\begin{align*}
B_2(u_h,v_h)&:=\sum_{F\in\Eib}\int_F\bavg{\frac{\partial(\Delta u_h)}{\partial n_{F}}}\jump{v_h}-\avg{\Delta u_h}\bjump{\frac{\partial v_h}{\partial n_{F}}},\\
J_2(u_h,v_h)&:=\sum_{F\in\Eib}\frac{\eta^2_F}{\tilde{h}_F^3}\langle\jump{u_h},\jump{v_h}\rangle_F+\frac{\eta^3_F}{\tilde{h}_F}\langle\jump{\nabla u_h\cdot n_F},\jump{\nabla v_h\cdot n_F}\rangle_F+\frac{\eta^4_F}{\tilde{h}_F}\langle\jump{\nabla_\Ta u_h},\jump{\nabla_\Ta v_h}\rangle_F,
\end{align*}
where $\eta_F^2,\eta_F^3,\eta_F^4$ are positive face dependent terms to be provided,
and $C:V_{2,h}\times V_{2,h}\to\mathbb R$ is defined by~(\ref{Cdef}).
Note that~(\ref{A2def:14:47}) is obtained by replacing $\sum_{K\in\Th}\langle\Delta u_{2,h},\Delta v_h\rangle_K$ with $\sum_{K\in\Th}\langle D^2 u_{2,h},D^2 v_h\rangle_K$\\$+C(u_{2,h},v_h)+C(v_h,u_{2,h})$, and including a tangential gradient penalty term in $J_2$, in the definition of $A_2$ given by~(\ref{Akdef}), and that the $C(v_h,u_{2,h})$ term results in $A_2$ defined by~(\ref{A2def:14:47})  being symmetric (whilst preserving the consistency of the scheme).
\begin{remark}[Comparison to the method of~\cite{MR2142191}]
In the method given by~(\ref{A2def:14:47}), we have taken $V_{2,h}$ to be the space of discontinuous piecewise polynomials. However, if we instead use the space $V_{h,c}:=V_{2,h}\cap H^1_0(\X)$, and assume that $\X$ is polygonal, then the scheme~(\ref{A2def:14:47}) coincides with the $C^0$-interior penalty method proposed in~\cite{MR2142191} (without the second order term, see (4.15) of~\cite{MR2142191}).
\end{remark}
\subsection{Stability estimates}
We now provide stability estimates for the finite element methods~(\ref{A1def:14:47}) and~(\ref{A2def:14:47}), which yield existence and uniqueness of a numerical solution. Let us recall the definitions of the norms that the stability estimates will hold in. We define the norms $\|\cdot\|_{h,k}:V_{k,h}\to\mathbb[0,\infty)$ as follows:
\begin{equation*}
\|v_h\|_{h,k}^2:=|v_h|_{H^k(\X;\Th)}^2+C_{*,k}J_k(v_h,v_h)
\end{equation*}
for positive constants $C_{*,k}$ are to be determined.
\begin{remark}[Trace estimates for jumps and averages]
We note that for $F\in\Ei$, $v_h\in V_{k,h}$, $k=1,2$, we have that $F=\overline{K^+}\cap\overline{K^-}$ for some $K^+,K^-\in\Th$, and denoting $v_h^\pm:= v_h|_{K^\pm}$ that 
$$\avg{D^\alpha v_h} = \frac{1}{2}(D^\alpha v_h^+|_F+D^\alpha v_h^-|_F)\quad\mbox{and}\quad\jump{D^\alpha v_h} = D^\alpha v_h^+|_F-D^\alpha v_h^-|_F,$$
where the multi-index $\alpha$ satisfies $|\alpha|\le m-1$, assuming that $\X$ is piecewise $C^m$. Then, momentarily denoting $\{\cdot\}$ to be either the jump or average operator, it follows that for any $s\ge 0$
\begin{equation}\label{javgbnd1in}
\begin{aligned}
\sum_{F\in\Ei}\tilde{h}_F^s\|\{D^\alpha v_h\}\|_{L^2(F)}^2&\le2\sum_{F\in\Ei}\tilde{h}_F^s\|D^\alpha v^+|_F\|_{L^2(F)}^2+\tilde{h}_F^s\|D^\alpha v^-|_F\|_{L^2(F)}^2\\
&\le2\sum_{F\in\Ei}\tilde{h}_F^s\|D^\alpha v^+|_{\partial K^+}\|_{L^2(\partial K^+)}^2+\tilde{h}_F^s\|D^\alpha v^-|_{\partial K^-}\|_{L^2(\partial K^-)}^2\\
&\le 2\sum_{F\in\Ei}\sum_{K\in\Th:F\subset\partial K}\tilde{h}_F^s\|D^\alpha v\|_{L^2(\partial K)}^2\\
&\le 2C_{\operatorname{Tr}}\sum_{F\in\Ei}\sum_{K\in\Th:F\subset\partial K}\tilde{h}_F^s(h_K^{-1}\|D^\alpha v\|_{L^2(K)}^2+h_K|D^\alpha v|_{H^1(K)}^2)\\
&\le 2C_{\operatorname{Tr}}\sum_{F\in\Ei}\sum_{K\in\Th:F\subset\partial K}h_K^s(h_K^{-1}\|D^\alpha v\|_{L^2(K)}^2+h_K|D^\alpha v|_{H^1(K)}^2),
\end{aligned}
\end{equation}
where the penultimate inequality follows from~(\ref{1simptrace}), and the final inequality is due to the fact that if $F\subset\partial K$, then $\tilde{h}_F\le h_K$. Furthermore, if $F\in\Eb$, then $F\subset\partial K$ for some $K\in\Th$, and
$$\avg{D^\alpha v_h} = \jump{D^\alpha v_h} = D^\alpha v_h|_F,\quad|\alpha|\le m-1,$$
and similarly we obtain
\begin{equation}\label{javgbnd1ext}
\sum_{F\in\Eb}\tilde{h}_F^s\|\{D^\alpha v_h\}\|_{L^2(F)}^2\le C_{\operatorname{Tr}}\sum_{F\in\Eb}\sum_{K\in\Th:F\subset\partial K}h_K^s(h_K^{-1}\|D^\alpha v\|_{L^2(K)}^2+h_K|D^\alpha v|_{H^1(K)}^2).
\end{equation}
Then, since the number of faces that make up the boundary of a simplex is bounded in terms of the dimension, $d$, the estimates~(\ref{javgbnd1in}) and~(\ref{javgbnd1ext}) yield the following (note that for simplicity we absorb the constant $2$ into the constant $C(d)$):
\begin{align*}
\sum_{F\in\Eib}\tilde{h}_F^s\|\{D^\alpha v_h\}\|_{L^2(F)}^2&\le 2C_{\operatorname{Tr}}\sum_{F\in\Eib}\sum_{K\in\Th:F\subset\partial K}h_K^s(h_K^{-1}\|D^\alpha v\|_{L^2(K)}^2+h_K|D^\alpha v|_{H^1(K)}^2)\\
&\le C_{\operatorname{Tr}}C(d)\sum_{K\in\Th}h_K^s(h_K^{-1}\|D^\alpha v\|_{L^2(K)}^2+h_K|D^\alpha v|_{H^1(K)}^2).
\end{align*}
In the sequel we shall utilise the above estimate several times, for various orders of $\alpha$, as it simplifies the exposition of the proofs.
\end{remark}
\begin{lemma}
Assume that $\X\subset\mathbb R^d$ is piecewise $C^2$, and that $(\Th)_{h>0}$ is a regular of order $1$ family of triangulations on $\overline{\X}$ satisfying Assumption~\ref{Meshconds}. Then, for any $\kappa_1\in(0,1)$ there exists a constant $\sigma_1>0$ depending on $\kappa_1$, $C_{\operatorname{Tr}}$, $C_I$, and $d$, such that if
\begin{equation}\label{etacond1}
\eta^1_F\ge\sigma_1\quad\forall F\in\Eib,
\end{equation}
then
\begin{equation}\label{A1stab}
A_1(v_h,v_h)\ge\kappa_1\|v_h\|_{h,1}^2\quad\forall v_h\in V_{1,h}.
\end{equation}
Thus, there exists a unique $u_{1,h}\in V_{1,h}$ that satisfies~(\ref{A1def:14:47}).
\end{lemma}
\emph{Proof:} Utilising the trace and inverse estimates~(\ref{1simptrace}) and~(\ref{local:invineq}), as well as the Cauchy--Schwarz inequality with a parameter, we obtain for any $\delta_1>0$, and any $v_h\in V_{1,h}$
\begin{equation*}
\begin{aligned}
A_1(v_h,v_h)&\ge|v_h|^2_{H^1(\X;\Th)}+\sum_{F\in\Eib}\left[(\eta^1_F-\delta_1^{-1})\tilde h_F^{-1}\|\jump{v_h}\|_{L^2(F)}^2-\delta_1\tilde{h}_F\left\|\bavg{\frac{\partial v_h}{\partial n_F}}\right\|_{L^2(F)}^2\right]\\
&\ge|v_h|^2_{H^1(\X;\Th)}-\delta_1C_{\operatorname{Tr}}C(d)\sum_{K\in\Th}[|v_h|_{H^1(K)}^2+h_K|v_h|_{H^2(K)}^2]+\sum_{F\in\Eib}(\eta^1_F-\delta_1^{-1})\tilde h_F^{-1}\|\jump{v_h}\|_{L^2(F)}^2\\
&\ge\left(1-\delta_1(1+C_I)C_{\operatorname{Tr}}C(d)\right)|v_h|^2_{H^1(\X;\Th)}+\sum_{F\in\Eib}(\eta^1_F-\delta_1^{-1})\tilde h_F^{-1}\|\jump{v_h}\|_{L^2(F)}^2.
\end{aligned}
\end{equation*}
For a given $\kappa_1\in(0,1)$, one can choose $\delta_1$ sufficiently small, such that $1-\delta_1(1+C_I)C_{\operatorname{Tr}}C(d)\ge\kappa_1$, which yields
$$A_1(v_h,v_h)\ge\kappa_1|v_h|^2_{H^1(\X;\Th)}+\sum_{F\in\Eib}(\eta^1_F-\delta_1^{-1})\tilde h_F^{-1}\|\jump{v_h}\|_{L^2(F)}^2.$$
Then, if $\eta^1_F$ satisfies~(\ref{etacond1}) with $\sigma_1>\delta_1^{-1}$, setting $C_{*,1}:=(\sigma_1-\delta_1^{-1})/\kappa_1>0$, we obtain
$$A_1(v_h,v_h)\ge\kappa_1\|v_h\|_{h,1}^2,$$
as desired.$\quad\quad\square$
\begin{lemma}
Assume that $\X\subset\mathbb R^d$ is piecewise $C^4$, and that $(\Th)_{h>0}$ is a regular of order $3$ family of triangulations on $\overline{\X}$ satisfying Assumption~\ref{Meshconds}. Then, for any $\kappa_2\in(0,1)$ there exists a constant $\sigma_2>0$ depending on $\kappa_2$, $C_{\operatorname{Tr}}$, $C_I$, $C_P$, and $d$, such that if
\begin{equation}\label{etacond2}
\eta^2_F,\eta_F^3\ge\sigma_2\quad\forall F\in\Eib,
\end{equation}
then
\begin{equation}\label{A2stab2}
A_2(v_h,v_h)\ge\kappa_2\|v_h\|_{h,2}^2\quad\forall v_h\in V_{2,h},
\end{equation}
Thus, there exists a unique $u_{2,h}\in V_{2,h}$ that satisfies~(\ref{A2def:14:47}).
\end{lemma}
\emph{Proof:} Utilising the trace and inverse estimates~(\ref{1simptrace}) and~(\ref{local:invineq}), and the Cauchy--Schwarz inequality with a parameter, we obtain for any $\delta_2>0$, and any $v_h\in V_{2,h}$
\begin{equation}\label{B2bound16:37}
\begin{aligned}
|B_2(v_h,v_h)|&\le\frac{1}{2}\sum_{F\in\Eib}\left[\delta_2\tilde{h}_F^3\left\|\bavg{\frac{\partial (\Delta v_h)}{\partial n_F}}\right\|_{L^2(F)}^2+\delta_2\tilde{h}_F\|\avg{\Delta v_h}\|_{L^2(F)}^2\right.\\
&~~~~~~~~~~~~~~~~~~~~~~~~~~~~~~~~~~~~~+\left.(\delta_2\tilde{h}_F^3)^{-1}\|\jump{v_h}\|_{L^2(F)}^2+(\delta_2\tilde{h}_F)^{-1}\left\|\bjump{\frac{\partial v_h}{\partial n_F}}\right\|_{L^2(F)}^2\right]\\
&\le\delta_2(1+C_I)C_{\operatorname{Tr}}C(d)(|v_h|_{H^1(\X;\Th)}^2+|v_h|_{H^2(\X;\Th)}^2)\\
&~~~~~~~~~~~~~~~~~~~~~~~~~~~~~~~~~~~~~+\frac{1}{2\delta_2}\sum_{F\in\Eib}\left[\tilde{h}_F^{-3}\|\jump{v_h}\|_{L^2(F)}^2+\tilde{h}_F^{-1}\left\|\bjump{\frac{\partial v_h}{\partial n_F}}\right\|_{L^2(F)}^2\right],
\end{aligned}
\end{equation}
Furthermore, a further application of the Cauchy--Schwarz inequality with a parameter yields% from~(\ref{1.1in})--(\ref{18:43ddd}) we obtain %(denoting for $F\in\Ei$, $F=\overline{K^+}\cap\overline{K^-}$,)$v_h^\pm:=v_h|_{K^\pm}$
\begin{align*}
|C(v_h,v_h)|&\le\frac{\delta_2}{2}\sum_{F\in\Ei}\tilde{h}_F\left[\|\Delta_\Ta v_h\|_{L^2(F)}^2+\mc_F^2\left\|\frac{\partial v_h}{\partial n_F}\right\|_{L^2(F)}^2+\left\|\nabla_\Ta\bavg{\frac{\partial v_h}{\partial n_F}}\right\|_{L^2(F)}^2\right.\\
&\left.~~~~~~~~~~~~~~~~~~~~~~~~+\|\nabla n_F^T\|_{L^\infty(F)}^2\left(\|\nabla_\Ta\avg{v_h}\|_{L^2(F)}^2+\left\|\bavg{\frac{\partial v_h}{\partial n_F}}\right\|_{L^2(F)}^2\right)\right]\\
&~~~~~~~~~~~+\frac{1}{2\delta_2}\sum_{F\in\Ei}\tilde{h}_F^{-1}\left[2\left\|\bjump{\frac{\partial v_h}{\partial n_F}}\right\|_{L^2(F)}^2+3\|\jump{\nabla_\Ta v_h}\|_{L^2(F)}^2\right].
\end{align*}
Then, applying~(\ref{1.1in})--(\ref{18:43ddd}) in combination with the trace and inverse estimates~(\ref{1simptrace}) and~(\ref{local:invineq}) to the above estimate, we obtain
\begin{equation}\label{Cbound16:37}
\begin{aligned}
|C(v_h,v_h)|&\le \delta_2(1+C_I)C_{\operatorname{Tr}}C(d)(|v_h|_{H^1(\X;\Th)}^2+|v_h|_{H^2(\X;\Th)}^2)\\
&~~~~~~~~~+\frac{1}{2\delta_2}\sum_{F\in\Ei}\tilde{h}_F^{-1}\left[2\left\|\bjump{\frac{\partial v_h}{\partial n_F}}\right\|_{L^2(F)}^2+3\|\jump{\nabla_\Ta v_h}\|_{L^2(F)}^2\right].
\end{aligned}
\end{equation}
From~(\ref{B2bound16:37}) and~(\ref{Cbound16:37}), and the discrete Poincar\'{e}--Friedrichs' inequality~(\ref{GPF:est:disc}) it then follows that for any $\delta_2>0$ and any $v_h\in V_{2,h}$
\begin{align*}
A_h(v_h,v_h)&\ge|v_h|_{H^2(\X;\Th)}^2-2|B(v_h,v_h)|-2|C(v_h,v_h)|\\
&\ge(1-4(1+C_I)C_{\operatorname{Tr}}C(d)\delta_2)|v_h|_{H^2(\X;\Th)}^2-4(1+C_I)C_{\operatorname{Tr}}C(d)\delta_2|v_h|_{H^1(\X;\Th)}^2\\
&~~~~~~~-\frac{1}{2\delta_2}\sum_{F\in\Eib}6\tilde{h}_F^{-1}\left[\left\|\bjump{\frac{\partial v_h}{\partial n_F}}\right\|_{L^2(F)}^2+\|\jump{\nabla_\Ta v_h}\|_{L^2(F)}^2\right]+2\tilde{h}_F^{-3}\|\jump{v_h}\|_{L^2(F)}^2\\
&~~~~~~~~~~~~~+\sum_{F\in\Eib}\frac{\eta^2_F}{\tilde{h}_F^3}\|\jump{v_h}\|_{L^2(F)}^2+\frac{\eta^3_F}{\tilde{h}_F}\left\|\bjump{\frac{\partial v_h}{\partial n_F}}\right\|_{L^2(F)}^2+\frac{\eta_F^4}{\tilde{h}_F}\|\jump{\nabla_\Ta v_h}\|_{L^2(F)}^2\\
\end{align*}
\begin{align*}
&\ge(1-4(1+C_I)(1+C_P)C_{\operatorname{Tr}}C(d)\delta_2)|v_h|_{H^2(\X;\Th)}^2\\
&~~~~~~~~~~~~~+\sum_{F\in\Eib}\left(\frac{\eta^2_F}{\tilde{h}_F^3}-\frac{1}{\delta_2\tilde{h}_F^3}-\frac{4(1+C_I)C_{\operatorname{Tr}}C_P\delta_2}{\tilde{h}_F}\right)\|\jump{v_h}\|_{L^2(F)}^2\\
&~~~~~~~~~~~~~+\sum_{F\in\Eib}\frac{1}{\tilde{h}_F}\left(\eta^3_F-\frac{3}{\delta_2}-4(1+C_I)C_{\operatorname{Tr}}C_P\delta_2\right)\left\|\bjump{\frac{\partial v_h}{\partial n_F}}\right\|_{L^2(F)}\\
&~~~~~~~~~~~~~+\sum_{F\in\Eib}\frac{1}{\tilde{h}_F}\left(\eta^4_F-\frac{3}{\delta_2}\right)\|\jump{\nabla_\Ta v_h}\|_{L^2(F)}^2.
\end{align*}
Then, for any given $\kappa_2\in(0,1)$, we may choose $\delta_2$ sufficiently small so that $$1-4(1+C_I)(1+C_P)C_{\operatorname{Tr}}C(d)\delta_2\ge\kappa_2,$$ and, for such $\delta_2$, if $\eta^2_F,\eta^3_F,\eta^4_F$ satisfy~(\ref{etacond2}), with
$$\sigma_2>\frac{3}{\delta_2}+4(1+C_I)C_{\operatorname{Tr}}C_P\delta_2,$$
setting $C_{*,2}:=(\sigma_2-((3/\delta_2)+4(1+C_I)C_{\operatorname{Tr}}C_P\delta_2))/\kappa_2>0$, we obtain
$$A_1(v_h,v_h)\ge\kappa_2\|v_h\|_{h,2}^2,$$
as desired.$\quad\quad\square$
\end{section}
\section{Error analysis}\label{sec:5}
We now use the consistency of the two schemes to prove optimal a priori error estimates for the numerical solutions $u_{k,h}$, $k=1,2$, assuming that $\X$ is piecewise $C^{m+1}$, $m\in\mathbb N$, and that $u_k\in H^k_0(\X)\cap H^{2k}(\X)\cap H^{\mathbf{s}_k}(\X;\Th)$, where $\mathbf{s}_k = (s_K^k)_{K\in\Th}$ and each $s_K^k\ge2k$, are the true solutions of~(\ref{2and4}) for $k=1,2$. That is, we prove the following estimate
\begin{equation*}
\|u_{k,h}-u_k\|_{h,k}\le C_k\left(\sum_{K\in\Th}h_K^{2t_K^k-2k}\|u_k\|_{H^{s_K^k}(K)}^2\right)^{1/2},
\end{equation*}
where $t_K^k := \max\{p+1,m+1,s_K^k\}$.

Let us first recap on the approach we shall take, since this will in fact shorten the upcoming proofs. Let us take $z_{k,h}\in V_{k,h}$, $k=1,2$ to be arbitrary, denoting $\xi_{k,h}:=u_k-z_{k,h}$ and $\psi_{k,h}:=z_{k,h}-u_{k,h}$, we see that
$$\|u_k-u_{k,h}\|_{h,k}\le\|\xi_{k,h}\|_{h,k}+\|\psi_{k,h}\|_{h,k}.$$
Let us first estimate $\|\xi_{k,h}\|_{h,k}$. Due to the interpolation estimate~(\ref{opt:interp:2}) (since the choice of $z_{k,h}\in V_{k,h}$ is arbitrary) one can see that
\begin{equation}\label{xibound15:32}
\|\xi_{k,h}\|_{h,k}^2 = |\xi_{k,h}|_{H^k(\X;\Th)}^2+C_{*,k}J_k(\xi_{k,h},\xi_{k,h})\le C\sum_{K\in\Th}h_K^{2t_K^k}\|u_k\|_{H^{s_K^k}(K)}^2+C_{*,k}J_k(\xi_{k,h},\xi_{k,h}).
\end{equation}
Since there are constants $C_{\eta_k}$, $k=1,2$, satisfying $\eta_F^1\le C_{\eta_1}$ and $\eta_F^j\le C_{\eta_2}$, $j=2,3,4$, for all $F\in\Eib$, we have that
\begin{align*}
J_1(\xi_{1,h},\xi_{1,h}) & = \sum_{F\in\Eib}\frac{\eta^1_F}{\tilde{h}_F}\|\jump{\xi_{1,h}}\|_{L^2(F)}^2\\
&\le C_{\eta_1} C_{\operatorname{Tr}}\sum_{F\in\Eib}\sum_{K\in\Th:F\subset\partial K}\tilde{h}_F^{-1}(h_K^{-1}\|\xi_{1,h}\|_{L^2(K)}^2+h_K|\xi_{1,h}|_{H^1(K)}^2)\\
&\le C_{\eta_1}C_{\operatorname{Tr}}C_{\mathcal{T}}C(d)\sum_{K\in\Th}(h_K^{-2}\|\xi_{1,h}\|_{L^2(K)}^2+|\xi_{1,h}|_{H^1(K)}^2).
\end{align*}
Furthermore,
\begin{align*}
J_2(\xi_{2,h},\xi_{2,h}) & = \sum_{F\in\Eib}\frac{\eta^2_F}{\tilde{h}_F^3}\|\jump{\xi_{2,h}}\|_{L^2(F)}^2+\frac{\eta_F^3}{\tilde{h}_F}\left\|\bjump{\frac{\partial\xi_{2,h}}{\partial n_F}}\right\|_{L^2(F)}^2+\frac{\eta_F^4}{\tilde h_F}\|\jump{\nabla_\Ta\xi_{2,h}}\|_{L^2(F)}^2\\
&\le C_{\eta_2} C_{\operatorname{Tr}}C_{\mathcal{T}}C(d)\sum_{K\in\Th}(h_K^{-4}\|\xi_{2,h}\|_{L^2(K)}^2+h_K^{-2}|\xi_{2,h}|_{H^1(K)}^2+|\xi_{2,h}|_{H^2(K)}^2).
\end{align*}
Thus, from these two estimates, and an application of the interpolation estimate~(\ref{opt:interp:2}), we obtain
\begin{equation}\label{Jbound15:33}
\begin{aligned}
J_k(\xi_{k,h},\xi_{k,h})&\le C_{\eta_k} C_{\operatorname{Tr}}C_{\mathcal{T}}C(d)\sum_{j=0}^k\sum_{K\in\Th}h_K^{2(j-k)}|\xi_h|_{H^j(K)}^2\\
&\le CC_{\eta_k} C_{\operatorname{Tr}}C_{\mathcal{T}}C(d)\sum_{j=0}^k\sum_{K\in\Th}h_K^{2(j-k)}h_K^{2(t_K^k-j)}\|u_k\|_{H^{s_K^k}(K)}^2\\
&\le CC_{\eta_k} C_{\operatorname{Tr}}C_{\mathcal{T}}C(d)\sum_{K\in\Th}h_K^{2t_K^k-2k}\|u_k\|_{H^{s_K^k}(K)}^2.
\end{aligned}
\end{equation}
Applying~(\ref{Jbound15:33}) to~(\ref{xibound15:32}) and taking square roots, we obtain
\begin{equation}\label{xibnduseful}
\|\xi_{k,h}\|_{h,k}\le C_k\left(\sum_{K\in\Th}h_K^{2t_K^k-2k}\|u_k\|_{H^{s_K^k}(K)}^2\right)^{1/2},
\end{equation}
and so, it remains to estimate $\|\psi_{k,h}\|_{h,k}$, which relies upon the consistency of the schemes, and will be the objective of the next two proofs.
\begin{lemma}\label{errorest1lem}
Assume that $\X\subset\mathbb R^d$ is piecewise $C^{m+1}$ for some $m\in\mathbb N$, and that $(\Th)_{h>0}$ is a regular of order $m$ family of triangulations on $\overline{\X}$ satisfying Assumption~\ref{Meshconds}. Moreover, assume that $u_1\in H^1_0(\X)\cap H^{2}(\X)\cap H^{\mathbf{s}_1}(\X;\Th)$, where $\mathbf{s}_1 = (s_K^1)_{K\in\Th}$ and each $s_K^1\ge2$, is the true solution of~(\ref{2and4}) for $k=1$. Furthermore, let $\eta^1_F$ satisfy~(\ref{etacond1}) such that~(\ref{A1stab}) holds for some $\kappa_1\in(0,1)$. Then, the following estimate holds
\begin{equation}\label{regerrorest1}
\|u_{1,h}-u_1\|_{h,1}\le C_1\left(\sum_{K\in\Th}h_K^{2t_K^1-2}\|u_k\|_{H^{s_K^1}(K)}^2\right)^{1/2},
\end{equation}
where $t_K^1:=\min\{p+1,m+1,s_K^1\}$, $u_{1,h}\in V_{1,h}$ is the unique solution of~(\ref{A1def:14:47}), and the constant $C_1$ depends upon the shape-regularity constant of the mesh, $C_I$, $C_{\operatorname{Tr}}$, $C_\mathcal{T}$, $d$, $\mathbf{s}_K^1$, and $p$, but not upon $h_K$.
\end{lemma}
\emph{Proof:} Let us take $z_h\in V_{1,h}$, to be arbitrary, denoting $\xi_h:=u_1-z_h$ and $\psi_h:=z_h-u_{1,h}$, we see that
\begin{equation}\label{16:42::}
\|u_1-u_{1,h}\|_{h,1}\le\|\xi_h\|_{h,1}+\|\psi_h\|_{h,1}.
\end{equation}
Furthermore, from~(\ref{xibnduseful}) for $k=1$, we see that
\begin{equation}\label{otherxibnduseful}
\|\xi_h\|_{h,1}\le C_1\left(\sum_{K\in\Th}h_K^{2t_K^1-2}\|u_k\|_{H^{s_K^1}(K)}^2\right)^{1/2},
\end{equation}
and so it remains to estimate $\|\psi_h\|_{h,1}$. Due to the stability estimate~(\ref{A1stab}), and the consistency of the scheme, since $\psi_h\in V_{1,h}$, we see that
\begin{equation}\label{16:39}
\begin{aligned}
\|\psi_h\|_{h,1}^2&\le\kappa_1^{-1}A_1(\psi_h,\psi_h)\\
&=\kappa_1^{-1}(A_1(z_h,\psi_h)-A_1(u_{1,h},\psi_h))\\
&=\kappa_1^{-1}(A_1(z_h,\psi_h)-\ell(\psi_h))\\
&=\kappa_1^{-1}(A_1(z_h,\psi_h)-A_1(u_1,\psi_h))\\
&=\kappa_1^{-1}A_1(\xi_h,\psi_h).
\end{aligned}
\end{equation}
Then,
applying the Cauchy--Schwarz inequality for vectors in $\mathbb R^N$, estimate~(\ref{Jbound15:33}), and the interpolation estimate~(\ref{opt:interp:2}), we obtain
\begin{align*}
A_1(\xi_h,\psi_h) & = \sum_{K\in\Th}\langle\nabla \xi_h,\nabla \psi_h\rangle_K+B_1( \xi_h,\psi_h)+B_1( \psi_h, \xi_h)+J_1(\xi_h,\psi_h)\\
&\le|\xi_h|_{H^1(\X;\Th)}|\psi_h|_{H^1(\X;\Th)}+B_1( \xi_h,\psi_h)+B_1( \psi_h, \xi_h)+J_1(\xi_h,\xi_h)^{1/2}J_1(\psi_h,\psi_h)^{1/2}\\
&\le C\max\{1,C_{*,1}^{-1/2}\}\left(\sum_{K\in\Th}h_K^{2t_K^1-2}\|u_k\|_{H^{s_K^1}(K)}^2\right)^{1/2}\|\xi_h\|_{h,1}+B_1(\xi_h,\psi_h)+B_1( \psi_h, \xi_h).
\end{align*}
We again
apply the Cauchy--Schwarz inequality for vectors in $\mathbb R^N$, estimate~(\ref{Jbound15:33}), the interpolation estimate~(\ref{opt:interp:2}), the trace estimate~(\ref{1simptrace}), and the inverse estimate~(\ref{local:invineq}) yielding
\begin{align*}
&B_1(\xi_h,\psi_h)+B_1(\psi_h,\xi_h) \le\left(\sum_{F\in\Eib}\tilde{h}_F\left\|\bavg{\frac{\partial \xi_h}{\partial n_F}}\right\|_{L^2(F)}^2\right)^{1/2}\left(\sum_{F\in\Eib}\tilde{h}_F^{-1}\|\jump{\psi_h}\|_{L^2(F)}^2\right)^{1/2}\\
&~~~~~~~~~~~~~~~~~~~~~~~~~~~~~~~~~~~~~+\left(\sum_{F\in\Eib}\tilde{h}_F\left\|\bavg{\frac{\partial \psi_h}{\partial n_F}}\right\|_{L^2(F)}^2\right)^{1/2}\left(\sum_{F\in\Eib}\tilde{h}_F^{-1}\|\jump{\xi_h}\|_{L^2(F)}^2\right)^{1/2}\\
&\le C\left(\sum_{K\in\Th}h_K^{2t_K^1-1}\|u_1\|^2_{H^{s^k_K}(K)}\right)^\frac{1}{2}J_1(\psi_h,\psi_h)^\frac{1}{2}+\left(\sum_{K\in\Th}|\psi_h|_{H^1(K)}^2+h_K|\psi_h|_{H^2(K)}^2\right)^\frac{1}{2}J_1(\xi_h,\xi_h)^\frac{1}{2}\\
&\le C\left(\sum_{K\in\Th}h_K^{2t_K^1-2}\|u_1\|^2_{H^{s^k_K}(K)}\right)^\frac{1}{2}\|\psi_h\|_{h,1}.
\end{align*}
It then follows that
$$A_1(\xi_h,\psi_h)\le C_1\left(\sum_{K\in\Th}h_K^{2t_K^1-2}\|u_1\|^2_{H^{s^k_K}(K)}\right)^\frac{1}{2}\|\psi_h\|_{h,1}.$$
Applying the above estimate to~(\ref{16:39}), we obtain
$$\|\psi_h\|_{h,1}^2\le\kappa_1^{-1}C_1\left(\sum_{K\in\Th}h_K^{2t_K^1-2}\|u_1\|^2_{H^{s^k_K}(K)}\right)^\frac{1}{2}\|\psi_h\|_{h,1},$$
and so
$$\|\psi_h\|_{h,1}\le\kappa_1^{-1}C_1\left(\sum_{K\in\Th}h_K^{2t_K^1-2}\|u_1\|^2_{H^{s^k_K}(K)}\right)^\frac{1}{2}.$$
Applying the above estimate and~(\ref{otherxibnduseful}) to~(\ref{16:42::}) yields the desired result.$\quad\quad\square$
\begin{lemma}\label{errorest2lem}
Assume that $\X\subset\mathbb R^d$ is piecewise $C^{m+1}$ for some $m\in\mathbb N$, $m\ge3$, and that $(\Th)_{h>0}$ is a regular of order $m$ family of triangulations on $\overline{\X}$ satisfying Assumption~\ref{Meshconds}. Moreover, assume that $u_2\in H^2_0(\X)\cap H^{4}(\X)\cap H^{\mathbf{s}_2}(\X;\Th)$, where $\mathbf{s}_2 = (s_K^2)_{K\in\Th}$ and each $s_K^2\ge4$, is the true solution of~(\ref{2and4}) for $k=2$. Furthermore, let $\eta_F^2,\eta_F^3,\eta_F^4$ satisfy~(\ref{etacond2}) such that~(\ref{A2stab2}) holds for some $\kappa_2\in(0,1)$, and assume that $p\ge 3$. Then, the following estimate holds
\begin{equation}\label{regerrorest1}
\|u_{2,h}-u_2\|_{h,2}\le C_2\left(\sum_{K\in\Th}h_K^{2t_K^2-4}\|u_k\|_{H^{s_K^2}(K)}^2\right)^{1/2},
\end{equation}
where $t_K^2:=\min\{p+1,m+1,s_K^2\}$, $u_{2,h}\in V_{2,h}$ is the unique solution of~(\ref{A2def:14:47}), and the constant $C_2$ depends upon the shape-regularity constant of the mesh, $C_P$, $C_I$, $C_{\operatorname{Tr}}$, $C_\mathcal{T}$, $d$, $\mathbf{s}_K^2$, and $p$, but not upon $h_K$.
\end{lemma}
\emph{Proof:} Let us take $z_h\in V_{2,h}$, to be arbitrary, denoting $\xi_h:=u_2-z_h$ and $\psi_h:=z_h-u_{2,h}$, we see that
\begin{equation}\label{16:42::2}
\|u_2-u_{2,h}\|_{h,1}\le\|\xi_h\|_{h,1}+\|\psi_h\|_{h,2}.
\end{equation}
Furthermore, from~(\ref{xibnduseful}) for $k=2$, we see that
\begin{equation}\label{otherxibnduseful}
\|\xi_h\|_{h,1}\le C_2\left(\sum_{K\in\Th}h_K^{2t_K^2-4}\|u_2\|_{H^{s_K^2}(K)}^2\right)^{1/2},
\end{equation}
and so it remains to estimate $\|\psi_h\|_{h,1}$. Due to the stability estimate~(\ref{A2stab}), and the consistency of the scheme, since $\psi_h\in V_{2,h}$, we see that
\begin{equation}\label{16:392}
\begin{aligned}
\|\psi_h\|_{h,1}^2&\le\kappa_2^{-1}A_2(\psi_h,\psi_h)\\
&=\kappa_2^{-1}(A_2(z_h,\psi_h)-A_2(u_{2,h},\psi_h))\\
&=\kappa_2^{-1}(A_2(z_h,\psi_h)-\ell(\psi_h))\\
&=\kappa_2^{-1}(A_2(z_h,\psi_h)-A_2(u_2,\psi_h))\\
&=\kappa_2^{-1}A_2(\xi_h,\psi_h).
\end{aligned}
\end{equation}
Then,
applying the Cauchy--Schwarz inequality for vectors in $\mathbb R^N$, estimate~(\ref{Jbound15:33}), and the interpolation estimate~(\ref{opt:interp:2}), we obtain
\begin{equation}\label{Abound18:45}
\begin{aligned}
&A_2(\xi_h,\psi_h)  = \sum_{K\in\Th}\langle D^2 \xi_h,D^2 \psi_h\rangle_K+C( \xi_h,\psi_h)+C( \psi_h,\xi_h)+B_2( \xi_h,\psi_h)+B_2( \psi_h, \xi_h)+J_2(\xi_h,\psi_h)\\
&\le|\xi_h|_{H^2(\X;\Th)}|\psi_h|_{H^2(\X;\Th)}+C( \xi_h,\psi_h)+C( \psi_h,\xi_h)+B_2( \xi_h,\psi_h)+B_2( \psi_h, \xi_h)+J_2(\xi_h,\xi_h)^\frac{1}{2}J_2(\psi_h,\psi_h)^\frac{1}{2}\\
&\le C\max\{1,C_{*,2}^{-1/2}\}\left(\sum_{K\in\Th}h_K^{2t_K^2-4}\|u_2\|_{H^{s_K^2}(K)}^2\right)^{1/2}\|\xi_h\|_{h,1}\\
&~~~~~~~~~~~~~~~~~~~~~~~~~~~~~~~~~~~~~~~~~~~~~+C( \xi_h,\psi_h)+C( \psi_h,\xi_h)+B_2(\xi_h,\psi_h)+B_2( \psi_h, \xi_h).
\end{aligned}
\end{equation}
Furthermore, applying the Cauchy--Schwarz inequality for vectors in $\mathbb R^N$, the interpolation estimate~(\ref{opt:interp:2}), and estimates~(\ref{1.1in})--(\ref{18:43ddd}), we obtain
\begin{equation}\label{1:18:45}
\begin{aligned}
&C(\xi_h,\psi_h)\le\left(\sum_{F\in\Ei}\tilde{h}_F\|\Delta_\Ta\avg{\xi_h}\|_{L^2(F)}^2+\tilde{h}_F\mc_F^2\left\|\bavg{\frac{\partial\xi_h}{\partial n_F}}\right\|_{L^2(F)}^2\right)^\frac{1}{2}\left(\sum_{F\in\Ei}2\tilde{h}_F^{-1}\left\|\bjump{\frac{\partial\psi_h}{\partial n_F}}\right\|_{L^2(F)}^2\right)^\frac{1}{2}\\
&~~~~+\left(\sum_{F\in\Ei}\tilde{h}_F\left\|\nabla_\Ta\bavg{\frac{\partial\xi_h}{\partial n_F}}\right\|_{L^2(F)}^2+\tilde{h}_F\|\nabla_\Ta n_F^T\|_{L^\infty(F)}^2\left(\|\nabla_\Ta\avg{\xi_h}\|_{L^2(F)}^2+\left\|\bavg{\frac{\partial\xi_h}{\partial n_F}}\right\|_{L^2(F)}^2\right)\right)^\frac{1}{2}\\
&~~~~~~~~~~~~~~~~~~~~~~~~~~~~~~~~~~~~~~~~~~~~~~~~~~~~~~~~~~~~~~~~~~~~~~~~~~~~~~~~~~~\times\left(\sum_{F\in\Ei}3\tilde{h}_F^{-1}\|\jump{\nabla_\Ta\psi_h}\|_{L^2(F)}^2\right)^\frac{1}{2}\\
&~~~~~~~~~~~~~\le C\left(\sum_{K\in\Th}h_K^{2t_K^2-4}\|u_2\|_{H^{s_K^2}(K)}^2\right)^{1/2}J_2(\psi_h,\psi_h)^\frac{1}{2}\\
&~~~~~~~~~~~~~\le CC_{*,2}^{-1/2}\left(\sum_{K\in\Th}h_K^{2t_K^2-4}\|u_2\|_{H^{s_K^2}(K)}^2\right)^{1/2}\|\psi_h\|_{h,2}.
\end{aligned}
\end{equation}
We then apply the Cauchy--Schwarz inequality for vectors in $\mathbb R^N$, the interpolation estimate~(\ref{opt:interp:2}), the trace estimate~(\ref{1simptrace}), the inverse estimate~(\ref{local:invineq}), estimates~(\ref{1.1in})--(\ref{18:43ddd}), estimate~(\ref{Jbound15:33}), and the Poincar\'{e}--Friedrichs' inequality~(\ref{GPF:est:disc}), yielding
\begin{equation}\label{2:18:45}
\begin{aligned}
&C(\psi_h,\xi_h)\le\left(\sum_{F\in\Ei}\tilde{h}_F\|\Delta_\Ta\avg{\psi_h}\|_{L^2(F)}^2+\tilde{h}_F\mc_F^2\left\|\bavg{\frac{\partial\psi_h}{\partial n_F}}\right\|_{L^2(F)}^2\right)^\frac{1}{2}\left(\sum_{F\in\Ei}2\tilde{h}_F^{-1}\left\|\bjump{\frac{\partial\xi_h}{\partial n_F}}\right\|_{L^2(F)}^2\right)^\frac{1}{2}\\
&~~~~+\left(\sum_{F\in\Ei}\tilde{h}_F\left\|\nabla_\Ta\bavg{\frac{\partial\psi_h}{\partial n_F}}\right\|_{L^2(F)}^2+\tilde{h}_F\|\nabla_\Ta n_F^T\|_{L^\infty(F)}^2\left(\|\nabla_\Ta\avg{\psi_h}\|_{L^2(F)}^2+\left\|\bavg{\frac{\partial\psi_h}{\partial n_F}}\right\|_{L^2(F)}^2\right)\right)^\frac{1}{2}\\
&~~~~~~~~~~~~~~~~~~~~~~~~~~~~~~~~~~~~~~~~~~~~~~~~~~~~~~~~~~~~~~~~~~~~~~~~~~~~~~~~~~~\times\left(\sum_{F\in\Ei}3\tilde{h}_F^{-1}\|\jump{\nabla_\Ta\xi_h}\|_{L^2(F)}^2\right)^\frac{1}{2}\\
&\le C\left(\sum_{K\in\Th}|\psi_h|_{H^1(K)}^2+(\tilde{h}_F+1)|\psi_h|_{H^2(K)}^2+\tilde{h}_F|\psi_h|_{H^3(K)}^2\right)^{1/2}J_2(\xi_h,\xi_h)^{1/2}\\
&\le C\left(\sum_{K\in\Th}h_K^{2t_K^2-4}\|u_2\|_{H^{s_K^2}(K)}^2\right)^{\frac{1}{2}}\left(|v|_{H^2(\X;\Th)}^2+\sum_{F\in\Ei}\tilde{h}_F^{-1}\left\|\bjump{\frac{\partial\psi_h}{\partial n_F}}\right\|_{L^2(F)}^2+\sum_{F\in\Eib}\tilde{h}_F^{-1}\|\jump{\psi_h}\|_{L^2(F)}^2\right)^{\frac{1}{2}}\\
&\le C\max\left\{1,C_{*,2}^{-1/2}\right\}\left(\sum_{K\in\Th}h_K^{2t_K^2-4}\|u_2\|_{H^{s_K^2}(K)}^2\right)^{\frac{1}{2}}\|\psi_h\|_{h,2}.
\end{aligned}
\end{equation}
Now, by the Cauchy--Schwarz inequality for vectors in $\mathbb R^N$, and the interpolation estimate~(\ref{opt:interp:2}), we find that
\begin{equation}\label{3:18:45}
\begin{aligned}
B_2(\xi_h,\psi_h)&\le\left(\sum_{F\in\Eib}\tilde{h}_F^3\left\|\bavg{\frac{\partial(\Delta\xi_h)}{\partial n_F}}\right\|_{L^2(F)}^2\right)^{1/2}\left(\sum_{F\in\Eib}\tilde{h}_F^{-3}\|\jump{\psi_h}\|_{L^2(F)}^2\right)^{1/2}\\
&~~~~~~~~~~~~~+\left(\sum_{F\in\Eib}\tilde{h}_F\|\avg{\Delta\xi_h}\|_{L^2(F)}^2\right)^{1/2}\left(\sum_{F\in\Eib}\tilde{h}_F^{-1}\left\|\bjump{\frac{\partial\psi_h}{\partial n_F}}\right\|_{L^2(F)}^2\right)^{1/2}\\
&~~~~~~~~~~\le C\left(\sum_{K\in\Th}h_K^{2t_K^2-4}\|u_2\|_{H^{s_K^2}(K)}^2\right)^{1/2}J_2(\psi_h,\psi_h)^{1/2}\\
&~~~~~~~~~~\le CC_{*,2}^{-1/2}\left(\sum_{K\in\Th}h_K^{2t_K^2-4}\|u_2\|_{H^{s_K^2}(K)}^2\right)^{1/2}\|\psi_h\|_{h,2}.
\end{aligned}
\end{equation}
Finally, after applying the Cauchy--Schwarz inequality for vectors in $\mathbb R^N$, the interpolation estimate~(\ref{opt:interp:2}), the trace estimate~(\ref{1simptrace}), the inverse estimate~(\ref{local:invineq}), estimate~(\ref{Jbound15:33}), and the Poincar\'{e}--Friedrichs' inequality~(\ref{GPF:est:disc}), yielding
\begin{equation*}
\begin{aligned}
&B_2(\psi_h,\xi_h)\le\left(\sum_{F\in\Eib}\tilde{h}_F^3\left\|\bavg{\frac{\partial(\Delta\psi_h)}{\partial n_F}}\right\|_{L^2(F)}^2\right)^{1/2}\left(\sum_{F\in\Eib}\tilde{h}_F^{-3}\|\jump{\xi_h}\|_{L^2(F)}^2\right)^{1/2}\\
&~~~~~~~~~~~~~~~~~~~~~~~~~~+\left(\sum_{F\in\Eib}\tilde{h}_F\|\avg{\Delta\psi_h}\|_{L^2(F)}^2\right)^{1/2}\left(\sum_{F\in\Eib}\tilde{h}_F^{-1}\left\|\bjump{\frac{\partial\xi_h}{\partial n_F}}\right\|_{L^2(F)}^2\right)^{1/2}\\
\end{aligned}
\end{equation*}
\begin{equation}\label{4:18:45}
\begin{aligned}
&\le C\left(\sum_{K\in\Th}|\psi_h|_{H^2(K)}^2+h_K^2|\psi_h|_{H^3(K)}^2+h_K^4|\psi_h|_{H^4(K)}^2\right)^{1/2}J_2(\xi_h,\xi^2_h)^{1/2}\\
&\le C\left(\sum_{K\in\Th}h_K^{2t_K^2-4}\|u_2\|_{H^{s_K^2}(K)}^2\right)^{\frac{1}{2}}\left(|v|_{H^2(\X;\Th)}^2+\sum_{F\in\Ei}\tilde{h}_F^{-1}\left\|\bjump{\frac{\partial\psi_h}{\partial n_F}}\right\|_{L^2(F)}^2+\sum_{F\in\Eib}\tilde{h}_F^{-1}\|\jump{\psi_h}\|_{L^2(F)}^2\right)^{\frac{1}{2}}\\
&\le C\max\left\{1,C_{*,2}^{-1/2}\right\}\left(\sum_{K\in\Th}h_K^{2t_K^2-4}\|u_2\|_{H^{s_K^2}(K)}^2\right)^{\frac{1}{2}}\|\psi_h\|_{h,2}.
\end{aligned}
\end{equation}
Applying estimates~(\ref{1:18:45})--(\ref{4:18:45}) to~(\ref{Abound18:45}), and applying the resulting inequality to~(\ref{16:392}), we obtain
$$\|\psi_h\|_{h,2}^2\le\kappa_2^{-1}C_2\left(\sum_{K\in\Th}h_K^{2t_K^2-4}\|u_2\|_{H^{s_K^2}(K)}^2\right)^{\frac{1}{2}}\|\psi_h\|_{h,2},$$
and so
$$\|\psi_h\|_{h,2}\le\kappa_2^{-1}C_2\left(\sum_{K\in\Th}h_K^{2t_K^2-4}\|u_2\|_{H^{s_K^2}(K)}^2\right)^{\frac{1}{2}}.$$
Applying the above estimate, and estimate~(\ref{otherxibnduseful}), to~(\ref{16:42::2}) yields the desired result.$
\quad\quad\square$
\section{Numerical results}\label{sec:6}
We now briefly discuss the implementation of the numerical methods in Firedrake~\cite{Rathgeber2016,Luporini2015}
\subsection{Implementation}\label{sec:5:imp}
\emph{Software and code:} The experiments in this Chapter have been implemented in the most recent version of the Firedrake software~\cite{Rathgeber2016,Luporini2015} (as of 3rd July 2018), which interfaces directly with PETSc~\cite{petsc-user-ref,petsc-efficient} running through a Python interface~\cite{Dalcin2011,Chaco95}. There are two Firedrake scripts, Curved-Dirichlet-DGFEM.py (applicable to~(\ref{A1def:14:47})), and DGFEM-curved-biharmonic.py (applicable to~(\ref{A2def:14:47})) used to generate the experiments of this section is available in the Github repository:\\ https://github.com/ekawecki/Firedrake\_Poisson\_Biharmonic.

\emph{Linear systems and condition numbers:}
The evaluation of $A_k(u_h,v_h)$ for $u_h,v_h\in\dg$, where the  bilinear forms $A_k$ are defined by~(\ref{A1def:14:47}) for $k=1$, and~(\ref{A2def:14:47}) for $k=2$ involves the integration of products of first and second order partial derivatives, respectively. This typically leads to the matrix $\mathbf{A}_k$, describing the linear system given by~(\ref{A1def:14:47}) for $k=1$, and~(\ref{A2def:14:47}) for $k=2$, to have a Euclidean norm condition number of order $h^{-2k}$. This can pose difficulties when applying iterative methods to solve the linear system (particularly for $k=2$) and thus to ensure that we solve the linear system with sufficiently high accuracy as the mesh size $h$ decreases, we apply the Iterative refinement algorithm, i.e., Algorithm 1.1 of~\cite{carson2017new}. We implement the Iterative refinement algorithm by using the following choices in the Firedrake ``solve" function.
\begin{lstlisting}[language=Python]
        solve(A == L, U,
              solver_parameters = {
              "snes_type": "newtonls",
              "ksp_type": "preonly",
              "pc_type": "lu",
              "snes_monitor": False,
              "snes_rtol": 1e-16,
              "snes_atol": 1e-25})
\end{lstlisting}
\emph{Two-dimensional curved boundary approximation:} When implementing curved finite elements, we use a piecewise quadratic polynomial mapping to obtain a higher order approximation of the domain boundary. This is implemented in Firedrake by first using Gmsh~\cite{geuzaine2008gmsh} (version 3.0.1) to generate an affine triangulation $\X_h$ that approximates $\X$. We then define the \emph{continuous} Lagrange finite element space $\mathbb V:=\{v\in C(\overline{\X_h};\mathbb R^2):v\in\mathbb P^2(K;\mathbb R^2)\,\forall K\in\X_h\}$.
Then, we take $\psi_i:\omega_i\to\mathbb R^2$, $\omega_i\subset\mathbb R$, $i=1,\ldots,n$, to be the collection of charts that locally describe $\partial\X$, and denote $\{x_j\}_{j=1}^N$ to be the degrees of freedom of $\mathbb V$. We partition the collection of degrees of freedom by defining $J_{\operatorname{ext}} = \{j\in\{1,\ldots,N\}:x_j\in\partial\X_h\}$, and $J_{\operatorname{int}} = \{1,\ldots,N\}\setminus J_{\operatorname{ext}}$, and so $\{x_j\}_{j=1}^N = \{x_j\}_{j\in J_{\operatorname{int}}}\cup\{x_j\}_{j\in J_{\operatorname{ext}}}$. We then
 define the the function $T\in\mathbb V$ by
\begin{equation}\label{T:def}
\left\{
\begin{aligned}
T(x_j) &= x_j,\quad j\in J_{\operatorname{int}},\\
T(x_j) &=\psi_i(x_j),\quad j\in J_{\operatorname{ext}},\quad i\in\{1,\ldots,n\}\mbox{ such that }x_j\in\omega_i.
\end{aligned}
\right.
\end{equation}
Finally, we define our computational finite element space $\dgcomp:= \{v\in L^2(\X):v\circ T^{-1}\in\mathbb P^p(\hat K)\}$. 
This procedure is implemented in Firedrake, in the code snippet below, utilising the Firedrake ``Mesh" function. In this case $\X$ is the unit disk, and so there is only one chart, $\psi:=x/|x|$. Furthermore, when we refine the mesh in our experiments, the meshes at each refinement level are not related to one another. That is, there is no hierarchical mesh structure, i.e., at each refinement level, we ``remesh". A collection of the meshes used for the computations of this paper can be found in the folder ``Meshes" in the Github repository: https://github.com/ekawecki/Firedrake\_Poisson\_Biharmonic.
\begin{lstlisting}[language=Python]
        #Affine mesh of the unit disk, generated in Gmsh
        mesh = Mesh("quasiunifrefdisk.msh")
        # Implementing quadratic domain approximation
        V = FunctionSpace(mesh, "CG", 2)
        # Defining a function that identifies the curved portion of the boundary
        bdry_indicator = Function(V)
        bc = DirichletBC(V, Constant(1.0), 1)
        bc.apply(bdry_indicator)
        # Defining the continuous, piecewise quadratic vector-valued finite element space
        VV = VectorFunctionSpace(mesh, "CG", 2)
        T = Function(VV)
        T.interpolate(SpatialCoordinate(mesh))
        # Defining the function T given by (6.1)
        T.interpolate(conditional(abs(1-bdry_indicator) < 1e-5, T/sqrt(inner(T,T)), T))
        # Defining the curved mesh
        mesh = Mesh(T)
        # Defining the space V_{h,p}^{comp}
        FES = FunctionSpace(mesh,"DG",deg)
\end{lstlisting}
\subsection{Experiment 1}\label{exp1}
In this experiment, we consider the Poisson problem~(\ref{2and4}) (for $k=1$), with $\X=\{x\in\mathbb R^2:|x|<1\}$, and right-hand side function $f$ chosen so that the true solution 
$$u(x_1,x_2) = \frac{1}{4}\sin(\pi(x_1^2+x_2^2)).$$
We took the penalty parameter $\eta^1_F=10p^4$ (obtained experimentally), where $p$ is the polynomial degree of the space $\dgcomp$. For each polynomial degree $p=1,2,3$, we successively refined the mesh quasiuniformly. We observe the predicted optimal convergence rate $\|u-u_h\|_{h,1}=\mathcal{O}(h^{p})$, as well as the optimal rate $\|u-u_h\|_{L^2(\X)}=\mathcal{O}(h^{p+1})$, with the true values and EOCs in brackets provided in Tables~\ref{exp1h,1errs} and~\ref{exp1L2errs}, respectively.
\begin{table}[H]
\begin{center}
\begin{tabular}{|c|cc|cc|cc|}\hline
\cline{1-7} 
\multicolumn{1}{|c|}{Mesh size} & \multicolumn{2}{|c|}{$p=1$} &  \multicolumn{2}{|c|}{$p=2$} & \multicolumn{2}{|c|}{$p=3$}\\\hline
\cline{1-7}
$ 0.4981 $&$ 2.83 $&$ $&$ 4.34\times 10^{-1} $&$ $&$ 2.83\times 10^{-1} $&$ $\\
$ 0.2828 $&$ 1.51 $&$( 1.11 )$&$ 9.99\times 10^{-2} $&$( 2.59 )$&$ 5.10\times 10^{-2} $&$( 3.03 )$\\
$ 0.1627 $&$ 7.06\times 10^{-1} $&$( 1.38 )$&$ 5.58\times 10^{-2} $&$( 1.06 )$&$ 1.00\times 10^{-2} $&$( 2.94 )$\\
$ 0.0973 $&$ 3.66\times 10^{-1} $&$( 1.28 )$&$ 2.01\times 10^{-2} $&$( 1.98 )$&$ 1.58\times 10^{-3} $&$( 3.59 )$\\
$ 0.0508 $&$ 1.68\times 10^{-1} $&$( 1.20 )$&$ 5.37\times 10^{-3} $&$( 2.03 )$&$ 2.08\times 10^{-4} $&$( 3.12 )$\\
$ 0.0269 $&$ 8.51\times 10^{-2} $&$( 1.06 )$&$ 1.45\times 10^{-3} $&$( 2.06 )$&$ 2.69\times 10^{-5} $&$( 3.21 )$\\\hline
\end{tabular}
\caption{Error values in the $\|\cdot\|_{h,1}$-norm and EOCs for Experiment 1.}\label{exp1h,1errs}
\end{center}
\end{table}
\begin{table}[H]
\begin{center}
\begin{tabular}{|c|cc|cc|cc|}\hline
\cline{1-7} 
\multicolumn{1}{|c|}{Mesh size} & \multicolumn{2}{|c|}{$p=1$} &  \multicolumn{2}{|c|}{$p=2$} & \multicolumn{2}{|c|}{$p=3$}\\\hline
\cline{1-7}
$ 0.4981 $&$ 1.24\times 10^{-1} $&$ $&$ 3.57\times 10^{-2} $&$ $&$ 5.18\times 10^{-3} $&$ $\\
$ 0.2828 $&$ 5.17\times 10^{-2} $&$( 1.55 )$&$ 3.28\times 10^{-3} $&$( 4.21 )$&$ 1.44\times 10^{-3} $&$( 2.26 )$\\
$ 0.1627 $&$ 1.70\times 10^{-2} $&$( 2.01 )$&$ 1.25\times 10^{-3} $&$( 1.75 )$&$ 1.71\times 10^{-4} $&$( 3.85 )$\\
$ 0.0973 $&$ 5.74\times 10^{-3} $&$( 2.12 )$&$ 2.49\times 10^{-4} $&$( 3.14 )$&$ 1.55\times 10^{-5} $&$( 4.68 )$\\
$ 0.0508 $&$ 1.38\times 10^{-3} $&$( 2.19 )$&$ 3.38\times 10^{-5} $&$( 3.07 )$&$ 1.11\times 10^{-6} $&$( 4.05 )$\\
$ 0.0269 $&$ 3.70\times 10^{-4} $&$( 2.07 )$&$ 4.69\times 10^{-6} $&$( 3.10 )$&$ 7.40\times 10^{-8} $&$( 4.25 )$\\\hline
\end{tabular}
\caption{Error values in the $\|\cdot\|_{L^2(\X)}$-norm and EOCs for Experiment 1.}\label{exp1L2errs}
\end{center}
\end{table}
\subsection{Experiment 2}\label{exp2}
In this experiment, we consider the biharmonic problem~(\ref{2and4}) (for $k=2$), with $\X=\{x\in\mathbb R^2:|x|<1\}$, and right-hand side function $f$ chosen so that the true solution 
$$u(x_1,x_2) = \sin^2(\pi(x_1^2+x_2^2)).$$
We took the penalty parameter $\eta^2_F=c_pp^6$, and $\eta^3_F,\eta^4_F=c_pp^4$  where $p$ is the polynomial degree of the space $\dgcomp$ (the order of these parameters with respect to $p$ were guided by the choice of penalty parameters in Section 6 of~\cite{suli2007hp}), and $c_p=0.1$ for $p=2$, and $c_p=10$ for $p=3,4$. For each polynomial degree $p=2,3,4$, we successively refined the mesh quasiuniformly. We observe the optimal convergence rate $\|u-u_h\|_{h,2}=\mathcal{O}(h^{p-1})$ for $p=2,3,4$, confirming the estimate of Lemma~\ref{errorest2lem}. We also observe the optimal rate $|u-u_h|_{H^1(\X;\Th)}=\mathcal{O}(h^{p})$, for $p=2,3,4$. We provide the error values, and EOCs (in brackets) in the $\|\cdot\|_{h,2}$-norm and $|\cdot|_{H^1(\X;\Th)}$-seminorm in Tables~\ref{exp2h,2errs} and~\ref{exp2H1errs}, respectively.
\begin{table}[H]
\begin{center}
\begin{tabular}{|c|cc|cc|cc|}\hline
\cline{1-7} 
\multicolumn{1}{|c|}{Mesh size} & \multicolumn{2}{|c|}{$p=2$} &  \multicolumn{2}{|c|}{$p=3$} & \multicolumn{2}{|c|}{$p=4$}\\\hline
\cline{1-7}
$ 0.4981 $&$ 1.05\times 10^2 $&$ $&$ 5.70\times 10^1 $&$ $&$ 4.60\times 10^1 $&$ $\\
$ 0.2828 $&$ 7.33\times 10^1 $&$( 0.63 )$&$ 3.30\times 10^1 $&$( 0.97 )$&$ 1.12\times 10^1 $&$( 2.49 )$\\
$ 0.1627 $&$ 4.74\times 10^1 $&$( 0.79 )$&$ 1.36\times 10^1 $&$( 1.60 )$&$ 2.64 $&$( 2.62 )$\\
$ 0.0973 $&$ 2.60\times 10^1 $&$( 1.17 )$&$ 5.05 $&$( 1.93 )$&$ 6.32\times 10^{-1} $&$( 2.78 )$\\
$ 0.0508 $&$ 1.28\times 10^1 $&$( 1.09 )$&$ 1.24 $&$( 2.16 )$&$ 9.31\times 10^{-2} $&$( 2.95 )$\\
$ 0.0269 $&$ 6.61 $&$( 1.04 )$&$ 4.07\times 10^{-1} $&$( 1.75 )$&$ 1.39\times 10^{-2} $&$( 2.98 )$\\\hline
\end{tabular}
\caption{Error values in the $\|\cdot\|_{h,2}$-norm and EOCs for Experiment 2.}\label{exp2h,2errs}
\end{center}
\end{table}
\begin{table}[H]
\begin{center}
\begin{tabular}{|c|cc|cc|cc|}\hline
\cline{1-7} 
\multicolumn{1}{|c|}{Mesh size} & \multicolumn{2}{|c|}{$p=2$} &  \multicolumn{2}{|c|}{$p=3$} & \multicolumn{2}{|c|}{$p=4$}\\\hline
\cline{1-7}
$ 0.4981 $&$ 3.89 $&$ $&$ 5.08 $&$ $&$ 3.41 $&$ $\\
$ 0.2828 $&$ 2.78 $&$( 0.60 )$&$ 1.67 $&$( 1.97 )$&$ 3.35\times 10^{-1} $&$( 4.10 )$\\
$ 0.1627 $&$ 1.49 $&$( 1.12 )$&$ 3.45\times 10^{-1} $&$( 2.85 )$&$ 5.38\times 10^{-2} $&$( 3.31 )$\\
$ 0.0973 $&$ 6.07\times 10^{-1} $&$( 1.75 )$&$ 8.06\times 10^{-2} $&$( 2.83 )$&$ 7.58\times 10^{-3} $&$( 3.81 )$\\
$ 0.0508 $&$ 1.67\times 10^{-1} $&$( 1.98 )$&$ 9.89\times 10^{-3} $&$( 3.23 )$&$ 5.45\times 10^{-4} $&$( 4.05 )$\\
$ 0.0269 $&$ 4.71\times 10^{-2} $&$( 1.99 )$&$ 1.64\times 10^{-3} $&$( 2.82 )$&$ 4.18\times 10^{-5} $&$( 4.03 )$\\\hline
\end{tabular}
\caption{Error values in the $|\cdot|_{H^1(\X;\Th)}$-seminorm and EOCs for Experiment 2.}\label{exp2H1errs}
\end{center}
\end{table}
\section{Conclusion}
In the setting of curved finite elements, we have successfully reviewed several key estimates from theory of finite elements on polytopal domains, such as trace estimates, inverse estimates, discrete Poincar\'e--Friedrichs' inequalities, and optimal interpolation estimates in noninteger Hilbert-Sobolev norms, that are well known in the case of polytopal domains. Furthermore, we have proven curvature bounds for curved simplices, and utilised all of these estimates by providing stability, and a priori error analysis, of the IPDG method for the Poisson problem, orginally introduced in~\cite{douglas1976interior}, and for a variant of the $h$-version of the $hp$-DGFEM for the biharmonic problem introduced in~\cite{suli2007hp}.

In Section~\ref{sec:6}, we have provided numerical experiments for both the Poisson and biharmonic problem, where the domain is taken to be the unit disk. We implement a polynomial approximation of the domain, validating the a priori error estimates of Lemmas~\ref{errorest1lem} and~\ref{errorest1lem}. The estimates proven as part of this paper should serve useful for future applications to second- and fourth-order (as well as higher order) elliptic problems on curved domains, in particular, nondivergence form second-order elliptic equations.
\bibliographystyle{plain}
\bibliography{toskaweckiBIB}

Ellya L. Kawecki, Department of Mathematics and Center for Computation and Technology, Louisiana State University, Baton Rouge, LA 70803, USA

\emph{E-mail address}: \tt{ekawecki@cct.lsu.edu}
\end{document}